\documentclass[12pt]{article}
\usepackage[margin=1.0in]{geometry}
\usepackage[utf8]{inputenc}
\usepackage{tikz}
\usepackage{mathtools}
\usepackage{mathrsfs,amssymb,amsfonts,amsthm,amsmath,amsbsy,bbm}
\usepackage[english,french]{babel}
\usepackage[T1]{fontenc}
\usepackage{hyperref}
\usepackage{float}
\usepackage{mdframed}
\usepackage{amssymb}
\usepackage{verbatim}
\usepackage{dsfont}
\usepackage{eso-pic}
\usepackage{authblk,hyperref}
\usepackage{stmaryrd}
\usepackage{soul}
\usepackage{titlesec} 
\usepackage{listings}
\usepackage{lipsum,enumerate}
\usepackage[export]{adjustbox}
\usepackage[
backend=biber,
style=numeric,maxnames=99
]{biblatex}
\newtheorem{theorem}{Theorem}[section]
\newtheorem{proposition}[theorem]{Proposition} 
\newtheorem{corollary}[theorem]{Corollary}
\newtheorem{lemma}[theorem]{Lemma} 
 
\newtheorem{assumption}[theorem]{Assumption}

\newtheorem{definition1}[theorem]{Definition}

\newtheorem{remark1}[theorem]{Remark}
\newenvironment{remark}{\begin{remark1}\rm}{\hfill$\square$\end{remark1}}
\newtheorem{fact}[theorem]{Fact}

\newcommand{\somme}[3]{\sum_{#1}^{#2} #3}
\newcommand{\indi}[1]{\mathds{1}_{#1}}
\newcommand{\Pf}{\mathbb{P}}

\newcommand{\ent}[1]{\lfloor #1 \rfloor}

\DeclarePairedDelimiter\floor{\lfloor}{\rfloor}

\newcommand{\Lt}{\widetilde{L}}
\newcommand{\Ht}{\widetilde{H}}

\newcommand{\eps}{\varepsilon}

\newcommand{\R}{\mathbb{R}}
\newcommand{\rooot}{\rho}
\newcommand{\Chn}{\C_{H=n}}
\newcommand{\Chgn}{\C_{H \geq n}}
\newcommand{\Chgen}{\C_{H \geq (1-\eps) n}}

\newcommand{\rk}[1]{\left(#1\right)}

\newcommand{\C}{\mathcal{C}}
\newcommand{\F}{\mathcal{F}}
\newcommand{\Cge}[1]{\mathcal{C}_{\ge #1}}
\newcommand{\T}{\mathcal{T}}
\newcommand{\bT}{\mathbf{T}}
\newcommand{\Pbt}{\mathbb{P}_{\mathbf{T}}}
\newcommand{\Pb}{\mathbb{P}_{\alpha}}
\newcommand{\bPb}{\mathbf{P}_{\alpha}}
\newcommand{\Ta}{\mathcal{T}_{\alpha}}

\newcommand{\cadlag}{c\`adl\`ag }
\newcommand{\CC}{C_{\alpha}}
\newcommand{\K}{K_{\alpha}}
\newcommand{\bW} {\mathbf{W}} 

\newcommand{\lr}{\leftrightarrow}
\newcommand{\oo}{\mathbf{o}}

\newcommand{\lrpc}{\overset{p_c}{\longleftrightarrow }}

\newcommand{\pr}[1]{\mathbb{P}_{\alpha}\!\left(#1\right)}

\newcommand{\prcond}[3]{\mathbb{P}_{#3}\!\left(#1\;\middle\vert\;#2\right)}

\newcommand{\prb}[1]{\mathbf{P}_{\alpha}\!\left(#1\right)}
\newcommand{\Eb}[1]{\mathbf{E}_{\alpha}\!\left[#1\right]}

\newcommand{\econdb}[2]{\mathbf{E}_{\alpha}\!\left[#1\;\middle\vert\;#2\right]}

\newcommand{\Varb}[1]{\mathrm{\mathbf{Var}}_{\alpha}\!\left(#1\right)}
\newcommand{\PfT}{\mathbb{P}_{\mathbf{T}}}

\newcommand{\dGHP}[1]{d_{GHP}\!\left(#1\right)}

\newcommand{\Ito}{It\^o }
\newcommand{\Height}{\textsf{Height}}

\newcommand{\Pt}[1]{\mathbb{P}_{\mathbf{T}}
\!\left(#1\right)}
\newcommand{\Et}[1]{\mathbb{E}_{\mathbf{T}}\!\left[#1\right]}

\newcommand{\Ptcond}[3]{\mathbb{P}_{#3}\!\left(#1\;\middle\vert\;#2\right)}

\title{\textbf{Quenched scaling limit of critical percolation clusters on Galton-Watson trees}}

\author{Eleanor Archer\thanks{Universit\'e Paris-Dauphine, archer@ceremade.dauphine.fr} \ and Tanguy Lions 
\thanks{ENS Lyon, tanguy.lions@ens-lyon.fr}}

\date{}

\begin{document}

\maketitle
\selectlanguage{english}

\begin{abstract}
We consider quenched critical percolation on a supercritical Galton--Watson tree with either finite variance or $\alpha$-stable offspring tails for some $\alpha \in (1,2)$. We show that the Gromov-Hausdorff-Prokhorov (GHP) scaling limit of a quenched critical percolation cluster on this tree is the corresponding $\alpha$-stable tree, as is the case in the annealed setting. As a corollary we obtain that a simple random walk on the cluster also rescales to Brownian motion on the stable tree. Along the way, we also obtain quenched asymptotics for the tail of the cluster size, which completes earlier results obtained in Michelen (2019) and Archer-Vogel (2024).
\end{abstract}

\begin{center}
   \begin{figure}[H]
    \centering
\includegraphics[scale = 0.6]{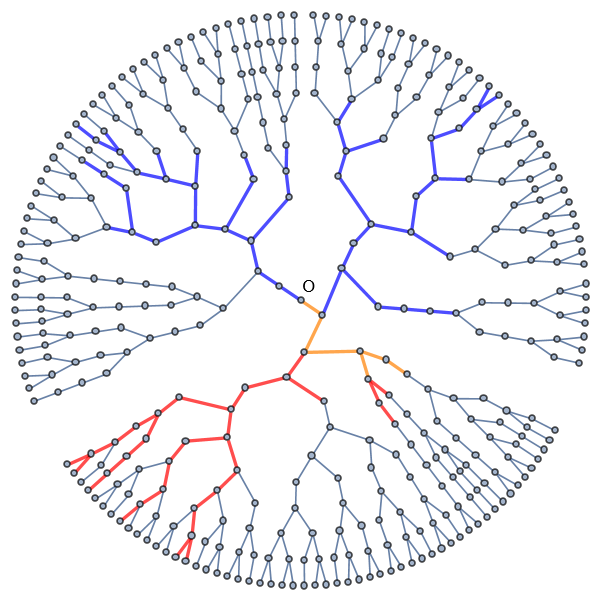}
    \label{}
    \caption{A supercritical Galton-Watson tree cut at level $11$. The blue and red parts are two independent critical percolation clusters containing $\rooot$ on the tree. The orange part is the intersection of the two percolation clusters.}
\end{figure} 
\end{center}

%
%\textcolor{red}{Liste des corrections effectuées.\\
%referee 1 (2 page report):
%\begin{enumerate}
%	\item Done
%	\item Seems specific to the AOAP template
%	\item Seems specific to the AOAP template
%	\item Seems specific to the AOAP template
%	\item Done
%	\item Done
%	\item Done
%	\item Done 
%	\item Done 
%	\item Done
%\end{enumerate}
%referee 2 (5 page report):
%\begin{itemize}
%	\item[$\bullet$] I treated the section 4 called "questions". In particular I made Section~\ref{sctn:contour limit} more concise by combining  4 Propositions in just one. The referee seemed to feel this part was not important to spend time on so I made it simple. Otherwise, I modified Proposition~\ref{technical_proposition}, it should be correct now.
%\end{itemize}
%\begin{enumerate}
%	\item Done
%	\item Done
%	\item Done, I removed the $b \in (0,1)$. We do not use it.
%	\item  Done
%	\item Done (removed first citation)
%	\item Done
%	\item Done
%	\item Done
%	\item  Done
%\end{enumerate}
%}

\tableofcontents

\section{Introduction}

Let $\mathbf{T}$ be a supercritical Galton-Watson tree, with root $\rooot$. We suppose that its offspring distribution has mean $\mu > 1$, is supported on $\{1,2, \ldots \}$ and that it is either in the domain of attraction of a stable law with parameter $\alpha \in (1,2)$, or has finite variance. In the latter case we set $\alpha=2$. Given $\alpha$, we let $\mathbf{P}_{\alpha}$ denote the law of $\bT$. It was shown by Lyons \cite[Theorem $6.2$ and Proposition $6.4$]{Lyons_critical_percolation} that $\bPb$ almost-surely, the critical (Bernoulli) percolation threshold on $\mathbf{T}$ is $\frac{1}{\mu}$. The aim of this paper is to obtain a quenched Gromov-Hausdorff-Prokhorov (GHP) scaling limit of a critical percolation cluster on $\bT$, conditioned on its exact height or size. In addition, we obtain quenched convergence of a simple random walk on the cluster to Brownian motion on the limiting fractal tree.

Conditionally on $\mathbf{T}$, we let $\mathbb{P}_{\mathbf{T}}$ denote the law of critical Bernoulli percolation on $\mathbf{T}$. The \emph{annealed} law $\Pb$ is defined by $\mathbb{P}_{\alpha} = \bPb \circ \mathbb{P}_{\mathbf{T}}$. Under $\mathbb{P}_{\alpha}$, the root cluster (henceforth denoted by $\C$) has the law of a \textit{critical} Galton-Watson tree with offspring distribution in the domain of attraction of an $\alpha$-stable law or with finite variance. We root it at $\rooot$. Consequently, it is known that under $\mathbb{P}_{\alpha}$, the GHP scaling limit of the critical cluster conditioned to be large is an $\alpha$-stable Lévy tree or the continuum random tree (CRT). In particular, if $\C_{=n}$ denotes the cluster conditioned to have size equal to $n$, $d_n$ the intrinsic graph metric on $\C_{=n}$, and $\nu_n$ the counting measure on its vertices, then there exists a random rooted metric-measure space $(\mathcal{T}^{= 1}_{\alpha},d_{\mathcal{\mathcal{T}_{\alpha}}},\nu_{\alpha},\rho_{\alpha})$ (whose law depends only on $\alpha$) and an explicit constant $\gamma \in (0, \infty)$ (depending on the offspring law of $\bT$) such that
   \begin{equation}\label{eqn:annealed GHP}
  (\C_{=n}, \gamma n^{-\left(1-\frac{1}{\alpha}\right)}d_n, n^{-1}\nu_n, \rooot) \overset{(d)}{\to} (\mathcal{T}^{= 1}_{\alpha},d_{\mathcal{\mathcal{T}_{\alpha}}},\nu_{\alpha},\rho_{\alpha})
  \end{equation}
as $n \to \infty$. See \cite{AldousCRTiii, legall2005random, duquesne2003contourlimit}. The space $(\mathcal{T}^{= 1}_{\alpha},d_{\mathcal{\mathcal{T}_{\alpha}}},\nu_{\alpha},\rho_{\alpha})$ is known as the \textit{$\alpha$-stable tree} (conditioned to have size exactly $1$, as indicated by the superscript) in the case $\alpha \in (1,2)$, and more commonly as the \textit{Brownian tree} or the \textit{CRT} in the case $\alpha=2$. The main result of this paper is Theorem~\ref{theorem_exactconditioning}. It states that the same scaling limit result is true under $\Pbt$, for $\bPb$-almost every $\mathbf{T}$. This result is proved in Section~\ref{section_conditioning_on_the_size}.

We will work under the following assumption throughout the paper.

\begin{assumption}\label{assumption}
    Assume that the offspring distribution of $\bT$ is supported on $\{1, 2, \ldots\}$, its mean is given by $\mu>1$, and that one of the following conditions holds.
    \begin{enumerate}[(a)]
        \item (Finite variance case.) The offspring distribution of $\bT$ has finite variance $\sigma^2$. In this case set $\alpha=2$.%, $\CC = \frac{2\Eb{\abs{\bT_1}}^2}{\Eb{\abs{\bT_1}\rk{\abs{\bT_1}-1}}}$ and ${\Eb{\abs{\bT_1}\rk{\abs{\bT_1}-1}}}$.
        \item (Stable case.) The offspring distribution of $\bT$ has infinite variance with stable (power-law) tails, meaning that there exist $c \in (0, \infty)$ and $\alpha \in (1,2)$ such that $\prb{|\bT_1| \geq x} \sim c x^{-\alpha}$ as $x \to \infty$ (and in this case we will use the subscript $\alpha$ to denote the dependence on $\alpha$).% In this case, we set
%        \begin{equationj
%            \CC=c^{-\beta}_1\mu^{\alpha\beta}\Gamma(1-\alpha)^{-\beta}\beta^\beta\qquad\text{and}\qquad \phi (\theta)=1-\CC^{-1}\theta (1+(\CC^{-1}\theta)^{\alpha-1})^{-\beta}\, .
%        \end{equation}
    \end{enumerate}
\end{assumption}

\begin{remark}
We exclude the case $\alpha=2$ from case (b) above for ease of reading as this necessitates adding various logarithmic scaling corrections to all of the results. However, in this setting the annealed scaling limit is again the CRT (this follows for example from \cite[Theorem 4.17]{kallenberg1997foundations} and \cite[Theorem 2.1.1, Theorem 2.3.1, Theorem 2.3.2]{duquesne2002random}) and our proof should apply in the quenched setting too. Similarly, the proof should still work in exactly the same way on allowing for slowly-varying corrections to the offspring tails (and carrying these through the proofs). In addition, we anticipate that the assumption that $\bT$ has no leaves could be removed using the Harris decomposition of supercritical Galton--Watson trees, which decomposes such a tree conditioned to survive into a supercritical core that has no leaves (to which our main theorem applies) to which finite Galton--Watson trees are attached (see \cite[Proposition 5.28]{lyons2017probability}). However transferring the result would require some work and we decided not to pursue this in the present paper.
\end{remark}

Before stating the results, we recall an important result about $\bT$ itself, which plays a role in the first theorem. Let $\bT_n$ be the set of vertices at generation $n$ in $\bT$. It is well-known that there is a random variable $\bW$ such that
\begin{equation*}
    \bW_n := \frac{|{\bT_n}|}{\mu^n}\to \bW\, ,
\end{equation*}
as $n\to\infty$ almost surely and in $\mathrm{L}^p$ if $\Et{|{\bT_1}|^p}<\infty$, see \cite[Theorems 0 and 5]{bingham1974asymptotic} (in particular this holds whenever $p<\alpha$).

We start with a result on the quenched convergence of the law of the total size of $\C$. This fills a gap from \cite{archer2023quenched} and is an ingredient in the proof of our main GHP convergence result.

In the annealed setting, it is well-known (see for example \cite[Lemma A.3(i)]{curien2015percolation}) that there exists a constant $\K \in (0, \infty)$ such that, as $n \to \infty$,
\begin{equation}\label{eqn:annealed size prob convergence}
    n^{\frac{1}{\alpha}}\pr{|\C|>n}\to \K\, .
\end{equation}
Our first result shows that this is also true in the quenched setting, up to a small dependence on the tree $\mathbf{T}$.

\begin{theorem}\label{thm:ConvergenceOfSize}
For $\bPb$-almost every $\bT$, we have that
\begin{equation*}
    n^{\frac{1}{\alpha}}\Pt{|\C|>n}\to \K \bW
\end{equation*}
as $n \to \infty$.
\end{theorem}

To state our main results, we first introduce some notation. We let $\C_{=n}$ (respectively $\C_{\ge n}$) denote the root percolation cluster conditioned on having total size exactly $n$ (respectively at least $n$), $d_n$ the intrinsic graph metric on $\C_{=n}$ (respectively $\C_{\ge n}$), $\nu_n$ the counting measure on its vertices (likewise for $\C_{\ge n}$), and $\rho_n$ its root (which coincides with the root $\rooot$ of $\bT$). We similarly let $\Chn$ (respectively $\Chgn$) denote the root percolation cluster conditioned on having total height exactly $n$ (respectively height at least $n$), and similarly for $d_n, \nu_n, \rho_n$ as above. We recall that $\mathcal{T}^{= 1}_{\alpha}$ denotes the $\alpha$-stable tree conditioned to have size $1$. We similarly denote by $\mathcal{T}_{\alpha}^{H=1}$ the $\alpha$-stable tree conditioned to have height equal to $1$, $\mathcal{T}^{\geq 1}_{\alpha}$ the $\alpha$-stable tree conditioned to have size at least $1$, and $\mathcal{T}^{H\geq 1}_{\alpha}$ the $\alpha$-stable tree conditioned to have {height} at least one. We state the first main theorem that gives the scaling limits for the clusters conditioned to have size or height at least $n$. The pointed Gromov-Hausdorff-Prokhorov topology will be defined in Section \ref{sctn:GHP topology}.

\begin{theorem}\label{main_theorem}
%Calling $\mathcal{C}_n$ the root percolation cluster conditioned on having height at least $n$ we have :
%\begin{itemize}
%    \item[$\bullet$]$\bPb$-almost surely,  $\displaystyle (\mathcal{C}_n,\frac{d_n}{\textcolor{red}{\gamma} n},\frac{ \nu_n}{\textcolor{red}{\gamma^{'}}n^2},\rho_n) \underset{n \to +\infty}{\overset{(d)}{\longrightarrow}} (\mathcal{T},d_{\mathcal{\mathcal{T}}},\nu,\rho)$.
%    \item[$\bullet$] $\forall \alpha \in (1,2), \mathbf{P}_{\alpha}$-almost surely, $\displaystyle (\mathcal{C}_n,\frac{d_n}{\textcolor{red}{\gamma} n},\frac{\nu_n}{\textcolor{red}{\gamma^{'}}n^{\frac{\alpha}{\alpha - 1}}},\rho_n) \underset{n \to +\infty}{\overset{(d)}{\longrightarrow}} (\mathcal{T}_{\alpha},d_{\mathcal{\mathcal{T}_{\alpha}}},\nu_{\alpha},\rho_{\alpha})$. 
%\end{itemize}
%Here, the topology for the distribution convergence is the pointed Gromov-Hausdorff-Prokhorov topology. The trees $(\mathcal{T},d_{\mathcal{\mathcal{T}}},\nu,\rho)$ and $(\mathcal{T}_{\alpha},d_{\mathcal{\mathcal{T}_{\alpha}}},\nu_{\alpha},\rho_{\alpha})$ are respectively the Continuum Random Tree and the stable tree with index $\alpha$.
Take $\gamma$ as in \eqref{eqn:annealed GHP}. Then, for $\bPb$-almost every $\bT$, the following convergence holds in law under $\mathbb{P}_{\bT}$:
\begin{align}
(\mathcal{C}_{\ge n}, \gamma n^{-\left(1-\frac{1}{\alpha}\right)}d_n, n^{-1}\nu_n,\rho_n) &\underset{n \to +\infty}{\overset{(d)}{\longrightarrow}} (\mathcal{T}^{\geq 1}_{\alpha},d_{\mathcal{\mathcal{T}_{\alpha}}},\nu_{\alpha},\rho_{\alpha}) \label{eqn:first geq conv size}\\
(\Chgn, n^{-1}d_n,{(\gamma n)^{-\frac{\alpha}{\alpha - 1}}} \nu_n,\rho_n) &\underset{n \to +\infty}{\overset{(d)}{\longrightarrow}} (\mathcal{T}^{H\geq 1}_{\alpha},d_{\mathcal{\mathcal{T}_{\alpha}}},\nu_{\alpha},\rho_{\alpha})\label{eqn:first geq conv height}
\end{align}
with respect to the pointed Gromov-Hausdorff-Prokhorov topology.
\end{theorem}

We state the second main theorem that gives the scaling limits for the clusters conditioned to have size or height exactly $n$.

\begin{theorem}\label{theorem_exactconditioning}
Take $\gamma$ as in \eqref{eqn:annealed GHP}. Then, for $\bPb$-almost every $\bT$, the following convergence holds in law under $\mathbb{P}_{\bT}$:
\begin{align*}
&(\mathcal{C}_{=n}, \gamma n^{-\left(1-\frac{1}{\alpha}\right)}d_n, n^{-1}\nu_n,\rho_n) \underset{n \to +\infty}{\overset{(d)}{\longrightarrow}} (\mathcal{T}_{\alpha}^{=1},d_{\mathcal{\mathcal{T}_{\alpha}}},\nu_{\alpha},\rho_{\alpha}),\\
 &(\Chn, n^{-1}d_n,{(\gamma n)^{-\frac{\alpha}{\alpha - 1}}} \nu_n,\rho_n) \underset{n \to +\infty}{\overset{(d)}{\longrightarrow}} (\mathcal{T}_{\alpha}^{H=1},d_{\mathcal{\mathcal{T}_{\alpha}}},\nu_{\alpha},\rho_{\alpha}).
\end{align*}
with respect to the pointed Gromov-Hausdorff-Prokhorov topology.
\end{theorem}

Although Theorem \ref{main_theorem} may in fact be deduced from Theorems \ref{thm:ConvergenceOfSize} and \ref{main_theorem}, we state the two theorems in this order as this reflects our proof strategy. In particular Theorem \ref{main_theorem} is a crucial ingredient in the proof of Theorem \ref{theorem_exactconditioning}.

\begin{remark1}
The constant $\gamma$ appearing in \eqref{eqn:annealed GHP} and Theorem \ref{main_theorem} can be computed explicitly as a function of the offspring law of $\bT$. In particular, in the finite variance case $\gamma = \frac{\widetilde{\sigma}}{2}$ where $\widetilde{\sigma}^2= \frac{\sigma^2}{\mu^2}+1-\mu^{-1}$ is the variance of the annealed law. In the stable case, $\gamma = \left(|\Gamma (1-\alpha)|c\right)^{1/\alpha}$ (see \cite[Chapter 8.3]{bingham1989regular} for background).
\end{remark1}

We add to this with convergence of the simple random walk on $\C_{=n}$ to Brownian motion on $\Ta^{= 1}$. In what follows we let $(X^{(n)}_m)_{m \geq 0}$ denote a discrete time simple random walk on $\C_{=n}$, with quenched law $\mathcal{P}_n$, and we let $(B_t)_{t \geq 0}$ denote Brownian motion on $\Ta^{= 1}$, with quenched law $\mathcal{P}$. These objects will be introduced in Section \ref{sctn:SRW conv}. Of course, the corollary could equally be stated on ${\C}_{\geq n}$, $\Chn$ or $\Chgn$ (with appropriately updated scaling exponents).

\begin{corollary}\label{cor:SRW convergence}
$\bPb$-almost surely, there exists a probability space $(\Omega_{\bT}, \F_{\bT}, \Pbt)$ on which the convergence of Theorem \ref{main_theorem} holds almost surely. Then, on this space:
\[
\mathcal{P}_n \left(\left(X^{(n)}_{\lfloor \gamma^{-1} n^{\frac{2\alpha-1}{\alpha}} t \rfloor}\right)_{t \geq 0} \in \cdot \right) \to \mathcal{P} \left(\left(B_t \right)_{t \geq 0} \in \cdot \right)
\]
weakly as probability measures on the space of \cadlag paths equipped with the Skorokhod-$J_1$ topology.
\end{corollary}

\begin{remark1}
Corollary \ref{cor:SRW convergence} can also be stated formally as a joint convergence with that of Theorem \ref{main_theorem} using a topology constructed in \cite{Khezeli}. We will give more details about this in Section \ref{sctn:SRW conv}.
\end{remark1}

\paragraph{Physical motivation and applications.} Percolation has well-known applications in the study of a variety of physical systems, such as fluid flow through porous media, spread of disease and magnetism. The study of percolation at criticality is especially delicate and often gives rise to anomalous behaviour. Moreover, the study of the associated simple random walks (also known as the problem of the ant in the labyrinth) has been an important research area since the seminal work of Kesten \cite{kesten1986incipient} who first established subdiffusive behaviour of random walks on critical percolation clusters. The convergence of these random walks to Brownian motion on the continuum random tree is expected to be a universal phenomenon in appropriate high-dimensional settings and establishing this in some generality is an active area of research. See for example \cite{ben2019scaling, arous2019scaling, heydenreich2014random} for some results in this direction.

Random trees serve as a good proxy for many physical models (for example, the Alexander-Orbach conjecture, proved in high dimensions by Kozma and Nachmias \cite{kozma2009alexander} in 2009, states that the key random walk exponents for various percolation models agree with those for a random walk on a critical tree). Consequently, the study of statistical physics models or particle systems on random trees can give insight into the behaviour of the same models on more complicated physical structures. Since most real-world systems are intrinsically random, it is moreover a natural question to study these models in the quenched setting and understand how and why the behaviour may deviate from the system's typical (annealed) behaviour.

\paragraph{Sketch of proof.}
The proof of Theorem \ref{thm:ConvergenceOfSize} follows a similar strategy to that used to establish an analogous result for the height of $\C$ in \cite{archer2023quenched}: in particular, we choose $m$ of order $\log n$ such that, with high probability on the event $\{|\C|\geq n\}$, there is a single vertex at generation $m$ in $\bT$ that connects to the root and in addition has a large percolation cluster in the subtree emanating from it. The result of the theorem is then essentially obtained by averaging over the choice of this vertex.

To prove Theorem \ref{main_theorem}, rather than working directly with $\C$, we work with its so-called \textit{height function} (see Section \ref{sctn:coding of forests} for a definition). In particular, we show that the height function $X$ coding a sequence of i.i.d. samples of $\C$ (under $\Pbt$) converges to the height function coding an i.i.d. forest of stable trees. For technical reasons it is also helpful to keep track of the local time of $X$ at $0$, which we will denote by $\Lambda$. From this it is fairly classical to deduce the result of Theorem \ref{main_theorem} (by restricting to the first tree in the forest with size or height exceeding $n$).

The strategy to obtain convergence of height functions is a second moment argument on the quantity $\mathbb{E}_{\mathbf{T}}\bigg[F\left(\left(n^{-(1-\frac{1}{\alpha})}X_{\floor{nt}}, n^{-\frac{1}{\alpha}}\Lambda_{\floor{nt}}\right)_{0 \le t \le T}\right)\bigg]$, where $F$ is a bounded Lipschitz function $F : \C[0,T]^2 \longrightarrow \R$, which will then be sufficient to apply a Borel-Cantelli argument. To bound the variance of this quantity, the intuition is roughly as follows: conditionally on $\mathbf{T}$, we take two independent copies of $\C$, which we denote by $\mathcal{C}^{1}$ and $\mathcal{C}^{2}$. Note that $\mathcal{C}^{1}\cap \mathcal{C}^{2}$ (formed from the unconditioned clusters) is a subcritical root percolation cluster, and hence the clusters $\mathcal{C}^{1}$ and $\mathcal{C}^{2}$ visit disjoint parts of $\mathbf{T}$ as soon as they are not too close to the origin. The same logic applies on sampling further copies of $\C$ and this essentially breaks the dependence between different trees in the forest coded by $X$. To conclude one just has to argue that the parts of the clusters near the origin are small and average out over large scales.

Finally, let us describe the strategy to prove Theorem~\ref{theorem_exactconditioning}. We only describe the strategy to prove the convergence of $\C_{=n}$ since the strategy for $\Chn$ is very similar. The starting point is the following remark: the tree $\C_{=n}$ can be sampled by first sampling $\C_{\ge (1-\varepsilon)n}$ and then additionally conditioning this tree to have size exactly $n$. Under this conditioning we approximate $\C_{=n}$ by $\C_{n,\eps}$, the subtree obtained from $\C_{\ge (1-\varepsilon)n}$ by trimming it at the first generation where the total mass up to that point exceeds $(1-\varepsilon)n$, which is proved to lead to a decent approximation of $\C_{=n}$ (see point (c) of Proposition~\ref{estimateA_i}). Moreover, the behaviour of $\C_{n,\varepsilon}$ is captured by Theorem~\ref{main_theorem} (see Theorem~\ref{main_theorem with generation size}). All that there remains to show is that, in some sense, conditionally on $\C_{n, \eps}$, we have $\Pt{\#\C = n \mid \#\C \ge (1-\varepsilon)n} \sim \Pb(\#\C = n \mid \#\C \ge (1-\varepsilon)n)$, i.e. proving that the quenched probability behaves as the annealed one asymptotically. To do this we import a result of \cite{archer2023quenched} which allows us to additionally control the final generation size of $\C_{n, \eps}$, jointly with the convergence of Theorem \ref{main_theorem}. This final generation size essentially determines the conditional probabilities and the asymptotic can then be obtained using a similar second moment strategy to the one detailed above for the proof of Theorem~\ref{main_theorem}. This strategy is outlined in more detail in Section \ref{sctn:strategy exact cond}.

Corollary \ref{cor:SRW convergence} is a direct consequence of Theorem \ref{main_theorem} and a general result of Croydon \cite{croydon2018scaling}. \\

We mention that the strategy to upgrade Theorem \ref{main_theorem} to Theorem \ref{theorem_exactconditioning} is somewhat inspired by the proof of \cite{kortchemski2013simple}, which was in turn inspired by ideas of \cite[Sections 6 and 7]{le2010ito}, to make the same upgrade in the annealed setting. Here the authors instead consider a depth-first exploration of the tree and cut it at the first moment at least $(1-\eps)n$ vertices have been explored. They then similarly show that this reduced tree captures the behaviour of the entire tree conditioned to have size $n$ (as $\eps \downarrow 0$); moreover the conditional probability of having size exactly $n$ can be written in terms of certain depth-first coding functions and converges to a limiting density expressed in terms of the limiting coding functions. We expect that this could also be achieved using our approach of cutting at heights: in this case the limiting density would be written in terms of a certain local time at the relevant level in the limiting fractal tree. Moreover, our general upgrade strategy of cutting at heights to compare quenched and annealed probabilities is fairly robust and should be more generally applicable to sequences of random trees for which GHP convergence, as well as convergence of the sequence of generation sizes, are both known to hold under conditioning the size or height to be at least $n$. 

We also remark that we expect similar results to be true for critical percolation on hyperbolic random planar maps, for which an annealed GHP scaling limit was obtained in \cite{archercroydon2023scaling}. However, the quenched analysis is more delicate due to the loss of tree structure, which means in particular that clusters can (in theory) be disjoint in some annulus and then merge again outside the annulus.

\paragraph{Organisation.} The paper is organised as follows. In Section \ref{sctn:background} we give background on quenched critical  percolation on Galton--Watson trees as well as scaling limits of the latter. In Section \ref{sctn:size conv} we prove Theorem \ref{thm:ConvergenceOfSize}. In Section \ref{sctn:contour limit} we establish the main ingredient to prove Theorem \ref{main_theorem}, namely the {convergence of the height functions} via a second moment estimate, postponing the proof of one technical proposition to the Appendix. In Section \ref{section_conditioning_on_the_size} we explain in detail how to deduce the first statement of Theorem \ref{theorem_exactconditioning} from Theorem \ref{main_theorem}, via a careful analysis of the law of $\# \C$, conditionally on $\C_{n, \eps}$. Again we postpone some technical details to the Appendix. In Section \ref{section_conditioning_on_the_height} we also give an outline of how the same approach works when conditioning on the height. Finally in Section \ref{sctn:SRW conv} we explain the framework needed to deduce Corollary \ref{cor:SRW convergence}.

\paragraph{Acknowledgements.} The research of EA was partially funded by ANR ProGraM (reference ANR-19-CE40-0025). We are also grateful to ENS Lyon and ENS Paris-Saclay for funding the research of Tanguy Lions.

\section{Background}\label{sctn:background}

As partly mentioned in the introduction, the random tree $\bT$ is defined on the probability space $(\mathbf{\Omega}, \mathcal{F}, \bPb)$. For $h \geq 0$, we let $\mathcal{F}_h \subset \mathcal{F}$ denote the sigma-algebra generated by the first $h$ generations of $\bT$. Given $\bT$, the critical cluster $\C$ is defined on the space $(\Omega_{\bT}, \mathcal{G}_{\bT}, \mathbb{P}_{\bT}$), where $\mathcal{G}_{\bT}$ denotes the canonical sigma algebra on subsets of edges of $\bT$ (generated by cylinder sets).

\subsection{Critical percolation on Galton--Watson trees}\label{sctn:background quenched perco}

Here we give a brief outline of known results about critical percolation on Galton--Watson trees under Assumption \ref{assumption}. 

We first recall an important result about $\bT$ itself, which plays a role in certain quenched results. Let $\bT_n$ be the set vertices in generation $n$ in $\bT$. It is well-known that there is a random variable $\bW$ such that
\begin{equation}\label{eqn:W convergence}
    \bW_n := \frac{|{\bT_n}|}{\mu^n}\to \bW\, ,
\end{equation}
as $n\to\infty$ almost surely and in $\mathrm{L}^p$ if $\Et{|{\bT_1}|^p}<\infty$; see \cite[Theorems 0 and 5]{bingham1974asymptotic} (in particular this holds whenever $p<\alpha$).

In the annealed setting, the critical cluster of $\bT$ is just a critical Galton--Watson tree with offspring distribution \textsf{Binomial}($Z, 1/\mu$) where $Z$ follows the offspring distribution of $\bT$, and $\mu$ is its mean. As a consequence, the large-scale behaviour of the cluster is essentially completely understood: asymptotic tails for its height and total size, and various scaling limits (see the left hand side of Table \ref{table:annealed quenched} for a full list).

We mention two annealed results to which we will make specific reference. The first concerns the asymptotics for the tails of the cluster size and has already been stated in \eqref{eqn:annealed size prob convergence}. 
% First, letting $\C$ be the critical cluster of the origin, it is well-known (see for example \cite[Lemma A.3(i)]{curien2015percolation}) that there exists a constant $\K \in (0, \infty)$ such that, as $n \to \infty$,
% \begin{equation}\label{eqn:annealed size prob convergence}
%     n^{\frac{1}{\alpha}}\pr{|\C|>n}\to \K\, .
% \end{equation}
% \red{already in (2)}
Similarly, letting \textsf{Height}($\C$) denote the height of $\C$, that is, $\sup \{n \geq 0: \bT_n \cap \C \neq \emptyset\}$, it was shown in \cite{Slack} that there exists a constant $\CC \in (0, \infty)$ such that, as $n \to \infty$,
\begin{equation}\label{eqn:slack prob convergence}
    n^{\frac{1}{\alpha-1}}\pr{\text{\textsf{Height}}(\C) \geq n}\to \CC\, .
\end{equation}
% In addition, we also know from \cite[Lemma 2.5]{archer2020random} that the function $x \mapsto \pr{\text{\textsf{Height}}(\C) =x}$ is non-increasing in $x$ (we can take $H_Y=\text{\textsf{Height}}(\C) $ there). Hence for all sufficiently large $n$ we have
% \begin{align*}
% \pr{\text{\textsf{Height}}(\C)  = n} \leq 2n^{-1}\pr{\text{\textsf{Height}}(\C)  \in \left(\frac{n}{2},n\right]} \leq 2^{\frac{\alpha}{\alpha-1}} \CC n^{-\frac{\alpha}{\alpha-1}}.
% \end{align*}
% \red{can remove}
For any $n \ge 0$, we denote by $Y_n$ the size of the generation of $\mathcal{C}$ at height $n$. In the quenched setting, the first relevant result for us is that of Lyons \cite{Lyons_critical_percolation} who showed that $p_c (\bT) = 1/\mu$ almost surely. The problem was later studied by Michelen \cite{michelen}, who showed, under some moment conditions, quenched convergence of connection probabilities as well as quenched convergence of the rescaled law of $Y_n$, conditioned on survival (in particular he established the well-known \textit{Yaglom limit}). The moment conditions were relaxed in \cite{archer2023quenched} and the scaling limit result was extended to prove convergence of the entire sequence of generation sizes to a continuous state branching process.

These results are listed in Table \ref{table:annealed quenched}. The notable gap in the previous results is GHP convergence of the cluster, which is addressed in the present work. Along the way we also obtained quenched convergence of the cluster size tails.

\begin{center}
\begin{table}[h]
\centering
	\begin{tabular}{c|c}
\textbf{Annealed} & \textbf{Quenched} \\
\hline
$p_c = 1/\mu$ & $p_c (T) = 1/\mu$ a.s.   \\

$\pr{\oo \lrpc T_n} \sim \CC n^{-\frac{1}{\alpha-1}} $ & $\Pt{\oo \lrpc T_n} \sim {\bW} \cdot \CC n^{-\frac{1}{\alpha-1}} $ \\

$\pr{|C| \geq n} \sim \K n^{-1/\alpha} $ & $\Pt{|C| \geq n} \sim {\bW} \cdot \K n^{-1/\alpha} $ ($\ast$) \\

Given $Y_n>0$: $n^{-\frac{1}{\alpha-1}}Y_n \overset{(d)}{\to} Y$ &  {$n^{-\frac{1}{\alpha-1}}Y^T_n \overset{(d)}{\to} Y$} a.s. \\

$\left(n^{-\frac{1}{\alpha-1}}Y_{n(1+t)}\right)_{t \geq 0} \overset{(d)}{\to} (Y_t)_{t\geq 0}$ & {$\left(n^{-\frac{1}{\alpha-1}}Y^T_{n(1+t)}\right)_{t \geq 0} \overset{(d)}{\to} (Y_t)_{t\geq 0}$} a.s. \\

Given $Y_{\infty}>0$: $n^{-\frac{1}{\alpha-1}}Y_n \overset{(d)}{\to} Y^*$ &  {$n^{-\frac{1}{\alpha-1}}Y^T_n \overset{(d)}{\to} Y^*$} \\

$(C_{\geq n},\gamma n^{-(1-1/\alpha)} d_n, \frac{1}{n}\mu_n) \underset{\text{GHP}}{\overset{(d)}{\to}} \mathcal{T}^{\geq 1}_{\alpha}$ & {$(C^T_{\geq n}, n^{-(1-1/\alpha)} d_n, \frac{1}{n}\mu_n) \underset{\text{GHP}}{\overset{(d)}{\to}} \mathcal{T}^{\geq 1}_{\alpha}$} a.s.  ($\ast$) \\

$(C_{= n}, \gamma n^{-(1-1/\alpha)} d_n, \frac{1}{n}\mu_n) \underset{\text{GHP}}{\overset{(d)}{\to}} \mathcal{T}^{= 1}_{\alpha}$ & {$(C^T_{= n}, n^{-(1-1/\alpha)} d_n, \frac{1}{n}\mu_n) \underset{\text{GHP}}{\overset{(d)}{\to}} \mathcal{T}^{= 1}_{\alpha}$} a.s.  ($\ast$) \\
\end{tabular}
\caption{Summary of annealed vs. quenched results. In this paper we prove the results labelled ($\ast$). (We also prove analogues of the last two statements conditioned on the height.)}\label{table:annealed quenched}
\end{table}
\end{center}

\subsection{Gromov-Hausdorff-type topologies}\label{sctn:GHP topology}

We now introduce the pointed Gromov-Hausdorff-Prokhorov (GHP) topology under which Theorem \ref{main_theorem} is stated. To this end, let $\mathbb{K}_c$ denote the set of quadruples $(K,d,\mu,\rho)$ such that $(K,d)$ is a compact metric space, $\mu$ is a locally-finite Borel measure on $K$, and $\rho$ is a distinguished point of $K$. Suppose that $(K,d,\mu,\rho)$ and $(K',d',\mu',\rho')$ are elements of $\mathbb{K}_c$. Given a metric space $(M, d_M)$, and isometric embeddings $\phi$ and $\phi'$ of $(K,d)$ and $(K',d')$ respectively into $(M, d_M)$, we define $d_{M}\big((K,d,\mu,\rho, \phi), (K',d',\mu',\rho', \phi')\big)$ to be equal to
\begin{align*}
d_M^H(\phi (K), \phi' (K')) + &d_M^P(\mu \circ \phi^{-1}, \mu' \circ {\phi'}^{-1} ) + d_M(\phi (\rho), \phi' (\rho')).
\end{align*}
Here $d_M^H$ denotes the Hausdorff distance between two sets in $(M, d_M)$, and $d_M^P$ denotes the Prokhorov distance between two measures (see for example \cite[Chapter 1]{BillsleyConv} for a definition). The pointed Gromov-Hausdorff-Prokhorov distance between $(K,d,\mu,\rho)$ and $(K',d',\mu',\rho')$ is then given by
\begin{equation}\label{eqn:GHP def}
\dGHP{(K,d,\mu,\rho), (K',d',\mu',\rho')} = \inf_{\phi, \phi', M} d_{M}\big((K,d,\mu,\rho, \phi), (K',d',\mu',\rho', \phi')\big)
\end{equation}
where the infimum is taken over all isometric embeddings $\phi, \phi'$ of $(X,d)$ and $(X',d')$ into a common metric space $(M, d_M)$. This defines a metric on the space of equivalence classes of $\mathbb{K}_c$ (see \cite[Theorem 2.5]{AbDelHoschNoteGromov}), where we say that two spaces $(K,d,\mu,\rho)$ and $(K',d',\mu',\rho')$ are equivalent if there is a measure and root-preserving isometry between them. Moreover, $\mathbb{K}_c$ is a Polish space with respect to the topology induced by $d_{GHP}$ (again, see \cite[Theorem 2.5]{AbDelHoschNoteGromov}).

Later, in order to pass from the convergence of the full space in Theorem \ref{main_theorem} to the balls of a certain radius in Section \ref{section_conditioning_on_the_size} we will use the following deterministic result. It can be proved straightforwardly using the definition of the GHP topology; we leave the proof to the reader. (The constants on the right hand side are not necessarily optimal.)

\begin{lemma}\label{lem:GHP conv on balls}
Suppose that $\dGHP{(\widetilde{X}, \widetilde{d}, \widetilde{\mu}, \widetilde{\rho}), (X, d, \mu, \rho) } \leq \eps$. Then, for all $R>0$,
\begin{align*}
\dGHP{(B_R(\widetilde{X}), \widetilde{d}|_R, \widetilde{\mu}|_R, \widetilde{\rho}), (X|_R, d|_R, \mu|_R, \rho) } \leq 2\eps \vee  \mu\left(B_{R+3\eps}({X}) \setminus B_{R-\eps}({X})\right).
\end{align*}
\end{lemma}

We also mention two extensions of this topology, that allow us to keep track of some extra information. Firstly, it will be useful in Section \ref{section_conditioning_on_the_size} to keep track of a certain generation size; for this we will work in the space $\mathbb{K}_c \times \R_{\geq 0}$, endowed with the metric $D((M,x), (M',x')) = \max (\dGHP{M,M'}, |x-x'|)$ (this induces the product topology).

Secondly, we will need the following extension, introduced in \cite{Khezeli}, which incorporates càdlàg paths on $K$. To this end, we let $\widetilde{\mathbb{K}}_c$ denote the set of quintuplets $(K,d,\mu,\rho,X)$, where $(K,d,\mu,\rho)\in \mathbb{K}_c$ and $X$ is a càdlàg path from $[0,\infty)$ to $K$. Similarly to above, given a metric space $(M, d_M)$, and isometric embeddings $\phi, \phi'$ of $(K,d)$ and $(K',d')$ respectively into $(M, d_M)$, we define $\widetilde{d}_{M}\big((K,d,\mu,\rho, X,\phi), (K',d',\mu',\rho',X', \phi')\big)$ to be equal to
\begin{align*}
d_M^H(\phi (K), \phi' (K')) + &d_M^P(\mu \circ \phi^{-1}, \mu' \circ {\phi'}^{-1} ) + d_M(\phi (\rho), \phi' (\rho'))+d_M^{J_1}(\phi (X), \phi' (X')),
\end{align*}
where $d_M^{J_1}$ is the metrisation of the Skorokhod $J_1$-topology for c\`{a}dl\`{a}g paths on $M$ described in \cite[Example 3.44]{Khezeli}. We then set
\[
d_{\widetilde{\mathbb{K}}_c}\left({(K,d,\mu,\rho,X), (K',d',\mu',\rho',X')}\right) = \inf_{\phi, \phi', M} \widetilde{d}_{M}\big((K,d,\mu,\rho,X, \phi), (K',d',\mu',\rho',X', \phi')\big),
\]
where again the infimum is taken over all isometric embeddings $\phi$ and $\phi'$ of $(X,d)$ and $(X',d')$ into a common metric space $(M, d_M)$, which yields a distance on $\widetilde{\mathbb{K}}_c$(see \cite{Khezeli}). Moreover, $d_{\widetilde{\mathbb{K}}_c}$ defines a metric on the space of equivalence classes of $\widetilde{\mathbb{K}}_c$, where we say that two spaces $(K,d,\mu,\rho,X)$ and $(K',d',\mu',\rho',X')$ are equivalent if there is a measure, root and c\`{a}dl\`{a}g path preserving isometry between them. As above, $\widetilde{\mathbb{K}}_c$ is a Polish space with respect to the topology induced by $d_{\widetilde{\mathbb{K}}_c}$. 

\subsection{Convergence of random forests}

In this section we discuss convergence of random forests formed from sequences of i.i.d. Galton--Watson trees, along with their coding functions. The results are all taken from \cite{duquesne2002random}.

We will restrict the following discussion to \textit{plane trees}, meaning that there is a distinguished root vertex, and that the set of offspring of each vertex comes pre-equipped with a left-right ordering. In pictures, the root will be drawn at the base of the tree.

We will fix a parameter $\alpha \in (1,2]$ and assume that the corresponding offspring law $\Pb$ has expectation equal to $1$, and satisfies the following.

\begin{assumption}\label{assn:offspring}
$\Pb$ is aperiodic and critical and one of the following two conditions hold:
\begin{enumerate}[(I)]
\item $\alpha=2$ and $\Pb$ has finite variance $\sigma^2$.
\item $\alpha \in (1,2)$ and if $X \sim \Pb$, then there exists a constant $c\in (0, \infty)$  such that
\[
\pr{X > x} \sim cx^{-\alpha}
\]
as $x \to \infty$.
\end{enumerate}
\end{assumption}

In the latter case we say that $X$ is in the domain of attraction of an $\alpha$-stable law. One can also incorporate slowly-varying functions into the tails; we have omitted this for ease of reading. The results we mention can also be easily adapted to hold in the periodic case, but for the results of this paper the aperiodic case suffices (even if our supercritical tree $\bT$ has a periodic offspring law, the offspring law for the critical cluster will still be aperiodic).

\subsubsection{Coding of forests}\label{section_coding_of_forests}

We let $(T_i)_{i=1}^{\infty}$ denote a sequence of finite plane trees (the canonical case to have in mind is a sequence of i.i.d. sequence of Galton--Watson trees, each with offspring distribution $\Pb$, supported on $\{0, 1, \ldots, \}$). We now explain how to code this forest by a walk. We start with the case of a single tree for simplicity.

Suppose that $\T$ is a plane tree with $|\T|=n+1$. 
%Here we define the \textit{contour function} which encodes its structure: it is illustrated in Figure \ref{fig:contourheightfns} and defined as follows. We consider the motion of a particle that starts on the left of the root $\emptyset$ at time zero, and then continuously traverses the boundary of $\mathcal{T}$ at speed one, in the clockwise direction, until returning to the left side of the root. $X^{\T}(t)$ is equal to the height (i.e. distance from the root) of the particle at time $t$. Since each edge is traversed twice in this process, the contour function is defined in this way up until time equal to $2n$.
We first define the {lexicographical ordering} of vertices: for this we consider the motion of a particle that starts on the left of the root $\emptyset$ at time zero, and then continuously traverses the boundary of $\mathcal{T}$ at speed one, in the clockwise direction, until returning to the left side of the root. The \textit{lexicographical ordering} of the vertices corresponds to the order in which the vertices are first visited by this process (with no repeats). The \textit{height function} $H^{\T}$ is then defined by considering the vertices $u_0, u_1, \ldots, u_n$ in this lexicographical order, and then setting $H^{\T}_i$ to be equal to the generation of vertex $u_i$. The height function is defined precisely up until time $n$. Note that $H^{\T}_0=0$ but the function is otherwise strictly positive.

%, and the height function will be strictly positive.

%\begin{figure}[h]
%\includegraphics[width=13cm, height=5cm]{heightcontour1.jpg}
%\centering
%\caption{Example of contour function and height function for the %given tree}\label{fig:contourheightfns}
%\end{figure}
\begin{figure}[H]
\begin{center}
\includegraphics[scale=0.4]{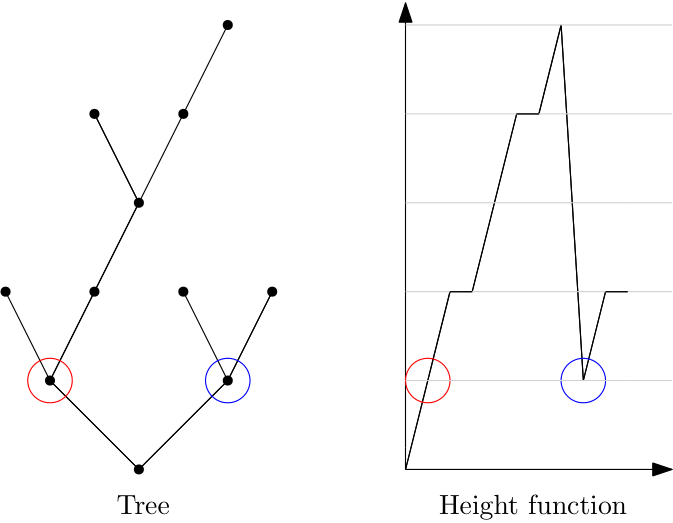}
\caption{Coding functions for the given tree. We have marked two vertices on the tree along with the points corresponding to these vertices in the excursions.}\label{fig:contourheightfns}
\end{center}
\end{figure}

The \textit{height} of $\T$ is equal to $\sup_{0 \leq j \leq n} H(j)$, and gives the maximal tree distance between any vertex and the root.

We can similarly encode \textit{forests} (that is, sequences of plane trees) by concatenating the corresponding height functions: formally, for $j \geq 0$, we set
\[
H(j) = H^{T_k}(j-\sum_{i=1}^k |T_i|) \qquad  \text{ if } \qquad  \sum_{i=1}^k |T_i| \leq j < \sum_{i=1}^{k+1} |T_i|.
\]
We then define $\tau_0=0$, and
\begin{align}\label{eqn:tau lambda def}
\tau_k = \inf\{j > \tau_{k-1}: H(j) = 0\}= \sum_{i=1}^k |T_i|, \hspace{1cm} \Lambda_j = \inf \{k : \tau_k > j\}.
\end{align}
Observe that the tree $T_i$ is coded by the interval $[\tau_{i-1}, \tau_i)$, and $\Lambda_j=i$ means that $j \in [\tau_{i-1}, \tau_i)$. The function $(\Lambda_j)_{j \geq 0}$ is known as the \textit{local time at zero} of $H$.
 %In the case of the contour and height functions, there is no ambiguity here since both functions start and end at $0$, but for the Lukasiewicz path we make the important convention that we concatenate using the \textit{increments}. In particular, if $T_i$ denotes the $i^{th}$ tree in the forest, then the part of the Lukasiewicz path coding $T_i$ is an excursion from $-(i-1)$ to $-i$. With this convention, it is straightforward to verify that \eqref{eqn:height Luk def} still holds (here the height of a vertex in $T_i$ is the distance to the root of $T_i$).
 
The importance of this concatenated height function is that it is actually in bijection with the forest. This provides an appealing way to construct scaling limits of random forests: we take scaling limits of the concatenated height functions, and then invert the bijection to construct a candidate for the limiting forest. One then just has to verify that the various operations are appropriately continuous. 

Since our eventual aim is to look at the scaling limit of a single tree conditioned on being large, sampled as the first tree in the forest satisfying that condition, it will also be important to keep track of the local time function $(\Lambda_j)_{j \geq 0}$. In particular it will be important that the first excursion of $H$ of length at least $n$ converges to the first excursion of length at least $1$ in its scaling limit. This may fail if the long discrete excursion comes arbitrarily close to zero in its interior, thus creating extra visits to zero in the scaling limit; this problem can be ruled out by additionally requiring that the local times converge.

\subsubsection{Scaling limits and continuum trees}\label{sctn:cont trees}

Duquesne and Le Gall \cite[Corollary 2.5.1]{duquesne2002random}, building on results of Le Gall and Le Jan \cite{LeGLeJ}, showed that the approach outlined above can be made precise and more specifically that one can construct a continuum height function $\Ht$, with associated local time at zero denoted by $\Lt$ such that the desired convergence of coding functions holds.

% (note that the stated corollary comprises the joint convergence of the first three functions below; we can add the local time due to \eqref{eqn:local time is min}, \eqref{eqn:local time is min cont}, monotonicity and since the process $(I_t)_{t \geq 0}$ is continuous).

\begin{proposition}\label{prop:annealed coding fn conv whole forest}
Under Assumption \ref{assn:offspring}, there exist (random) functions $\Ht$, $\Lt$ from $[0, \infty) \to \R$ and constants $c_2, c_3 \in (0, \infty)$ such that
%   \begin{align*}
%      ({c_1n^{-\frac{1}{\alpha}}}W_{\ent{nt}},{c_2n^{-\left(1-\frac{1}{\alpha}\right)}}H_{\ent{nt}},{c_2n^{-\left(1-\frac{1}{\alpha}\right)}}X_{\ent{2nt}},{c_3n^{-\frac{1}{\alpha}}}\Lambda_{\ent{nt}})_{t \ge 0} \underset{n \to +\infty}{\overset{(d)}{\rightarrow}} (\Wt_t, \Ht_t, \Ht_t, \Lt_t)_{t \ge 0},
%  \end{align*}
   \begin{align*}
      ({c_2n^{-\left(1-\frac{1}{\alpha}\right)}}H_{\ent{nt}}, {c_3n^{-\frac{1}{\alpha}}}\Lambda_{\ent{nt}})_{t \ge 0} \underset{n \to +\infty}{\overset{(d)}{\rightarrow}} (\Ht_t,  \Lt_t)_{t \ge 0},
  \end{align*}
jointly with respect to the uniform topology. Moreover the functions $\Ht$ and $\Lt$ are almost surely continuous and $\Lt_t$ corresponds to the local time of $\Ht$ at zero up until time $t$.

Under Assumption \ref{assn:offspring}(I), the function $\Ht$ is a reflected Brownian motion, $c_2 = \frac{2}{\sigma}$ and $c_3=\sigma$.
\end{proposition}

We refer to \cite{duquesne2002random} for the formal definitions of these processes.

Proposition \ref{prop:annealed coding fn conv whole forest} also suggests a natural way to define continuum trees. Notably, in the discrete setting, it is straightforward to verify that the graph metric $d_{\T}$ satisfies
\[
|d_{\T}(u_i, u_j) - d_{H^{\T}}(i,j) | \leq 2,
\]
%
%that $\T$ can be equivalently defined as (assuming $\T$ has $n+1$ vertices)
%\[
%(\{0, 1, \ldots, n\} / \sim_{{\T}}, d_{{\T}}),
%\]
where
\[
d_{H^{\T}}(i,j) = H^{\T}(i) + H^{\T}(j) - 2\min_{i \leq k \leq j} H^{\T}(k).
\]
%and $i \sim_{{\T}} j$ if and only if $d_{{\T}}(i,j)=0$. 
Clearly this discrepancy disappears in the scaling limit so in light of Proposition \ref{prop:annealed coding fn conv whole forest}, if the interval $[\beta_1, \beta_2]$ corresponds to an excursion of $\Ht$ above zero (this can be made sense of using excursion theory), then we define
\[
\Ta = ([\beta_1, \beta_2] / \sim_{\Ta}, d_{\Ta}),
\]
where 
\[
d_{\Ta} (s,t) = \Ht_s + \Ht_t - 2\inf_{s \leq r \leq t} \Ht_r
\]
and $s \sim_{\Ta} t$ if and only if $d_{\Ta} (s,t) =0$. Moreover, we equip $\Ta$ with the measure $\nu$, obtained as the image of Lebesgue measure on $[\beta_1, \beta_2]$ under the quotient operation. The root $\rho$ is equal to the projection of the point $\beta_1$. Note that the \textit{height} of $\Ta$ is defined as $\Height (\Ta) = \sup_{t \in [\beta_1, \beta_2]} d_{\Ta}(\beta_1, t)$, i.e. the maximal distance between any vertex and the root.

This can be defined formally using the \Ito excursion measure, the ``law'' under which excursions of $\Ht$ can be defined. It is in fact an infinite measure, but can be renormalised into a probability measure by conditioning the excursion to be large in an appropriate sense. In particular, the trees $\Ta^{\geq 1}$ and $\Ta^{H \geq 1}$ are respectively obtained by sampling an excursion of $\Ht$ conditioned to have a lifetime or height at least $1$, and applying the above construction. The trees $\Ta^{= 1}$ and $\Ta^{H= 1}$ are similarly respectively obtained by sampling an excursion of $\Ht$ conditioned to have a lifetime or height exactly equal $1$ - although this is a degenerate conditioning, this can also be formalised using excursion theory. We will not specifically need to use this excursion measure, so we refer to \cite[Chapter IV]{bertoin1996levy} for full details.

We also mention the notion of the local time at a certain level of $\Ta$. For $a>0$, Duquesne and Le Gall \cite{duquesne2002random} showed that one can construct a local time measure, supported on vertices at height $a$ in $\Ta$, and moreover such that the canonical measure $\nu$ on $\Ta$ satisfies
\begin{equation}\label{eqn:local a def}
\nu = \int_0^{\infty} \ell^a da,
\end{equation}
almost everywhere under the associated excursion measure. A vertex chosen according to $\ell^a$ can therefore be interpreted as a vertex chosen ``uniformly at level $a$'' in $\Ta$.

\subsubsection{Scaling limits of random trees}\label{sctn:coding of forests}
The significance of Proposition \ref{prop:annealed coding fn conv whole forest} is that this is enough to imply GHP convergence of individual trees conditioned to be large. We state the result below, and refer to \cite[Proposition 2.5.2]{duquesne2002random} for a proof. The reference in fact treats the case of conditioning a finite variance tree to have large height, but the proof of the more general statement below is the same, see the remark on \cite[page 62]{duquesne2002random}. We remark only that the key ingredient to replicate the proofs is the so-called local time support property for the limiting tree, which is well-known for $\Ta$ (see for example the remark of \cite[page 26]{duquesne2002random}).

\begin{proposition}\label{prop:GHP from contour}
Let $(T_i)_{i=1}^{\infty}$ be a sequence of trees, and let $H$ and $\Lambda$ denote their concatenated height and local time functions, as above. Suppose that there exist constants $c_2, c_3 \in (0, \infty)$ such that
   \begin{align}\label{eqn:forest contour conv}
      ({c_2n^{-\left(1-\frac{1}{\alpha}\right)}}H_{\ent{nt}},{c_3n^{-\frac{1}{\alpha}}}\Lambda_{\ent{nt}})_{t \ge 0} \underset{n \to +\infty}{\overset{(d)}{\rightarrow}} (\Ht_t, \Lt_t)_{t \ge 0},
  \end{align}
  jointly with respect to the uniform topology. Let $T_{\geq n}$ be the first tree in the sequence satisfying $|T| \geq n$, and let $T_{H \geq n}$ be the first tree in the sequence satisfying $\textsf{Height} (T) \geq n$. Then
  \[
  (T_{\geq n}, c_2n^{-\left(1-\frac{1}{\alpha}\right)}d_n, n^{-1}\nu_n, \rho_n) \overset{(d)}{\to} (\mathcal{T}^{\geq 1}_{\alpha},d_{\mathcal{\mathcal{T}_{\alpha}}},\nu_{\alpha},\rho_{\alpha})
  \]
  and 
   \[
  (T_{H \geq n}, n^{-1}d_n, (c_2n)^{-\left(\frac{\alpha}{\alpha-1}\right)}\nu_n, \rho_n) \overset{(d)}{\to} (\mathcal{T}^{H \geq 1}_{\alpha},d_{\mathcal{\mathcal{T}_{\alpha}}},\nu_{\alpha},\rho_{\alpha})
  \]
  with respect to the pointed Gromov-Hausdorff-Prokhorov topology, and where $\mathcal{T}^{\geq 1}_{\alpha}$ and $\mathcal{T}^{H \geq 1}_{\alpha}$ respectively denote the $\alpha$-stable tree conditioned to have total volume at least $1$ and height at least $1$.
\end{proposition}

In particular this applies under Assumption \ref{assn:offspring}. Moreover the conditioning can be made more precise.

\begin{proposition}
Under Assumption \ref{assn:offspring}, let $T_{= n}$ be a Galton--Watson tree with offspring law $\Pb$ conditioned to have exactly $n$ vertices, and let $T_{H = n}$ be a Galton--Watson tree with offspring law $\Pb$ conditioned on $\textsf{Height} (T) = n$.  Let $c_2$ be as in Proposition \ref{prop:annealed coding fn conv whole forest}. Then
  \[
  (T_{= n}, c_2n^{-\left(1-\frac{1}{\alpha}\right)}d_n, n^{-1}\nu_n, \rho_n) \overset{(d)}{\to} (\mathcal{T}^{= 1}_{\alpha},d_{\mathcal{\mathcal{T}_{\alpha}}},\nu_{\alpha},\rho_{\alpha})
  \]
  and 
   \[
  (T_{H =n}, n^{-1}d_n, (c_2n)^{-\left(\frac{\alpha}{\alpha-1}\right)}\nu_n, \rho_n) \overset{(d)}{\to} (\mathcal{T}^{H = 1}_{\alpha},d_{\mathcal{\mathcal{T}_{\alpha}}},\nu_{\alpha},\rho_{\alpha})
  \]
  with respect to the pointed Gromov-Hausdorff-Prokhorov topology, and where $\mathcal{T}^{= 1}_{\alpha}$ and $\mathcal{T}^{H = 1}_{\alpha}$ respectively denote the $\alpha$-stable tree conditioned to have total volume equal to $1$ and height equal to $1$.
\end{proposition}

Note that by comparing with \eqref{eqn:annealed GHP}, we see that $\gamma = c_2$.

\subsubsection{A useful fact}

We end this section with a useful lemma. We recall that under the annealed law $\mathbb{P}_{\alpha}$, the cluster $\C$ is just a critical Galton-Watson tree. It will later be useful to define the space $\C_{n, \eps}$ to be the ball of radius $k_{n, \eps}$ in $\C$, where $k_{n, \eps} = \inf \{r \geq 0: \sum_{i=0}^r Y_i \geq (1-\eps)n \}$, and similarly define $\C_{H, n, \eps}$ to be the ball of radius $(1-\eps)n$ in $\C$. The following result will be useful in order to refine the conditioning in Section \ref{section_conditioning_on_the_size}.

\begin{fact}\label{fact:annealed conv}
A consequence of the annealed pointed GHP convergence (under the conditioning $\# \mathcal{C} =n$ and $\Height (\C)=n$) is that, for any bounded Lipschitz function $F: \mathbb{K}_c \to \R$,
\begin{align*}
\lim_{\eps \downarrow 0} \sup_{n \geq 1} \bigg| \mathbb{E}_{\alpha}\bigg[F(\mathcal{C}_{n,\varepsilon}) \mid \# \mathcal{C} = n \bigg]  - \mathbb{E}_{\alpha}\bigg[F(\mathcal{C}) \mid \# \mathcal{C} =n \bigg] \bigg| &= 0, \\
\lim_{\eps \downarrow 0} \sup_{n \geq 1} \bigg| \mathbb{E}_{\alpha}\bigg[F(\mathcal{C}_{H, n,\varepsilon}) \mid \# \mathcal{C} = n \bigg]  - \mathbb{E}_{\alpha}\bigg[F(\mathcal{C}) \mid \# \mathcal{C} =n \bigg] \bigg| &= 0.
\end{align*}
\end{fact}

\section{The law of the total progeny}\label{sctn:size conv}

The aim of this section is to prove Theorem \ref{thm:ConvergenceOfSize}.%, or in other words that
% \begin{equation}
%     n^{\frac{1}{\alpha}}\Pt{|\C|>n}\to \K \bW\, ,
% \end{equation}
% $\bPb$-almost surely.

Before giving the proof, we note the following result that was proved but was not explicitly written in \cite{archer2023quenched} (it was written in a special case).

We set $\F_i=\sigma\rk{\bT_r\colon 0\le r\le i}$, i.e. the sigma algebra generated by the first $i$ levels of the tree. Conditionally on $\bT$ and given $u \in \bT_m$, let $T^{(u)}$ denote the subtree of $\bT$ emanating from and rooted at $u$. 

\begin{lemma}\label{lem:two vertices small prob}
Take $m \geq 1$. Conditionally on $\F_m$, let $(A_u)_{u \in \bT_m}$ be a sequence of events that are each respectively measurable with respect to $T^{(u)}$. For $u,v\in \bT_m$, we set $p_{u,v} =\Pt{A_u}\Pt{A_v}$, and $M=\sup_{u,v} \Eb{p_{u,v}}$. Then, for any $p < \frac{\alpha}{2} \leq 1$, there exists $C<\infty$ such that
\begin{align}\label{eqn:connection prob Markov}
 \Eb{\left(\sum_{\substack{u,v\in \bT_m \\ u\neq v}}\Pt{\rooot \leftrightarrow (u,v), A_u, A_v}\right)^p} 
        &\le C \mu^{m(1-p)}\Eb{{\overline{\bW}}^{2p}} M^p\, .
\end{align}
where $\overline{\bW}=\sup_{n \geq 0} \mu^{-n} |\bT_n|$.
\end{lemma}
\begin{proof}
The lemma was proved as part of the proof of \cite[Lemma 3.2]{archer2023quenched} in the case where $A_u$ is the event that $u$ connects to $\bT_n$ via a path of length $n-m$. The only proof ingredients are Jensen's inequality and Markov's inequality. Exactly the same proof works in the general case. Note that $\Eb{{\overline{\bW}}^{2p}}<\infty$ by Doob's $\mathrm{L}^p$ inequality and our choice of $p$.
\end{proof}

By monotonicity, it is sufficient to prove Theorem \ref{thm:ConvergenceOfSize} along a polynomial subsequence $n_k = \lfloor k^A \rfloor$, where $A$ is as large as we like.

\subsection{Lower bound in Theorem \ref{thm:ConvergenceOfSize}}
 For $u \in \bT$, let $C^{(u)} = \C \cap T^{(u)}$.
\begin{proposition}\label{prop:cluster size tail LB}
    We can choose $A$ large enough so that almost surely along the subsequence $(n_k)_{k \geq 1}$,
       \begin{equation}\label{eqn:prob size conv sub LB}
    \liminf_{k \to \infty} n_k^{\frac{1}{\alpha}} \Pt{|\C| \geq n_k} \geq \bW \K.
       \end{equation}
    Hence, by monotonicity,
    \begin{equation}\label{eqn:prob size conv all LB}
    \liminf_{n \to \infty} n^{\frac{1}{\alpha}} \Pt{|\C| \geq n} \geq \bW \K.
    \end{equation}
\end{proposition}

\begin{proof}
Note that \eqref{eqn:prob size conv all LB} is a straightforward consequence of \eqref{eqn:prob size conv sub LB} since if $n \in [n_{k-1}, n_{k}]$,
\[
n^{\frac{1}{\alpha}} \Pt{|\C| \geq n} \geq (1+o(1))n_{k}^{\frac{1}{\alpha}}\Pt{|\C| \geq n_{k}} \sim \bW \K.
\]
To prove \eqref{eqn:prob size conv sub LB}, we fix some small $\eps, \delta>0$ (we might reduce them later) and set $m=\lfloor \frac{1+\eps}{\alpha (\log \mu)} \log n \rfloor$ and write
\begin{align*}
\Pt{|\C| \geq n} &\geq \sum_{v \in \bT_m}  \Pt{\rooot \lr v, |C^{(v)}| \geq n} -  \sum_{u \neq v \in \bT_m} \Pt{  \rooot \lr (u,v), |C^{(u)}| \wedge |C^{(v)}| \geq n^{1-2\delta}} 
\end{align*}
(Note that the additional $\delta>0$ is not really necessary in the final probability above, but the bound we obtain will be useful for a later calculation.)

\textbf{First term.} We claim that, on rescaling by $n^{\frac{1}{\alpha}}$, the first term converges to $\bW\K$, $\bPb$-almost surely. To prove this, we first show that
\begin{align*}
S_m := \sum_{v \in \bT_m}  \left(\Pt{\rooot \lr v, |C^{(v)}| \geq n} - \mu^{-m} \pr{|\C| \geq n} \right) \to 0
\end{align*}
$\bPb$-almost surely. To see this, note that $\econdb{S_m}{\F_m}{} = 0$, $\bPb$-almost surely, and, provided that $n$ is sufficiently large,
\begin{align*}
\Varb{S_m | \F_m} = \mu^{-2m} \sum_{v \in \bT_m}  \Varb{\Pt{|C^{(v)}| \geq n}} &\leq \mu^{-2m} \sum_{v \in \bT_m}  \Eb{\Pt{|C^{(v)}| \geq n }} \\
&\leq 2 \mu^{-2m} |\bT_m| \K n^{-\frac{1}{\alpha}}.
\end{align*}
Combining with \eqref{eqn:W convergence}, we deduce that there exists a (random) constant $C<\infty$ such that, for all $n \geq 1$,
\begin{align*}
\Varb{S_m} = \Eb{\Varb{S_m | \F_m} } \leq C\K \mu^{-m} n^{-\frac{1}{\alpha}},
\end{align*}
and hence Chebyshev's inequality (and our choice of $m$) gives
\begin{align*}
\prb{|S_m|\geq \frac{1}{n^{\frac{1}{\alpha}}\log n}} \leq C\K n^{\frac{1}{\alpha}}(\log n)^2 \mu^{-m} = C\K n^{\frac{-\varepsilon}{\alpha}}(\log n)^2.
\end{align*}
Hence by Borel-Cantelli this goes to zero along the subsequence $(n_k)_{k \geq 1}$ provided we chose $A$ sufficiently large (which we can indeed do).

Applying \eqref{eqn:W convergence} and \eqref{eqn:annealed size prob convergence} it therefore follows that, along the subsequence $(n_k)_{k \geq 1}$,
\begin{align*}
n^{\frac{1}{\alpha}} \sum_{v \in \bT_m}  \Pt{\rooot \lr v, |C^{(v)}| \geq n} &= o(1) + n^{\frac{1}{\alpha}} \sum_{v \in \bT_m}  \mu^{-m} \pr{|C| \geq n} \\
&= o(1) + \mu^{-m} \bT_m n^{\frac{1}{\alpha}} \pr{|C| \geq n} \to \bW \K,
\end{align*}
$\bPb$-almost surely.

\textbf{Second term.} The second term can be dealt with using Lemma \ref{lem:two vertices small prob} and \eqref{eqn:annealed size prob convergence}, which imply that its $p^{th}$ moment (for $p \in (0,1)$) is upper bounded by
\begin{align*}
C \mu^{m(1-p)}\Eb{{\bW}^{2p}} n^{-\frac{2(1-2\delta)p}{\alpha}} = C n^{\frac{(1+\varepsilon)(1-p)}{\alpha}}\Eb{{\bW}^{2p}} n^{-\frac{2(1-2\delta)p}{\alpha}}.
\end{align*}
Take $\frac{1}{2} < p< \frac{\alpha}{2}$ and reduce $\varepsilon$ and $\delta$ if necessary so that $\kappa := (1-4\delta)p-(1+\varepsilon)(1-p)>0$. Then applying Markov's inequality (with the $p^{th}$ moment) gives
\begin{align*}
\prb{\Pt{\exists u \neq v \in T_m :  \rooot \lr (u,v), |C^{(u)}| \wedge |C^{(v)}| \geq n^{1-2\delta}} \geq \frac{1}{n^{\frac{1}{\alpha}} \log n}} \leq C' n^{-\kappa/2}.
\end{align*}
Hence by Borel-Cantelli this goes also to zero along the subsequence $(n_k)_{k \geq 1}$ provided we chose $A$ sufficiently large.
\end{proof}

\subsection{Upper bound in Theorem \ref{thm:ConvergenceOfSize}}

%Take some $\eps>0$ and some $\delta>0$ and set $m=\frac{1+\varepsilon}{\alpha} (\log n)/(\log \mu)$, so $\mu^m = n^{\frac{1+\varepsilon}{\alpha}}$. (For ease of reading, we will ignore complications caused by the fact that $m$ is not an integer, and note that everything follows through on taking floor and ceiling functions where appropriate.)

\begin{proposition}
    We can choose $A$ large enough so that almost surely along the subsequence $(n_k)_{k \geq 1}$
       \begin{equation}\label{eqn:prob size conv sub}
    \limsup_{k \to \infty} n_k^{\frac{1}{\alpha}} \Pt{|\C| \geq n_k} \leq \bW \K.
       \end{equation}
    Hence, by monotonicity,
    \begin{equation}\label{eqn:prob size conv all}
    \limsup_{n \to \infty} n^{\frac{1}{\alpha}} \Pt{|\C| \geq n} \leq \bW \K.
    \end{equation}
\end{proposition}
\begin{proof}
Again \eqref{eqn:prob size conv all} follows straightforwardly from \eqref{eqn:prob size conv sub}. We henceforth focus on proving \eqref{eqn:prob size conv sub}. Again set $m=\lfloor \frac{1+\eps}{\alpha (\log \mu)} \log n \rfloor$. For $u \in \bT$, we recall the notation $C^{(u)} = \C \cap T^{(u)}$, and let $N^{(u)}$ denote the number of siblings $v$ of $u$ to the left of $u$ that satisfy $|C^{(v)}| \geq n^{1-2\delta}$. Note that, by a union bound,
\begin{align*}
\Pt{|\C| \geq n} %&\leq \Pt{\exists v \in \bT_m : \rooot \lr v, |C^{(v)}| \geq n - n^{1-\delta}} \\
%&\qquad + \Pt{\sum_{v \in \bT_m} \mathbbm{1}\{ \rooot \lr v, |C^{(v)}| < n^{1-2\delta}, N^{(v)} \leq 1\} |C^{(v)}| \geq n^{1-2\delta}} \\
%&\qquad + \Pt{ |\C \cap (\bT_1 \cup \ldots \cup \bT_m)| \geq n^{1-2\delta}} \\
&\leq \Pt{ |\C \cap (\bT_1 \cup \ldots \cup \bT_m)| \geq n^{1-2\delta}} + \sum_{v \in \bT_m}  \Pt{\rooot \lr v, |C^{(v)}| \geq n - n^{1-\delta}} \\ 
&\qquad + \Pt{\sum_{v \in \bT_m} \mathbbm{1}\{ \rooot \lr v, |C^{(v)}| < n^{1-2\delta}, N^{(v)} \leq 1\} |C^{(v)}| \geq n^{1-\delta}} \\
&\qquad + \sum_{u \neq v \in \bT_m}\Pt{ \rooot \lr (u,v), |C^{(u)}| \wedge |C^{(v)}| \geq n^{1-2\delta}}.
\end{align*}
We will show that the second term concentrates on the desired quantity (up to an error of $o(n^{\frac{1}{\alpha}})$) and that the other terms are negligible (i.e. also $o(n^{\frac{1}{\alpha}})$) along the subsequence $(n_k)_{k\geq 1}$, provided we chose $A$ sufficiently large.

\textbf{First term.} For the first term note that by Markov's inequality we have 
\begin{align*}
\prb{\Et{ |\C \cap (\bT_1 \cup \ldots \cup \bT_m)|} \geq n^{\delta}} \leq n^{-\delta/2}
\end{align*}
and hence by Borel-Cantelli we can assume that $\Et{ |\C \cap (\bT_1 \cup \ldots \cup \bT_m)|} \leq n^{\delta}$ for all sufficiently large $n$ along the subsequence $(n_k)_{k \geq 1}$. On this latter event, the desired probability is upper bounded by $n^{-(1-3\delta)}$ by another application of Markov's inequality, which is $o(n^{1/\alpha})$ provided we took $\delta>0$ small enough in the first place.

\textbf{Second term.} The concentration of the second term follows exactly as in the proof of that of the first term in the proof of the lower bound (Proposition \ref{prop:cluster size tail LB}).

\textbf{Third term.} By Markov's inequality and Borel-Cantelli it is sufficient to show that 
\begin{align}\label{eqn:Markov sum geometric}
\econdb{\Pt{\sum_{v \in \bT_m} \mathbbm{1}\{ \rooot \lr v, |C^{(v)}| < n^{1-2\delta}, N^{(v)} \leq 1\} |C^{(v)}| \geq n^{1-\delta}}}{\F_m} \leq n^{-\frac{1+\delta}{\alpha}}.
\end{align}
The expectation in question is just the quantity
\begin{align*}
\sum_{V \subset \bT_m} \Pt{\C \cap \bT_m = V} 
\pr{\sum_{v \in V} \mathbbm{1}\{|C^{(v)}| < n^{1-2\delta}, N^{(v)} \leq 1\} |C^{(v)}| \geq n^{1-\delta}}.
\end{align*}
We will bound the latter probability uniformly over all choices of $V$. In particular, given any such $V$, consider the vertices of $V$ from left to right. We let $v_i$ denote the $i^{th}$ vertex in this ordering, and let $X_i$ denote the associated summand. Let ${N}$ denote the number of the terms in the sum and note that under $\Pb$, $N$ is stochastically dominated by the sum of two independent geometric random variables, and the sequence $(X_i)_{i \geq 1}$ is i.i.d.. Moreover, if $S_j$ denotes the partial sum with $j$ terms, i.e. $S_j = \sum_{i=1}^j X_i$, it follows that $S_{i+1}-S_i \leq n^{1-2\delta}$ for all $i$. Hence it follows from the memoryless property that for any $\lambda>2$,
\[
\pr{S_{N} \geq \lambda n^{1-2\delta}} \leq \pr{S_{N} \geq n^{1-2\delta}}^{\lambda/2}.
\]
In particular taking $\lambda = n^{\delta}$ and applying the tower property this easily implies \eqref{eqn:Markov sum geometric}, provided we can bound $\pr{S_{N} \geq n^{1-2\delta}}$ away from $1$. To do this, note that since $N$ can be upper bounded by the sum of two independent geometric random variables with parameter asymptotic to $\CC n^{\frac{1-2\delta}{\alpha}}$, it follows from standard results on scaling limits of stable variables that $n^{-(1-2\delta)}S_{N}$ converges in law to the value of a subordinator at a time $N$ which is equal in law to the sum of two independent \textsf{exp}($\CC$) random variables, and with jump measure proportional to $x^{1-\frac{1}{\alpha}}\mathbbm{1}\{x < 1\}$, and hence the probability in question converges to a constant in $(0,1)$.

\textbf{Fourth term.} The fourth term is the same as the second term in the proof of Proposition \ref{prop:cluster size tail LB}, hence we already showed it goes to $0$ under the rescaling.
% The fourth term can be dealt with using Lemma \ref{lem:two vertices small prob} and \eqref{eqn:annealed size prob convergence}, which imply that its $p^{th}$ moment (for $p \in (0,1)$) is upper bounded by
% \begin{align*}
% C \mu^{m(1-p)}\Eb{{\bW}^{2p}} n^{-\frac{2(1-2\delta)p}{\alpha}} = C n^{\frac{(1+\varepsilon)(1-p)}{\alpha}}\Eb{{\bW}^{2p}} n^{-\frac{2(1-2\delta)p}{\alpha}}.
% \end{align*}
% Take $\frac{1}{2} < p< \frac{\alpha}{2}$ and reduce $\varepsilon$ and $\delta$ if necessary so that $\kappa := (1-4\delta)p-(1+\varepsilon)(1-p)>0$. Then applying Markov's inequality (with the $p^{th}$ moment) gives
% \begin{align*}
% \prb{\Pt{\exists u \neq v \in T_m :  0 \lr (u,v), |C^{(u)}| \wedge |C^{(v)}| \geq n^{1-2\delta}} \geq \frac{1}{n^{\frac{1}{\alpha}} \log n}} \leq C' n^{-\kappa/2}.
% \end{align*}
% Hence by Borel-Cantelli this goes also to zero along the subsequence $(n_k)_{k \geq 1}$ provided we chose $A$ sufficiently large.

\end{proof}

\section{Scaling limit of the height function}\label{sctn:contour limit}

This section is dedicated to proving that the height function coding a forest of critical clusters converges under rescaling to its annealed limit (Theorem~\ref{convergence_contour_to_reflected_BM}).

We introduce $(\mathcal{C}^{i})_{i \ge 1}$ a family of random trees, such that conditionally on $\mathbf{T}$, the trees are i.i.d. and distributed as critical percolation clusters of the origin. We recall that the definitions for the height function and the local time were given in Section~\ref{section_coding_of_forests}. We denote by $X$ the height function associated to the random forest $(\mathcal{C}^{i})_{i \ge 1}$ and $\Lambda$ its local time at $0$. For $t \ge 0$, we introduce the notation
\begin{align*}%\label{renormalised_process}
   &(\beta^{\alpha,n}_t)_{t \ge 0} \text{ is the linear interpolation of } \frac{k}{n}\mapsto \frac{1}{n^{1-\frac{1}{\alpha}}}X_{k}, \\ &(\Upsilon^{\alpha,n}_t)_{t \ge 1} \text{ is the linear interpolation of } \frac{k}{n}\mapsto  \frac{1}{n^{\frac{1}{\alpha}}}\Lambda_{k}.
\end{align*}

The main theorem of the section is the following.

\begin{theorem}\label{convergence_contour_to_reflected_BM}
Assume that Assumption \ref{assumption} holds and let $\alpha \in (1,2]$ be as defined there. Then for $\mathbf{P}_{\alpha}$-almost every $\mathbf{T}$, under the quenched law $\mathbb{P}_{\mathbf{T}}$,
  \begin{align}\label{eqn:thm contour overall}
      (\beta^{\alpha,n}_t,\Upsilon^{\alpha,n}_t)_{t \ge 0} \underset{n \to +\infty}{\overset{(d)}{\rightarrow}} (c_2\Ht_t,c_3\mathbf{W}^{-1} \Lt_t)_{t \ge 0},
  \end{align}
where $\Ht,\Lt,c_2$ and $c_3$ are the processes and constants defined in \eqref{eqn:forest contour conv} (where the law of the Galton-Watson tree is given by $\Pf_{\alpha}$). This convergence holds jointly with respect to the uniform topology.
\end{theorem}

Let us give a short and intuitive explanation of the strategy to prove Theorem~\ref{convergence_contour_to_reflected_BM}. The first observation is that if we consider two percolation clusters $\mathcal{C}^1$ and $\mathcal{C}^2$ of $\mathbf{T}$ under the annealed law $\Pf_{\alpha}$, the intersection of the two clusters $\mathcal{C}^1 \cap \mathcal{C}^2$ is distributed as a subcritical Galton-Watson tree. It follows that, with probability close to $1$, if we cut the $n$ first clusters $\mathcal{C}^1,\cdots,\mathcal{C}^n$ at height $n^{\varepsilon}$ (where we chose $\varepsilon > 0$ arbitrarily small), and consider the subtrees emanating above this level, we obtain a family of independent trees all with law $\Pf_{\alpha}$. Thus, we first prove a version of Theorem~\ref{convergence_contour_to_reflected_BM} for a modified height process and local time associated to this family of subtrees obtained by cutting at level $n^{\varepsilon}$. This is the content of Proposition \ref{quenched_convergence_intermediate_process} in Section \ref{sctn:conv cutforest}.% In Proposition~\ref{CV_annealed} we prove an annealed version. Then, Proposition~\ref{variance_control_imply_convergence_distribution} and Proposition~\ref{convergence_of_series} extend the annealed to the quenched result.\\

Using this result, we can prove Theorem~\ref{convergence_contour_to_reflected_BM} using the fact that the part that we removed when cutting is small enough and by connecting the overall local time to the modified local time for the cutforest. Indeed, if one considers the first $n$ clusters $\mathcal{C}^1,\dots,\mathcal{C}^n$,  by Theorem~\ref{thm:ConvergenceOfSize} the number of edges is typically of order $n^{\alpha}$. However, the number of edges below level $n^{\varepsilon}$ is typically of order $n^{1+\varepsilon}$. Taking $\varepsilon$ small enough, we deduce that the removed part is small compared to the entire forest and so the overall height process should be typically close to the modified one. Concerning the local time, the idea is that the number of vertices at height $n^{\varepsilon}$ in the forest of trees $\mathcal{C}^1,\dots,\mathcal{C}^m$ should be typically of order $\mathbf{W}m$. Thus, one should be able to move from the modified local time to the local time of the height process by simply dividing by $\mathbf{W}$. This proof appears in Section \ref{sctn:convergence original height}.

{Note that the trees should technically be cut at an integer generation, so level $\lfloor n^{\eps} \rfloor$ rather than simply $n^{\eps}$. To ease notation (and since it has no effect on the argument), we have omitted this floor and ceiling notation throughout the section.}

\begin{remark}\label{rmk:Liipschitz single}
In order to prove \eqref{eqn:thm contour overall}, we will need to show that, $\bPb$-almost surely, for all non-negative bounded Lipschitz functions $F: \mathcal{C}([0,T],\mathbb{R})^2 \longrightarrow \mathbb{R}_{+}$,  
\[
 \mathbb{E}_{\bT}[F(\beta^{\alpha,n}_t,\Upsilon^{\alpha,n}_t)_{t \ge 0} ]\underset{n \to +\infty}{\overset{(d)}{\rightarrow}} \mathbb{E}_{\alpha}\left[{F(c_2\Ht_t,c_3\mathbf{W}^{-1} \Lt_t)_{t \ge 0}}\right].
\]
To do this, it will actually be sufficient to prove the convergence for a \textit{single} arbitrary non-negative bounded Lipschitz function $F$. Indeed, this allows us to prove tightness of the processes using an appropriate countable sequence of functions and a standard tightness criterion for the uniform topology. Once we establish that the process is tight, we can restrict to a compact subspace $\hat{\mathbb{K}} \subset \mathcal{C}([0,T],\mathbb{R})^2$. The space of bounded Lipschitz functions $F: \hat{\mathbb{K}} \longrightarrow \mathbb{R}_{+}$ is then separable, which means that the desired claim will again follow by testing only a countable number of functions $F$. This type of argument is written out in some detail in the proof of \cite[Lemma $4.1$]{Bolthausen} and we refer there for the details.
\end{remark}

\subsection{Convergence of the cutforest}\label{sctn:conv cutforest}

We start by giving the setup and some first observations. For $T$ a rooted tree and $k \ge 0$, we denote by $T^{\uparrow k}$ the subgraph induced by the vertices of $T$ with height at least $k$. Note this is not necessarily connected ; we also decompose $T^{\uparrow k} =\displaystyle  \bigsqcup_{j=1}^{Y_k}{T_{j}^{\uparrow k}}$ where $Y_k$ denotes the number of vertices at generation $k$ and $T_{j}^{\uparrow k}$ corresponds to the $j^{th}$-connected component of $T^{\uparrow k}$ from left to right. For a cluster $\C^i$ as above, we write $\mathcal{C}^{i,\uparrow k}_{j}$ in place of $(\mathcal{C}^{i})^{\uparrow k}_{j}$ and $Y^i_k$ to denote the size of generation $k$ in $\C^i$.

\begin{figure}[H]
    \centering
    \includegraphics[scale =0.2]{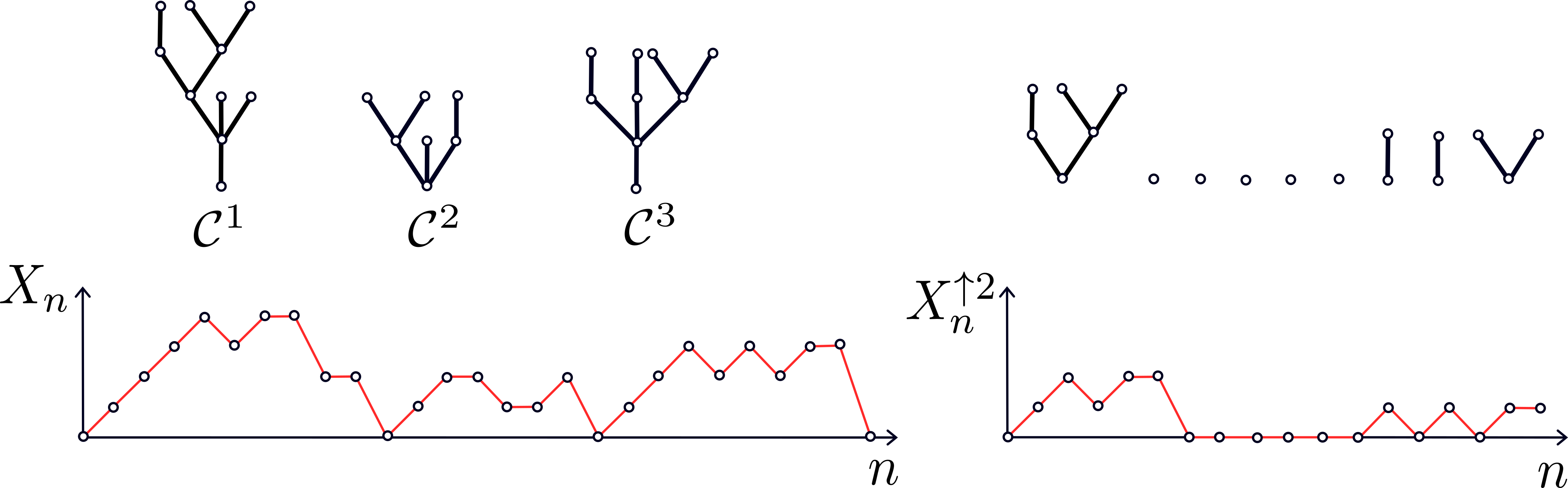}
    \caption{On the left side, we represent $\mathcal{C}^{1},\mathcal{C}^{2},\mathcal{C}^{3}$ on top and the concatenation of their height functions on the bottom. On the right side, we represent the trees $\mathcal{C}^{1,\uparrow 2}_1,\mathcal{C}^{1,\uparrow 2}_2,\mathcal{C}^{1,\uparrow 2}_3,\mathcal{C}^{2,\uparrow 2}_1,\mathcal{C}^{2,\uparrow 2}_2,\mathcal{C}^{2,\uparrow 2}_3,\mathcal{C}^{3,\uparrow 2}_1,\mathcal{C}^{3,\uparrow 2}_2,\mathcal{C}^{3,\uparrow 2}_3$ and their concatenated contour functions.}
\label{trees_cut_trees_and_contour_functions}
\end{figure}

 For $k \ge 0$, we define $X^{\uparrow k}$ to be the process which concatenates the height functions associated to the trees $((\mathcal{C}^{i,\uparrow k}_{j})_{1\le j \le Y^{i}_k}))_{i \ge 0}$ using the lexicographical order on $(i,j)$. In particular we have $X^{\uparrow 0} = X$. Similarly, we define $\Lambda^{\uparrow k}$ as the local time at $0$ associated to $X^{\uparrow k}$. See Figure \ref{trees_cut_trees_and_contour_functions} for an illustration.

For the rest of this section we fix $\varepsilon = \frac{1}{5}\frac{\alpha - 1}{\alpha}$. For $t \ge 0$, we introduce the notation
\begin{align}\label{renormalised_process_aboce_neeps}
\begin{split}
     &(\beta^{\alpha,\uparrow n^{\varepsilon}}_t)_{t \ge 0} \text{ is the linear interpolation of } \frac{k}{n}\mapsto \frac{1}{n^{1-\frac{1}{\alpha}}}X^{\uparrow n^{\varepsilon}}_{k}, \\ &(\Upsilon^{\alpha,\uparrow n^{\varepsilon}}_t)_{t \ge 0} \text{ is the linear interpolation of } \frac{k}{n}\mapsto  \frac{1}{n^{\frac{1}{\alpha}}}\Lambda^{\uparrow n^{\varepsilon}}_{k}.
\end{split}
\end{align}
\noindent For the rest of this section, we also fix a finite $r > \frac{3}{2} + \frac{1}{4(\alpha -1)}$ (its precise value is not important, but $n^r$ will be a convenient upper bound for the number of subtrees we need to consider). For $i,j, n \ge 1$, we introduce the events
\begin{align}\label{eqn:An def}
\mathcal{A}^{n}_{i,j} = \bigg\{\textsf{Height}(\mathcal{C}^i \cap \mathcal{C}^j) \le n^{\varepsilon} \bigg\} \text{ and } \mathcal{A}^{n} = \displaystyle \bigcap_{i\neq j \le n^r }\mathcal{A}_{i,j}^{n}.
\end{align}
For $i \ne j$, the tree $\mathcal{C}^{i}\cap \mathcal{C}^j$ is distributed as a subcritical Galton-Watson tree under the annealed law, thus we have the following bound
\begin{align}\label{bound_proba_An}
    \Pf_{\alpha}((\mathcal{A}^n)^{c}) \le Ce^{-cn^{\varepsilon}},
\end{align}
where $C,c > 0$ are constants that only depend on $\mathbf{P}_{\alpha}$.\\
We also introduce the event $\mathcal{B}^n = \displaystyle \left\{\#\{i \leq n^r : \textsf{Height}(\mathcal{C}^i)\ge n^{\frac{1}{4}}\} \geq n\right\}$. Using \cite[Theorem $1.2$]{archer2023quenched} and classical concentration inequalities for binomial random variables, it is easy to prove that for $\mathbf{P}_{\alpha}$-almost every $\mathbf{T}$, we have that
\begin{align}\label{proba_Bn_goes_to_1}
    \Pf_{\mathbf{T}}(\mathcal{B}^n) \underset{n \to +\infty}{\rightarrow} 1.
\end{align}
\begin{proposition}\label{quenched_convergence_intermediate_process}
Assume that Assumption \ref{assumption} holds and let $\alpha \in (1,2]$ be as defined there. For $\mathbf{P}_{\alpha}$-almost every $\mathbf{T}$, we have under the quenched law $\Pf_{\mathbf{T}}$,
\begin{align}\label{eqn:quenched cut contour convergence}
   (\beta^{\alpha,\uparrow n^{\varepsilon}}_t,\Upsilon_t^{\alpha,\uparrow n^{\varepsilon}})_{t \in [0,T]} \underset{n \to +\infty}{\overset{(d)}{\rightarrow}} (c_2\Ht_t,c_3\Lt_t)_{t \in [0,T]} ,
\end{align}
jointly with respect to the uniform topology.
\end{proposition}
\begin{proof}
The first step is to prove that \eqref{eqn:quenched cut contour convergence} holds under the annealed law (this is not completely trivial since distinct clusters are not independent under the annealed law), and then use a second moment argument to argue that the quenched process behaves like the annealed process.

\paragraph{Step 1: annealed convergence.} First, observe that on the event $\mathcal{A}^n$, and conditionally on the cut sizes $(Y^1_{n^{\varepsilon}},\dots,Y^{n^r}_{n^{\varepsilon}})$, the family of upper trees $((\mathcal{C}^{i,\uparrow n^{\varepsilon}}_{j})_{1\le j \le Y^{i}_{n^\varepsilon}})_{1\le i \le n^r}$ consists of i.i.d.\ variables distributed as $\mathcal{C}$ under $\Pf_{\alpha}$. Furthermore, on the event $\mathcal{B}^n$, the total size of these trees is at least $(n^{\frac{1}{4}}-n^{\varepsilon})n$, which ensures that the rescaled concatenated process covers the interval $[0,T]$ for large $n$.

Since $\Pf_{\alpha}(\mathcal{A}^n \cap \mathcal{B}^n) \to 1$, the truncated process $(\beta^{\alpha,\uparrow n}_t,\Upsilon_t^{\alpha,\uparrow n})_{t \in [0,T]}$ coincides with high probability with the height and local time process of a forest of i.i.d.\ Galton-Watson trees. Consequently, we obtain the annealed convergence:
\begin{align}\label{conv_annealed_truncated}
    (\beta^{\alpha,\uparrow n^{\varepsilon}}_t,\Upsilon_t^{\alpha,\uparrow n^{\varepsilon}})_{t \in [0,T]} \underset{n \to +\infty}{\overset{(d)}{\longrightarrow}} (c_2\Ht_t,c_3\Lt_t)_{t \in [0,T]} \quad \text{under } \Pf_{\alpha}.
\end{align}

\paragraph{Step 2: quenched convergence.} To upgrade this to a quenched convergence, we introduce an intermediate process denoted by $(U^n_t,V^n_t)_{t \in [0,T]}$. This process is constructed as follows:
\begin{itemize}
	\item[$\bullet$] First, construct the height function and local time for the forest consisting of the first $n^r$ truncated clusters $\mathcal{C}^{1,\uparrow n^{\varepsilon}},\dots,\mathcal{C}^{n^r,\uparrow n^{\varepsilon}}$, and then extending the forest with independent Galton-Watson trees distributed as $\mathcal{C}$ under $\Pf_{\alpha}$.
	\item[$\bullet$] Then, the process $(U^n_t,V^n_t)_{t \in [0,T]}$ is defined as the linear interpolation of the above pair of height function and local time, rescaled as in \eqref{renormalised_process_aboce_neeps}.
\end{itemize}
 Since $\Pf_{\mathbf{T}}(\mathcal{B}^n)\to 1$ and since the processes coincide on $\mathcal{B}^n$, the convergence \eqref{conv_annealed_truncated} implies that:
\begin{align}\label{eqn:conv annealed replace}
  (U^n_t,V^n_t)_{t \in [0,T]} \underset{n \to +\infty}{\overset{(d)}{\longrightarrow}} (c_2\Ht_t,c_3\Lt_t)_{t \in [0,T]} \quad \text{under } \Pf_{\alpha}.
\end{align}
We now aim to show that this convergence also holds in the quenched setting. To this end we take a non-negative bounded Lipschitz function $F: \mathcal{C}([0,T],\mathbb{R})^2 \longrightarrow \mathbb{R}_{+}$. By \eqref{eqn:conv annealed replace} we have that $\mathbf{E}_{\alpha}\left[ \mathbb{E}_{\mathbf{T}}\left[ F((U^n_t,V^n_t)_{t \in [0,T]})\right]\right] \to \mathbb{E}[F(c_2\Ht_t,c_3\Lt_t)_{t \in [0,T]}]$. We now control the variance of $\mathbb{E}_{\mathbf{T}}\left[ F((U^n_t,V^n_t)_{t \in [0,T]})\right]$. Let $(U^{n,1}_t,V^{n,1}_t)_{t \in [0,T]}$ and $(U^{n,2}_t,V^{n,2}_t)_{t \in [0,T]}$ be two independent copies of the process under the quenched measure $\Pf_{\mathbf{T}}$. Specifically, the first copy is generated using the clusters $\mathcal{C}^1,\dots,\mathcal{C}^{n^r}$ and the second using $\mathcal{C}^{n^r +1},\dots,\mathcal{C}^{2n^r}$. Then $\mathbf{Var}_{\alpha}\left( \mathbb{E}_{\mathbf{T}}\left[ F((U^n_t,V^n_t)_{t \in [0,T]})\right]\right) $ is equal to
\begin{align}\label{eqn:var decomp 1}
&\mathbb{E}_{\alpha}\left[  F((U^{n,1}_t,V^{n,1}_t)_{t \in [0,T]}) F((U^{n,2}_t,V^{n,2}_t)_{t \in [0,T]})\right] - \mathbb{E}_{\alpha}\left[ F((U^{n,1}_t,V^{n,1}_t)_{t \in [0,T]})\right]^2.
\end{align}
Crucially, on the event 
\[ \mathcal{E}_n = \bigcap_{1 \le i \le n^r < j \le 2n^r} \left\{ \mathrm{Height}(\mathcal{C}^i\cap \mathcal{C}^j) \le n^\varepsilon \right\}, \]
the trees used to construct the two copies explore disjoint parts of the underlying tree $\mathbf{T}$ above level $n^\varepsilon$. Therefore, under $\Pf_{\alpha}$, conditionally on $\mathcal{E}_n$, the two processes are independent. In particular the corresponding contribution to the left hand side of \eqref{eqn:var decomp 1} factorises and there is no net contribution to the variance.
Since the intersection of two independent critical percolation clusters is subcritical, we have the exponential bound
\begin{align*}
\Pf_{\alpha}(\mathcal{E}_n^c) \le \exp(-cn),
\end{align*}
for some $c > 0$, and hence
\begin{align*}
\mathbf{Var}_{\alpha}\left( \mathbb{E}_{\mathbf{T}}\left[ F((U^n_t,V^n_t)_{t \in [0,T]})\right]\right) \le C\exp(-cn),
\end{align*}
for some constant $C$ depending on $\alpha$ and $F$. Since this upper bound is summable, it follows from the Borel-Cantelli lemma (or \cite[Lemma $4.1$]{Bolthausen}) that for $\mathbf{P}_{\alpha}$-almost every $\mathbf{T}$, $ \mathbb{E}_{\mathbf{T}}\left[ F((U^n_t,V^n_t)_{t \in [0,T]})\right] \underset{n \to +\infty}{{\longrightarrow}}  \mathbb{E}_{\alpha}\left[ F( (c_2\Ht_t,c_3\Lt_t)_{t \in [0,T]} )\right]$ and hence that (see Remark \ref{rmk:Liipschitz single}):
\begin{align*}
   (U^n_t,V^n_t)_{t \in [0,T]} \underset{n \to +\infty}{\overset{(d)}{\longrightarrow}} (c_2\Ht_t,c_3\Lt_t)_{t \in [0,T]} \quad \text{under } \Pf_{\mathbf{T}}.
\end{align*}
Finally, using the fact that $(\beta^{\alpha,\uparrow n^{\varepsilon}}_t,\Upsilon_t^{\alpha,\uparrow n^{\varepsilon}})_{t \in [0,T]}$ and $(U^n_t,V^n_t)_{t \in [0,T]}$ coincide with high probability under $\Pf_{\mathbf{T}}$, we conclude that for almost every $\mathbf{T}$, $ \mathbb{E}_{\mathbf{T}}\left[ F((U^n_t,V^n_t)_{t \in [0,T]})\right]$:
\begin{align*}
 (\beta^{\alpha,\uparrow n^{\varepsilon}}_t,\Upsilon_t^{\alpha,\uparrow n^{\varepsilon}})_{t \in [0,T]} \underset{n \to +\infty}{\overset{(d)}{\longrightarrow}} (c_2\Ht_t,c_3\Lt_t)_{t \in [0,T]} \quad \text{under } \Pf_{\mathbf{T}}.\qedhere
\end{align*}
\end{proof}

\subsection{Proof of Theorem~\ref{convergence_contour_to_reflected_BM}}\label{sctn:convergence original height}

It remains to transfer the result back to the original sequence of trees, i.e. not cut at level $n^{\eps}$. We now turn to this, thus concluding the proof of Theorem~\ref{convergence_contour_to_reflected_BM}.

%For $k \ge 0$ and $i \ge 1$, we define $Y^{i}_k = \# \mathcal{C}^i_{k}$, the number of vertices at level $k$ of the $i^{th}$ tree $\C^i$. %moved

\begin{proof}[Proof of Theorem~\ref{convergence_contour_to_reflected_BM}]
Take $\alpha \in (1,2]$ as in the statement and fix and $T > 0$. Let us prove the convergence on $[0,T]$. Fix $\mathbf{T}$ such that \cite[Theorem 1.2,Theorem 1.3]{archer2023quenched}, Theorem~\ref{thm:ConvergenceOfSize} as well as the distributional convergence from Proposition~\ref{quenched_convergence_intermediate_process} all hold. Thus, we have
\begin{align}\label{distrib_Cv_to_use}
    (\beta^{\alpha,\uparrow n^{\varepsilon}}_{t},\Upsilon^{\alpha,\uparrow n^{\varepsilon}}_{t})_{0 \le t \le T} \underset{n \to +\infty}{\overset{(d)}{\rightarrow}} (c_2\Ht_t,c_3\Lt_t)_{0 \le t \le T}\quad \text{under } \Pf_{\mathbf{T}}.
\end{align}

We define events
\[
\mathcal{D}_n = \displaystyle \bigg\{\somme{i=1}{n^{1-2\varepsilon}}{\somme{k=1}{n^{\varepsilon}}{Y^i_k}} \le n^{1-\frac{\varepsilon}{2}}\bigg\}
\text{ and } \mathcal{E}_n = \displaystyle\bigg\{\somme{i=1}{n^{1-2\varepsilon}}{\# \mathcal{C}^{i} \ge n^{1+\frac{\varepsilon}{2}}}\bigg\}.
\]
We let the reader verify that using Theorem~\ref{thm:ConvergenceOfSize} and recalling that $\eps=\frac{\alpha-1}{5\alpha}$, we have 
\begin{align}\label{proba_D_cap_E_go_to_1}
    \Pf_{\mathbf{T}}(\mathcal{D}_n \cap \mathcal{E}_n) \underset{n \to +\infty}{\rightarrow}1,
\end{align}
$\bPb$-almost surely. For $n$ large enough, on $\mathcal{D}_n \cap \mathcal{E}_n$, there are fewer than $n^{1-\frac{\varepsilon}{2}}$ vertices with height less than $n^{\varepsilon}$ and more than $nT$ vertices with height more than $n^{\varepsilon}$ among the trees $\mathcal{C}^1,\cdots,\mathcal{C}^{n^{1-2\varepsilon}}$.\\

\noindent We will now couple $X$ and $X^{\uparrow n^{\varepsilon}}$ in the natural way. First, we claim that on $\mathcal{D}_n \cap \mathcal{E}_n$ we can define a function $\phi: \{0,\cdots,nT\} \to \{0,\cdots,nT\}$ such that
\begin{align}\label{trajectory_normal_to_nepsilon}
    \forall k\in \{0,\cdots,nT\},\hspace{0.3cm} |X^{}_k- X^{\uparrow n^{\varepsilon}}_{\phi(k)}| \le n^{\varepsilon} \text{ and }|k-\phi(k)| \le n^{1-\frac{\varepsilon}{2}}.
\end{align}
Indeed, for $0 \le k \le nT$, define $\phi'(k) = \inf\{i \ge k \text{ :}\text{ }X_i \ge n^{\varepsilon}\}$. Then on $\mathcal{D}_n \cap \mathcal{E}_n$, we have $\phi'(k) -k\le n^{1-\frac{\varepsilon}{2}}$ and there exists $ \phi(k) < \phi'(k)$ such that $X^{\uparrow n^{\varepsilon}}_{\phi(k)} = X_{\phi'(k)} - n^{\varepsilon}$ and $|\phi'(k)-\phi(k)|\le n^{1-\frac{\varepsilon}{2}}$. Indeed, denote by $x^k$ the $\phi'(k)^{th}$ vertex visited in the lexicographical exploration of the forest $(\mathcal{C}^i)_{i \ge 0}$. By definition of $\phi'(k)$, the vertex $x^k$ has height at least $n^{\varepsilon}$, thus it belongs to the forest $((\mathcal{C}^{i,\uparrow n^{\varepsilon}}_{j})_{1 \le j \le Y^i_k})_{i\ge 0}$. Then, one can choose $\phi(k)$ as the time $x^k$ is visited in the lexicographical exploration of the forest $((\mathcal{C}^{i,\uparrow n^{\varepsilon}}_{j})_{1 \le j \le Y^i_k})_{i\ge 0}$. On $\mathcal{D}_n \cap \mathcal{E}_n$, we have $|\phi(k) - \phi'(k)|\le n^{1-\frac{\varepsilon}{2}}$. One can use Figure~\ref{contour} for a visual explanation.\\
 Using a similar argument, on $\mathcal{D}_n \cap \mathcal{E}_n$ we can write
\begin{align}\label{connect_local_time}
    \forall k\in \{0,\cdots,nT\},\hspace{0.3cm} \bigg|\somme{m=1}{\Lambda_{k}}{Y^{m}_{n^{\varepsilon}}} - \Lambda^{\uparrow n^{\varepsilon}}_{\phi(k)}\bigg| \le Y^{\Lambda_k}_{n^{\varepsilon}}.
\end{align}
\\

\begin{figure}[H]
    \centering
    \includegraphics[scale =0.15]{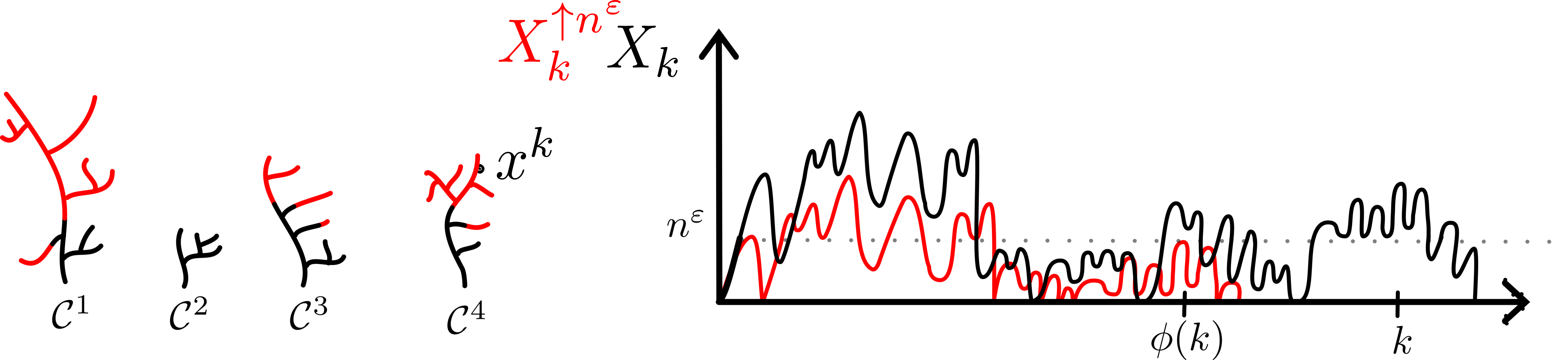}
    \caption{On the left we represent the trees $\mathcal{C}^i$ for $i \in \{1,2,3,4\}$. The red part corresponds to the vertices at height larger than $n^{\varepsilon}$. On the right we represent in black the height function $X$ of the trees $\mathcal{C}^i$. We represent in red the height function $X^{\uparrow n^{\varepsilon}}$ of the red part of the trees. On the event $\mathcal{D}_n \cap \mathcal{E}_n$, the size of the black part of the left trees is bounded by $n^{1-\varepsilon/2}$. This is sufficient to see that $|k - \phi(k)| \le n^{1-\varepsilon/2}$. It is similarly clear that $|\somme{m=1}{\Lambda_k}{Y^{m}_{n^{\varepsilon}}}-\Lambda^{\uparrow n^{\varepsilon}}_{\phi(k)}| \le Y^{\Lambda_k}_{n^{\varepsilon}}$.} 
    \label{contour}
\end{figure}
We will proceed in several steps to prove the theorem, starting with the height function, which is the simplest.
\paragraph{Convergence of the height function.}  On $\mathcal{D}_n \cap \mathcal{E}_n$, it follows from \eqref{trajectory_normal_to_nepsilon} that we have
\begin{align} \label{from_X_to_Xnepsilon}
   \sup_{0\le t \le T}|\beta^{\alpha,n}_t -\beta^{\alpha,\uparrow n^{\varepsilon}}_t | \le n^{\varepsilon-1+\frac{1}{\alpha}} + \sup_{\substack{0\le t,t^{'} \le 2T\\|t-t^{'}|\le n^{-\frac{\varepsilon}{2}}}}|\beta^{\alpha,\uparrow n^{\varepsilon}}_t -\beta^{\alpha,\uparrow n^{\varepsilon}}_{t^{'}} |.
\end{align}
Using \eqref{distrib_Cv_to_use}, it is clear that under $\Pf_{\mathbf{T}}$, the right side of the inequality tends to $0$ in probability as $n\to +\infty$. Together with \eqref{proba_D_cap_E_go_to_1}, this gives
\begin{align}\label{convergence_trajectory}
    \sup_{0\le t \le T}|\beta^{\alpha,n}_t -\beta^{\alpha,\uparrow n^{\varepsilon}}_t | \underset{n \to +\infty}{\overset{(\mathbb{P}_{\bT})}{\rightarrow}}0.
\end{align}
\paragraph{Convergence of the local time.} Now let us prove that
\begin{align}\label{to_control2}
    \sup_{0\le t \le T}|\mathbf{W}\Upsilon^{\alpha,n}_t - \Upsilon^{\alpha,\uparrow n^{\varepsilon}}_t| \underset{n \to +\infty}{\overset{(\mathbb{P}_{\bT})}{\rightarrow}}0.
\end{align}
The proof of this is quite involved and is divided into two steps: first we control the number of individuals appearing in generation $n^{\eps}$ in the first $\ent{n^{\frac{1}{\alpha}}t}$ subtrees. Then we compare this to the local time at zero over the same time period, and show that the two quantities are comparable.
\paragraph{Step 1: controlling the size of generation $n^{\eps}$.}
We show that for any $R > 0$, we have 
\begin{align*}% \label{Local_time_0_to_nepsilon}
   \displaystyle  \sup_{0 \le t \le R} \bigg[\mathbf{W}t - \frac{ \somme{m=1}{\ent{n^{\frac{1}{\alpha}}t}}{Y^{m}_{n^{\varepsilon}}}}{n^{\frac{1}{\alpha}}}\bigg]\underset{n \to +\infty}{\overset{(\Pf_{\bT})}{\to}}0.
\end{align*}
Set $\beta=\frac{1}{\alpha-1}$. We start by considering a fixed time $t >0$. We know from \cite[Theorem 1.3]{archer2023quenched} that we have
\begin{align*}
    (n^{-\beta}Y_n|Y_n > 0) \underset{n \to +\infty}{\overset{(d)}{\to}} Y,
\end{align*}
where $Y$ is an $\alpha$-stable random variable with expectation $C_{\alpha}^{-1}$ and Laplace transform $\psi$ given by
\begin{align*}
    \psi(\theta) = 1 - C_{\alpha}^{-1}\theta(1+(C_{\alpha}\theta)^{\alpha-1})^{-\beta},
\end{align*}
and where $C_{\alpha}$ is the constant defined in \eqref{eqn:slack prob convergence}. Let us introduce $(Z_{m,n})_{m \ge 0}$ a family of i.i.d. random variables distributed as $(n^{-\beta \varepsilon}Y_{n^{\varepsilon}}|Y_{n^{\varepsilon}} > 0)$. Then it is clear that we have
\begin{align}\label{eqn:rewrite sum Y}
     \frac{\displaystyle \somme{m=1}{\ent{n^{\frac{1}{\alpha}}t}}{Y^{m}_{n^{\varepsilon}}}}{n^{\frac{1}{\alpha}}} \overset{(d)}{=}  \frac{\displaystyle \somme{m=1}{\ent{n^{\frac{1}{\alpha}}t}}{\indi{Y^{m}_{n^{\varepsilon}}>0}}Z_{m,n}}{n^{\frac{1}{\alpha}-\beta \varepsilon}}.
\end{align}
Using \cite[Theorem 1.2]{archer2023quenched}, we let the reader verify (for example using a second moment argument) that we have
\begin{align}\label{binomial_concentration}
    \frac{\displaystyle\somme{m=1}{\ent{n^{\frac{1}{\alpha}}t}}{\indi{Y^{m}_{n^{\varepsilon}}>0}}}{\mathbf{W}C_{\alpha} n^{\frac{1}{\alpha}-\beta \varepsilon}t} \underset{n \to +\infty}{\overset{(\mathbb{P}_{\bT})}{\to}}1.
\end{align}
Introduce $\mathcal{I} = \{m\in \{1,\cdots,n^{\frac{1}{\alpha}}t\} \text{ } : \text{ } \indi{Y^{m}_{n^{\varepsilon}}>0} = 1\}$. Conditionally on $(\indi{Y^{m}_{n^{\varepsilon}}>0})_{m \ge 0}$, we write
\begin{align*} %\label{non_zero_term}
    \displaystyle \somme{m=1}{\ent{n^{\frac{1}{\alpha}}t}}{\indi{Y^{m}_{n^{\varepsilon}}>0}}Z_{m,n} = \somme{m \in \mathcal{I}}{}{Z_{m,n}} \overset{(d)}{=} \somme{m=1}{|\mathcal{I}|}{Z_{m,n}}.
\end{align*}
Using \cite[Theorem $1.2$]{archer2023quenched}, we see that the family $(Z_{n,m})_{n,m\ge 0}$ has bounded first moment. Recall also that $Z_{n,1} \underset{n \to +\infty}{\overset{(d)}{\rightarrow}} Y$ with  $\mathbb{E}_{\bT}[Y]=C_{\alpha}^{-1}$ and  $$\Et{Z_{n,1}} = \frac{n^{-\beta \varepsilon}\Et{Y_{n^{\varepsilon}}}}{\Pt{Y_{n^{\varepsilon}} > 0}} = \frac{n^{-\beta \varepsilon}\mu^{-n^{\varepsilon}}|\mathbf{T}_{n^{\varepsilon}}|}{\Pt{Y_{n^{\varepsilon}} > 0}}\underset{n\to +\infty}{\overset{(d)}{\rightarrow}} C_{\alpha}^{-1}.$$ By \eqref{binomial_concentration} we have $\displaystyle \frac{|\mathcal{I}|}{\mathbf{W}C_{\alpha} n^{\frac{1}{\alpha}-\beta \varepsilon}t} \underset{n \to +\infty}{\overset{(\mathbb{P}_{\mathbf{T}})}{\rightarrow}} 1$. Using Proposition~\ref{technical_proposition}, we deduce that
\begin{align*}
    \frac{\displaystyle \somme{m=1}{\ent{n^{\frac{1}{\alpha}}t}}{\indi{Y^{m}_{n^{\varepsilon}}>0}}Z_{m,n}}{n^{\frac{1}{\alpha}-\beta \varepsilon}} \underset{n \to +\infty}{\overset{(\mathbb{P}_{\mathbf{T}})}{\rightarrow}} \mathbf{W}t.
\end{align*}
\noindent By \eqref{eqn:rewrite sum Y}, we therefore have that
\begin{align*}
     \frac{\somme{m=1}{\ent{n^{\frac{1}{\alpha}}t}}{Y^{m}_{n^{\varepsilon}}}}{n^{\frac{1}{\alpha}}} \underset{n \to +\infty}{\overset{(\mathbb{P}_{\mathbf{T}})}{\to}}\mathbf{W}t.
\end{align*}
Using the fact that $\displaystyle \bigg(\frac{ \somme{m=1}{\ent{n^{\frac{1}{\alpha}}t}}{Y^{m}_{n^{\varepsilon}}}}{n^{\frac{1}{\alpha}}}\bigg)_{t \ge 0}$ is increasing in $t$, we deduce that for any $R \in (0, \infty)$ we have
\begin{align} \label{Local_time_0_to_nepsilon}
   \displaystyle  \sup_{0 \le t \le R} \bigg[\mathbf{W}t - \frac{ \somme{m=1}{\ent{n^{\frac{1}{\alpha}}t}}{Y^{m}_{n^{\varepsilon}}}}{n^{\frac{1}{\alpha}}}\bigg]\underset{n \to +\infty}{\overset{(\Pf_{\bT})}{\to}}0.
\end{align}
\paragraph{Step 2: relation to the local time at zero.} Now to prove \eqref{to_control2}, let us first observe that the law of $(\Upsilon^{\alpha,n}_{T})_{n \ge 0}$ is tight since for any $A > 0$ we can write
\begin{align*}
    \Pf_{\mathbf{T}}(\Upsilon^{\alpha,n}_{T}\ge A) = \Pf_{\mathbf{T}}(\Lambda^{}_{\ent{nT}}\ge An^{\frac{1}{\alpha}}) &\le \Pf_{\mathbf{T}}\left(\somme{m=1}{An^{\frac{1}{\alpha}}-1}{\# \mathcal{C}^i} \le nT\right).
\end{align*}
Using Theorem~\ref{thm:ConvergenceOfSize} and standard results on sums of random variables in the domain of attraction of a stable law (see for example \cite[Chapter 8.3]{bingham1989regular}), we see that this latter probability converges to $0$ as $A \to +\infty$, uniformly in $n$.
We now write
\begin{align}\label{dominate_distance_local_time}
     \sup_{0\le t \le T}|\mathbf{W} \Upsilon^{\alpha,n}_t - \Upsilon^{\alpha,\uparrow n^{\varepsilon}}_t|  &\le  \displaystyle  \sup_{0 \le t \le T} \bigg|\mathbf{W}\Upsilon^{\alpha,n}_t - \frac{ \displaystyle \somme{m=1}{\ent{n^{\frac{1}{\alpha}}\Upsilon^{\alpha,n}_t}}{Y^m_{n^{\varepsilon}}}}{n^{\frac{1}{\alpha}}}\bigg| \\&+ \displaystyle  \sup_{0 \le t \le T} \bigg|\frac{ \displaystyle \somme{m=1}{\ent{n^{\frac{1}{\alpha}}\Upsilon^{\alpha,n}_t}}{Y^m_{n^{\varepsilon}}}}{n^{\frac{1}{\alpha}}} - \Upsilon^{\alpha,\uparrow n^{\varepsilon}}_t\bigg|.\nonumber
\end{align}

\noindent We start with the second term in \eqref{dominate_distance_local_time}. On $\mathcal{D}_n \cap \mathcal{E}_n$, using \eqref{connect_local_time}, we have the bound
\begin{align}\label{eqn:local time second decomp}
     \displaystyle  \sup_{0 \le t \le T} \bigg|\frac{ \displaystyle \somme{m=1}{\ent{n^{\frac{1}{\alpha}}\Upsilon^{\alpha,n}_t}}{Y^m_{n^{\varepsilon}}}}{n^{\frac{1}{\alpha}}} - \Upsilon^{\alpha,\uparrow n^{\varepsilon}}_t\bigg| \le \displaystyle  \sup_{0 \le t \le T} \bigg| \Upsilon^{\alpha,\uparrow n^{\varepsilon}}_{t+n^{-\varepsilon/2}}- \Upsilon^{\alpha,\uparrow n^{\varepsilon}}_t\bigg| + \sup_{0\le t \le T}\bigg(\frac{Y^{\Lambda_{\ent{nt}}}_{n^{\varepsilon}}}{n^{\frac{1}{\alpha}}}\bigg).
\end{align}
By \eqref{distrib_Cv_to_use}, it is easy to verify that the first term on the right hand side of \eqref{eqn:local time second decomp} goes to zero in $\mathbb{P}_{\bT}$-probability, as $n \to \infty$.% we have
%\begin{align*}
%    \sup_{0 \le t \le T} \bigg| \Upsilon^{\alpha,\uparrow n^{\varepsilon}}_{t+n^{-\varepsilon/2}}- \Upsilon^{\alpha,\uparrow n^{\varepsilon}}_t\bigg|\underset{n \to +\infty}{\overset{(\mathbb{P}_{\bT})}{\rightarrow}}0.
%\end{align*}
 Moreover, for any $\delta,R > 0$, we can write
\begin{align*}
    \Pf_{\mathbf{T}}\left(\sup_{0\le t \le T}\bigg(\frac{Y^{\Lambda_{\ent{nt}}}_{n^{\varepsilon}}}{n^{\frac{1}{\alpha}}}\bigg) \ge \delta \right) &\le \Pf_{\mathbf{T}}\left(\sup_{0\le m \le Rn^{\frac{1}{\alpha}}}\bigg(\frac{Y^{m}_{n^{\varepsilon}}}{n^{\frac{1}{\alpha}}}\bigg) \ge \delta \right) + \Pf_{\mathbf{T}}\left( \Lambda_{\ent{nT}} \ge Rn^{\frac{1}{\alpha}}\right)\\&\le \Pf_{\mathbf{T}}\left( \sup_{0 \le t \le R} \bigg[\mathbf{W}t - \frac{ \somme{m=1}{\ent{n^{\frac{1}{\alpha}}t}}{Y^{m}_{n^{\varepsilon}}}}{n^{\frac{1}{\alpha}}}\bigg] \ge \delta/3 \right)+ \Pf_{\mathbf{T}}\left( \Lambda_{\ent{nT}} \ge Rn^{\frac{1}{\alpha}}\right).
\end{align*}
First letting $n \to +\infty$ and then letting $R \to +\infty$, the left term goes to $0$ as $n \to +\infty$ by \eqref{Local_time_0_to_nepsilon} and the right term as well by tightness of $(\Upsilon^{\alpha,n}_{T})_{n \ge 0}$. This shows that the second term on the right hand side of \eqref{eqn:local time second decomp} also goes to zero in $\mathbb{P}_{\bT}$-probability. %gives 
%\begin{align*}
%    \sup_{0\le t \le T}\bigg(\frac{Y^{\Lambda_{\ent{nt}}}_{n^{\varepsilon}}}{n^{\frac{1}{\alpha}}}\bigg) \underset{n \to +\infty}{\overset{(\mathbb{P}_{\mathbf{T}})}{\rightarrow}}0.
%\end{align*}
We deduce that the second term on the right hand side of \eqref{dominate_distance_local_time} goes to zero in $\mathbb{P}_{\bT}$-probability. The same is true for the first term using
%\begin{align*}
%    \displaystyle  \sup_{0 \le t \le T} \bigg|\frac{ \displaystyle \somme{m=1}{\ent{n^{\frac{1}{\alpha}}\Upsilon^{\alpha,n}_t}}{Y^m_{n^{\varepsilon}}}}{n^{\frac{1}{\alpha}}} - \Upsilon^{\alpha,\uparrow n^{\varepsilon}}_t\bigg|\underset{n \to +\infty}{\overset{(\mathbb{P}_{\mathbf{T}})}{\rightarrow}}0.
%\end{align*}
%For the first term in \eqref{dominate_distance_local_time}, combining
 the fact that $(\Upsilon^{\alpha,n}_{T})_{n \ge 0}$ is tight and \eqref{Local_time_0_to_nepsilon}.% we find 
%\begin{align*}
%    \sup_{0 \le t \le T} \bigg|\mathbf{W}\Upsilon^{\alpha,n}_t - \frac{ \displaystyle \somme{m=1}{\ent{n^{\frac{1}{\alpha}}\Upsilon^{\alpha,n}_t}}{Y^m_{n^{\varepsilon}}}}{n^{\frac{1}{\alpha}}}\bigg| \underset{n \to +\infty}{\overset{(\mathbb{P}_{\bT})}{\rightarrow}}0.
%\end{align*}

%Putting together these two last equations and \eqref{dominate_distance_local_time}, 
This establishes \eqref{dominate_distance_local_time} and thus concludes the proof of \eqref{to_control2}. Combining \eqref{convergence_trajectory}, \eqref{to_control2} and \eqref{distrib_Cv_to_use}, we conclude the desired result.
\end{proof}

By Proposition \ref{prop:GHP from contour}, Theorem \ref{main_theorem} is a direct consequence of Theorem \ref{convergence_contour_to_reflected_BM} (taking $\gamma = c_2$).

\section{Conditioning on the size}\label{section_conditioning_on_the_size}
 This section is dedicated to proving the following theorem, i.e. the first part of Theorem~\ref{theorem_exactconditioning}. Let us recall some notation. Conditionally on $\mathbf{T}$, for any $n \ge 0$, let us denote by $\mathcal{C}_{=n}$ (resp. $\mathcal{C}_{\ge n}$) the cluster $\mathcal{C}$ conditioned to have total size $n$ (resp. at least $n$) under $\PfT$. We denote by $\mathcal{T}_{\alpha}^{=1}$ the stable tree with parameter $\alpha$ of total mass $1$.
\begin{theorem}\label{theorem_exactconditioningsize}
Take $\gamma$ as in \eqref{eqn:annealed GHP}. Then, for $\bPb$-almost every $\bT$, the following convergence holds in law under $\mathbb{P}_{\bT}$:
\begin{align*}
(\mathcal{C}_{=n}, \gamma n^{-\left(1-\frac{1}{\alpha}\right)}d_n, n^{-1}\nu_n,\rho_n) &\underset{n \to +\infty}{\overset{(d)}{\longrightarrow}} (\mathcal{T}_{\alpha}^{=1},d_{\mathcal{\mathcal{T}_{\alpha}}},\nu_{\alpha},\rho_{\alpha}).
\end{align*}
with respect to the pointed Gromov-Hausdorff-Prokhorov topology.
\end{theorem}

The proof of this theorem is written in full detail. The argument to prove the analogous statement conditionally on the exact height will be given later in Section \ref{section_conditioning_on_the_height}. Since the latter argument is very similar, we will not give all of the details in Section \ref{section_conditioning_on_the_height}, and only explain the parts that are different. The reader may therefore like to keep in mind while reading that most of the estimates of Section \ref{section_conditioning_on_the_size} adapt straightforwardly under the exact height conditioning.

\subsection{Proof strategy}\label{sctn:strategy exact cond}

As in the proof of Theorem \ref{main_theorem}, we want to compare the law of $\mathcal{C}_{=n}$ under the quenched law $\PfT$ with the law of $\mathcal{C}_{=n}$ under the annealed law $\Pf_{\alpha}$. However, the contour function approach of the previous section fails in this setting, because the first cluster of size exactly $n$ is not captured on the timescale $[0, Tn]$. 

\begin{figure}[H]
\centering
\includegraphics[scale=0.35]{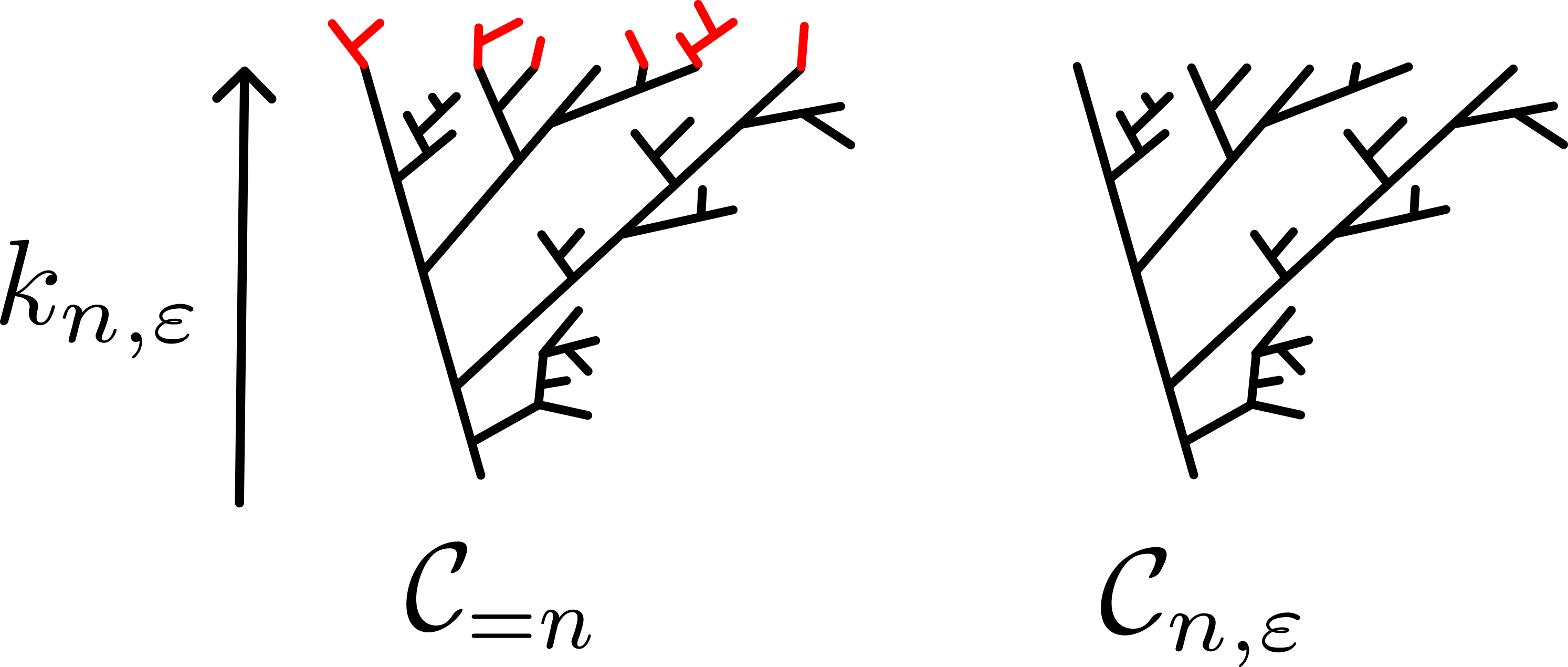}
\caption{On the left: the tree $\mathcal{C}_{=n}$ where we represent in red the part of the tree above level $k_{n,\varepsilon}$. On the right: the tree $\mathcal{C}_{n,\varepsilon}$ obtained from $\mathcal{C}_{=n}$ by cutting the red part.}
\label{tree_cut_levelnepsilon}
\end{figure}

Instead we will upgrade the result of Theorem \ref{main_theorem} via the following principal observation: the law of $\mathcal{C}_{=n}$ under $\PfT$ is the same as the law of $\Cge{(1- \varepsilon)n}$ conditioned to have size $n$. Thus we can first sample $\Cge{(1-\varepsilon)n}$ and then choose $k_{n, \eps}$ to be the minimal $k$ such that the first $k$ levels of $\Cge{(1-\varepsilon)n}$ have total mass at least $(1-\varepsilon)n$ (see Figure~\ref{tree_cut_levelnepsilon}). We denote the ball of radius $k_{n, \eps}$ in $\Cge{(1-\varepsilon)n}$ by $\mathcal{C}_{n,\varepsilon}$. It is natural to expect that under the extra conditioning and appropriate rescaling, we have $\mathcal{C}_{n,\varepsilon} \overset{d_{\mathrm{GHP}}}{\approx} \mathcal{C}_{=n}$, and moreover that, conditionally on $\C_{n, \eps}$, the probability of having exactly $n$ vertices is essentially determined by the number of individuals in generation $k_{n, \eps}$.

 Then, we want to use the following series of approximations:
\begin{align*}
\mathbb{E}_{\mathbf{T}}\bigg[F(\mathcal{C}) \text{ | }\# \mathcal{C} = n \bigg] &\approx\mathbb{E}_{\mathbf{T}}\bigg[F(\mathcal{C}_{n,\varepsilon}) \text{ | }\# \mathcal{C} = n \bigg]\\&= \frac{\mathbb{E}_{\mathbf{T}}\bigg[F(\mathcal{C}_{n,\varepsilon}) ,\# \mathcal{C} = n \text{ | }\# \mathcal{C} \ge (1-\varepsilon)n\bigg]}{\Pf_{\mathbf{T}}(\#\mathcal{C} = n\text{ | }\# \mathcal{C} \ge (1-\varepsilon)n)}\\&= \frac{\somme{i=0}{k}{\mathbb{E}_{\mathbf{T}}\bigg[F(\mathcal{C}_{n,\varepsilon}) ,\mathcal{C}_{n,\varepsilon} \in A_i,\# \mathcal{C} = n  \text{ | }\# \mathcal{C} \ge (1-\varepsilon)n \bigg]}}{\somme{i=0}{k}{\mathbb{P}_{\mathbf{T}}\bigg(\mathcal{C}_{n,\varepsilon} \in A_i,\# \mathcal{C} = n \text{ | }\# \mathcal{C} \ge (1-\varepsilon)n \bigg)}}.
\end{align*}
where $(A_i)_{0 \le i \le k}$ is a partition of the space of measured compact metric spaces $\mathbb{K}_c$ such that
\begin{itemize}
	\item[$\bullet$] For any $1 \le i \le k$, under the event $\mathcal{C}_{n,\varepsilon} \in A_i$, the value of $F(\mathcal{C}_{n,\varepsilon})$ is roughly constant, as are the size of the last generation $k_{n,\varepsilon}$ and the total size of $\mathcal{C}_{n,\varepsilon}$.
	\item[$\bullet$] We have $\Pf(\mathcal{C}_{n,\varepsilon} \in A_0 \text{ | }\# \mathcal{C} \ge (1-\varepsilon)n) \ll \varepsilon$.
\end{itemize}
The existence of such a partition will follow by combining the result of Theorem \ref{main_theorem with generation size} with a result of \cite{archer2023quenched} to control the final generation size.

Then, writing $F(A_i)$ for the value taken by $F(\mathcal{C}_{n,\varepsilon})$ when $\mathcal{C}_{n,\varepsilon} \in A_i$ and neglecting the term for $i=0$, we find that the last term of the previous equation is approximately 
\begin{align} \label{decomposing_on_Ai}
\frac{\somme{i=1}{k}{F(A_i)\Pf_{\mathbf{T}}\bigg(\mathcal{C}_{n,\varepsilon} \in A_i \text{ | }\# \mathcal{C} \ge (1-\varepsilon)n \bigg)\mathbb{P}_{\mathbf{T}}\bigg(\# \mathcal{C} = n  \text{ | }\# \mathcal{C} \ge (1-\varepsilon)n ,\mathcal{C}_{n,\varepsilon} \in A_i\bigg)}}{\somme{i=1}{k}{\Pf_{\mathbf{T}}\bigg(\mathcal{C}_{n,\varepsilon} \in A_i \text{ | }\# \mathcal{C} \ge (1-\varepsilon)n \bigg)\mathbb{P}_{\mathbf{T}}\bigg(\# \mathcal{C} = n  \text{ | }\# \mathcal{C} \ge (1-\varepsilon)n ,\mathcal{C}_{n,\varepsilon} \in A_i\bigg)}}.
\end{align}

Using Theorem \ref{main_theorem} we have
\begin{align*}
\forall i \in \{1,\cdots,k\}, \hspace*{0.2cm}\Pf_{\mathbf{T}}\bigg(\mathcal{C}_{n,\varepsilon} \in A_i \text{ | }\# \mathcal{C} \ge (1-\varepsilon)n\bigg) \sim \Pf_{\alpha}\bigg(\mathcal{C}_{n,\varepsilon} \in A_i \text{ | }\# \mathcal{C} \ge (1-\varepsilon)n\bigg).
\end{align*}
Now, it remains to prove 
\begin{align}\label{estimate_to_obtain}
\mathbb{P}_{\mathbf{T}}\bigg(\# \mathcal{C} = n  \text{ | }\# \mathcal{C} \ge (1-\varepsilon)n ,\mathcal{C}_{n,\varepsilon} \in A_i\bigg) \approx \mathbb{P}_{\alpha}\bigg(\# \mathcal{C} = n  \text{ | }\# \mathcal{C} \ge (1-\varepsilon)n ,\mathcal{C}_{n,\varepsilon} \in A_i\bigg).
\end{align}
This last approximation holds using a similar argument as the one used to prove Theorem~\ref{main_theorem}: the left-hand side $\mathbb{P}_{\mathbf{T}}\bigg(\# \mathcal{C} = n  \text{ | }\# \mathcal{C} \ge (1-\varepsilon)n ,\mathcal{C}_{n,\varepsilon} \in A_i\bigg)$ only depends on the part of the tree $\mathcal{C}$ above level $k_{n,\varepsilon}$. This should be very close to its expectation under $\mathbf{P}_{\alpha}$ by taking two independent clusters and using the fact that with very high probability the two clusters only intersect close to the origin and then evolve independently. Since the final generation size of $\C_{n, \eps}$ is roughly constant on $A_i$, this expectation is also very close to the annealed quantity (note this is not a priori automatic for conditional probabilities). This allows us to conclude that \eqref{decomposing_on_Ai} is well approximated by 
\begin{align} \label{decomposing_on_Ai_2}
\frac{\somme{i=1}{k}{F(A_i)\Pf_{\alpha}\bigg(\mathcal{C}_{n,\varepsilon} \in A_i \text{ | }\# \mathcal{C} \ge (1-\varepsilon)n \bigg)\mathbb{P}_{\alpha}\bigg(\# \mathcal{C} = n  \text{ | }\# \mathcal{C} \ge (1-\varepsilon)n ,\mathcal{C}_{n,\varepsilon} \in A_i\bigg)}}{\somme{i=1}{k}{\Pf_{\alpha}\bigg(\mathcal{C}_{n,\varepsilon} \in A_i \text{ | }\# \mathcal{C} \ge (1-\varepsilon)n \bigg)\mathbb{P}_{\alpha}\bigg(\# \mathcal{C} = n  \text{ | }\# \mathcal{C} \ge (1-\varepsilon)n ,\mathcal{C}_{n,\varepsilon} \in A_i\bigg)}},
\end{align}
which is in turn a good approximation of
\begin{align}\label{eqn:Ceps n good approx}
\mathbb{E}_{\alpha}\bigg[F(\mathcal{C}_{n,\varepsilon}) \text{ | }\# \mathcal{C} = n \bigg]\approx \mathbb{E}_{\alpha}\bigg[F(\mathcal{C}) \text{ | }\# \mathcal{C} = n \bigg].
\end{align}

The main technical inputs to run this argument are summarised in the following two key propositions, which will be proved in later subsections. The first of these constructs the events $(A_i)_{i=0}^N$ as outlined above.

\begin{proposition}\label{events_construction}
Fix $\alpha, \eps, \eta, \delta>0$. There exist $K < \infty$ depending only on $\alpha, \eps$ and $\eta$, $N< \infty$ (possibly depending on all of $\alpha, \eps, \eta, \delta$) and a sequence of sets $(A_i)_{i=0}^N$ with $A_i \subset \mathbb{K}_c \times \R_{\geq 0}$ for all $i$ such that, for all $n \geq 1$,
\begin{enumerate}[(a)]
\item $\Pt{\left((\mathcal{C}_{n,\eps}, \gamma n^{-\left(1-\frac{1}{\alpha}\right)}d_n, n^{-1}\nu_n,\rho_n), \frac{Y_{\lfloor k_{n, \eps}\rfloor}}{\gamma n^{1/\alpha}} \right) \in A_0} < \eta$ and \\ $\pr{\left((\mathcal{C}_{n,\eps}, \gamma n^{-\left(1-\frac{1}{\alpha}\right)}d_n, n^{-1}\nu_n,\rho_n), \frac{Y_{\lfloor k_{n, \eps} \rfloor}}{\gamma n^{1/\alpha}} \right) \in A_0} < \eta$.
\item $\mathbb{K}_c \times (K^{-1},K)^c \subset A_0$. 
\item For all $i=1, \ldots, N$, $\mathrm{diam}_D (A_i) \leq \delta$.
\item $(A_i)_{i=0}^N$ form a partition of $\mathbb{K}_c \times \R$.
\item For each $i \ge 1,$ we have $\Pt{\mathcal{C}_{n,\varepsilon} \in A_i \mid \#\mathcal{C}\ge (1-\varepsilon)n} \sim \mathbb{P}_{\alpha}(\mathcal{C}_{n,\varepsilon} \in A_i \mid \#\mathcal{C}\ge (1-\varepsilon)n) $.
%\item \red{ $\min_{i=1, \ldots, N} \Pt{B_{k_{\eps}}(\Ta^{\geq 1}) \in A_i}>0$?}\blue{Tanguy : I dont know if we really need that. I think we only use point $e$. Maybe that to prove $e)$ we need $f)$ ?}
\end{enumerate}
Moreover the above points also hold for the cluster conditioned on having height at least $(1-\eps)n$ (rescaled as in \eqref{eqn:joint conv height}). In this case the points (a) and (e) become (the other points do not change):
\begin{enumerate}
\item[($\widetilde{a}$)] $\Pt{\left((\widetilde{\mathcal{C}}_{n,\eps}, n^{-1}d_n, (\gamma n)^{-\frac{\alpha}{\alpha-1}}\nu_n,\rho_n), \frac{Y_{\lfloor (1-\eps)n \rfloor}}{\gamma^{\frac{\alpha}{\alpha-1}} n^{1/(\alpha-1)}} \right) \in A_0} < \eta$ and \\ $\pr{\left((\widetilde{\mathcal{C}}_{n,\eps}, n^{-1}d_n, (\gamma n)^{-\frac{\alpha}{\alpha-1}}\nu_n,\rho_n), \frac{Y_{\lfloor (1-\eps)n \rfloor}}{\gamma^{\frac{\alpha}{\alpha-1}} n^{1/(\alpha-1)}} \right) \in A_0} < \eta$.
\item[($\widetilde{e}$)] For each $i \ge 1,$ we have $\Pt{\mathcal{C}_{n,\varepsilon} \in A_i \mid \Height(\C) \ge (1-\varepsilon)n} \sim \mathbb{P}_{\alpha}(\mathcal{C}_{n,\varepsilon} \in A_i \mid  \Height(\C) \ge (1-\varepsilon)n) $.
%\item[($\widetilde{f}$)] \red{ $\min_{i=1, \ldots, N} \Pt{B_{1-\eps}(\Ta^{H \geq 1}) \in A_i}>0$?}
\end{enumerate}
\end{proposition}

The second proposition justifies the different approximations used in the proof. Before stating it, we clarify some notation. 

Throughout this section, we will fix $\eps, \eta, \delta>0$, the constants $K = K_{\alpha,\varepsilon,\eta}$, $N = N_{\alpha,\eps,\eta,\delta} $ and the family of events $(A_i)_{0 \le i \le N}$ as in Proposition~\ref{events_construction}. The reader should have in mind that the constants $\eps, \eta, \delta$ respect the ordering $0 < \delta \ll \eta \ll \varepsilon$ and will in the end be taken to zero in this order. In order to keep track of the relationships between these parameters, we will use big-O and little-o notation and we will add a subscript of $\delta, \eta$ or $\eps$ to indicate that the relevant multiplicative constants depend on these parameters. The asymptotic in the big-O always holds as $n \to \infty$ (the other parameters are viewed as fixed) but rate of convergence is allowed to depend on all three of the other parameters. Moreover, everything is allowed to depend on $\bT$. For example, $f(n) = \mathcal{O}_{\eps}(\eta)$ means that for any realisation of $\bT$, there exists a finite constant $C_{\eps}$, that depends on $\eps$ but \textit{not} on $\eta$ and \textit{not} on $\delta$, and moreover a natural number $N_{\delta, \eta, \eps}$, which may depend on all of $\delta, \eta$ and $\eps$, such that $f(n) \leq C_{\eps} \eta$ for all $n \geq N_{\delta, \eta, \eps}$. 

The second key proposition is as follows.

\begin{proposition}\label{estimateA_i}
Fix $\eps, \eta, \delta>0$, the constants $K = K_{\alpha,\varepsilon,\eta}$, $N = N_{\alpha,\eps,\eta,\delta} $ and the family of events $(A_i)_{0 \le i \le N}$ as in Proposition~\ref{events_construction}. Also let $F$ be a non-negative bounded Lipschitz function $\mathbb{K}_c \to \R$. 
\begin{enumerate}[(a)]
\item
$\mathbf{P}_{\alpha}$-almost surely,
 \begin{align*}
    \Pt{\mathcal{C}_{n,\varepsilon} \in A_0, \# \mathcal{C} = n \mid \#\mathcal{C}\ge (1-\varepsilon)n} = \mathcal{O}_{\varepsilon}(\eta n^{-1}).
 \end{align*}
 Similarly $\mathbb{P}_{\alpha} (\mathcal{C}_{n,\varepsilon} \in A_0,\#\mathcal{C}=n  \mid \# \mathcal{C} \ge (1-\varepsilon)n) = \mathcal{O}_{\varepsilon}(\eta n^{-1})$.
 \item $\mathbf{P}_{\alpha}$-almost surely, for all $i \in \{1,\dots,N\}$,
\begin{align*}
&\Pt{\mathcal{C}_{n,\varepsilon} \in A_i, \# \mathcal{C} = n \mid \#\mathcal{C}\ge (1-\varepsilon)n} \\
&\qquad = \Pf_{\alpha}(\mathcal{C}_{n,\varepsilon} \in A_i, \#\mathcal{C} = n \mid  \#\mathcal{C} \ge (1-\varepsilon)n)(1+\mathcal{O}_{\varepsilon,\eta}(\delta)) + o_{\eps}(n^{-1}).
\end{align*}
 \item $\mathbf{P}_{\alpha}$-almost surely,
 \begin{align*}
\mathbb{E}_{\mathbf{T}}\bigg[F(\mathcal{C}) \mid \# \mathcal{C} = n \bigg] =\mathbb{E}_{\mathbf{T}}\bigg[F(\mathcal{C}_{n,\varepsilon}) \mid \# \mathcal{C} = n \bigg] + \mathcal{O}(\varepsilon^{\frac{1}{2}(1-\frac{1}{\alpha})})+\mathcal{O}(\varepsilon) .
\end{align*}
\end{enumerate}
\end{proposition}

\begin{remark}
We note the following difference with the strategy of Section \ref{sctn:contour limit}. In both cases we want to prove that the relevant convergence statement holds for all non-negative bounded Lipschitz functions. In Section \ref{sctn:contour limit} this was achieved by proving convergence for a single such function and extending to all such functions using various countability and approximation arguments. The Lipschitz property was not strictly necessary for the proof; continuity would have sufficed. In this Section \ref{section_conditioning_on_the_size}, we take a different approach: we instead prove concentration of the quantities appearing in the statement of Proposition \ref{estimateA_i}. The convergence then extends essentially deterministically to all Lipschitz functions and crucially uses the Lipschitz property (as explained above).
\end{remark}

The rest of the section is organized as follows. We start by proving Theorem~\ref{theorem_exactconditioningsize} in Section~\ref{Section_proof_theorem} by following the previous proof strategy, using only Propositions \ref{events_construction} and \ref{estimateA_i} as inputs. All the required estimates are then proved in the later sections. In Section~\ref{sctn:the events Ai}, we define the events $(A_i)_{i=0}^{N}$ as described in the strategy, and prove Proposition \ref{events_construction}. In Section~\ref{Section_approximation_lemma} we make the approximations used in the strategy precise and prove Proposition \ref{estimateA_i}.

\subsection{Proof of Theorem~\ref{theorem_exactconditioningsize}, given Propositions \ref{events_construction} and \ref{estimateA_i}}
\label{Section_proof_theorem}
Fix $\eps, \eta, \delta>0$, the constants $K = K_{\alpha,\varepsilon,\eta}$, $N = N_{\alpha,\eps,\eta,\delta} $ and the family of events $(A_i)_{0 \le i \le N}$ as in Proposition~\ref{events_construction}. We follow the strategy outlined in the previous section.

\begin{proof}[Proof of Theorem~\ref{theorem_exactconditioningsize}, given Propositions \ref{events_construction} and \ref{estimateA_i}]
\label{Section_proof_theorem}
Let $F$ be a non-negative bounded Lipschitz function $\mathbb{K}_c \to \R$. Using point (c) of Proposition~\ref{estimateA_i}, we have
\begin{align*}
\mathbb{E}_{\mathbf{T}}\bigg[F(\mathcal{C}) \Big| \# \mathcal{C} = n \bigg] =\mathbb{E}_{\mathbf{T}}\bigg[F(\mathcal{C}_{n,\varepsilon}) \Big| \# \mathcal{C} = n \bigg] + \mathcal{O}(\varepsilon^{\frac{1}{2}(1-\frac{1}{\alpha})}) +\mathcal{O}(\varepsilon) .
\end{align*}
Expanding the conditional expectation, we write:
\begin{align*}
\mathbb{E}_{\mathbf{T}}\bigg[F(\mathcal{C}_{n,\varepsilon}) \mid \# \mathcal{C} = n \bigg] 
&= \frac{\mathbb{E}_{\mathbf{T}}\bigg[F(\mathcal{C}_{n,\varepsilon}) \mathbbm{1}\{\# \mathcal{C} = n\} \Big| \# \mathcal{C} \ge (1-\varepsilon)n\bigg]}{\Pf_{\mathbf{T}}(\#\mathcal{C} = n \big| \# \mathcal{C} \ge (1-\varepsilon)n)} \\
&= \frac{\sum_{i=0}^{N}\mathbb{E}_{\mathbf{T}}\bigg[F(\mathcal{C}_{n,\varepsilon})  \mathbbm{1}\{\mathcal{C}_{n,\varepsilon} \in A_i,\# \mathcal{C} = n \} \Big| \# \mathcal{C} \ge (1-\varepsilon)n \bigg]}{\sum_{i=0}^{N}\mathbb{P}_{\mathbf{T}}\bigg(\mathcal{C}_{n,\varepsilon} \in A_i,\# \mathcal{C} = n \mid \# \mathcal{C} \ge (1-\varepsilon)n \bigg)}.
\end{align*}
In both the numerator and denominator, isolating the term $i=0$ and using Proposition~\ref{estimateA_i}(a), the last equation can be rewritten as
\begin{align*}
\frac{\mathcal{O}_{\varepsilon}(\eta n^{-1})+\sum_{i=1}^{N}\mathbb{E}_{\mathbf{T}}\bigg[F(\mathcal{C}_{n,\varepsilon}) \mathbbm{1}\{\mathcal{C}_{n,\varepsilon} \in A_i,\# \mathcal{C} = n \} \Big| \# \mathcal{C} \ge (1-\varepsilon)n \bigg]}{\mathcal{O}_{\varepsilon}(\eta n^{-1}) +\sum_{i=1}^{N}\mathbb{P}_{\mathbf{T}}\bigg(\mathcal{C}_{n,\varepsilon} \in A_i,\# \mathcal{C} = n \big| \# \mathcal{C} \ge (1-\varepsilon)n \bigg)}.
\end{align*}

Using point $(c)$ in Proposition~\ref{events_construction}, for all $i \in \{1,\dots,N\}$, we have $\mathrm{diam}_D (A_i) \leq \delta$. Thus, for any $1 \le i \le N$, the random variable $F(\mathcal{C}_{n,\varepsilon})$ takes values in an interval of the form $[\alpha_i-C\delta,\alpha_i+C \delta ]$ on $\mathcal{C}_{n,\varepsilon} \in A_i$, where $\alpha_i \ge 0$ and $C > 0$ are constants depending only on $F$. Thus, we can rewrite the expression as:
\begin{align*}
\frac{\mathcal{O}_{\varepsilon}(\eta n^{-1})+(1+\mathcal{O}(\delta))\sum_{i=1}^{N}\alpha_i\mathbb{P}_{\mathbf{T}}\bigg(\mathcal{C}_{n,\varepsilon} \in A_i,\# \mathcal{C} = n  \mid \# \mathcal{C} \ge (1-\varepsilon)n \bigg)}{\mathcal{O}_{\varepsilon}(\eta n^{-1}) +\sum_{i=1}^{N}\mathbb{P}_{\mathbf{T}}\bigg(\mathcal{C}_{n,\varepsilon} \in A_i,\# \mathcal{C} = n \mid \# \mathcal{C} \ge (1-\varepsilon)n \bigg)}.
\end{align*}
Now using Proposition~\ref{estimateA_i}(b), this can be rewritten as
\begin{align*}
&\frac{\mathcal{O}_{\varepsilon}(\eta n^{-1}) + o_{\eps, \delta}(n^{-1}) +(1+\mathcal{O}_{\varepsilon,\eta}(\delta))\sum_{i=1}^{N}\alpha_i\mathbb{P}_{\alpha}\bigg(\mathcal{C}_{n,\varepsilon} \in A_i,\# \mathcal{C} = n  \mid \# \mathcal{C} \ge (1-\varepsilon)n \bigg)}{\mathcal{O}_{\varepsilon}(\eta n^{-1})+ o_{\eps, \delta}(n^{-1}) +(1+\mathcal{O}_{\varepsilon,\eta}(\delta))\sum_{i=1}^{N}\mathbb{P}_{\alpha}\bigg(\mathcal{C}_{n,\varepsilon} \in A_i,\# \mathcal{C} = n \mid \# \mathcal{C} \ge (1-\varepsilon)n \bigg)}.
\end{align*}
Using again the fact that $\mathbb{P}_{\alpha} (\mathcal{C}_{n,\varepsilon} \in A_0,\#\mathcal{C}=n  \mid \# \mathcal{C} \ge (1-\varepsilon)n) = \mathcal{O}_{\varepsilon}(\eta n^{-1})$ (see Proposition~\ref{estimateA_i}(a)),
%and setting $\alpha_0$ to be an arbitrary constant, the last quantity can be rewritten as
%\begin{align*}
%&\frac{\mathcal{O}_{\varepsilon}(\eta n^{-1})+(1+\mathcal{O}_{\varepsilon,\eta}(\delta))\sum_{i=0}^{N}\alpha_i\mathbb{P}_{\alpha}\bigg(\mathcal{C}_{n,\varepsilon} \in A_i,\# \mathcal{C} = n  \mid \# \mathcal{C} \ge (1-\varepsilon)n \bigg)}{\mathcal{O}_{\varepsilon}(\eta n^{-1})+(1+\mathcal{O}_{\varepsilon,\eta}(\delta))\sum_{i=0}^{N}\mathbb{P}_{\alpha}\bigg(\mathcal{C}_{n,\varepsilon} \in A_i,\# \mathcal{C} = n \mid \# \mathcal{C} \ge (1-\varepsilon)n \bigg)}.
%\end{align*}
and following the same logic as before, this can be rewritten as
\begin{align*}
\frac{\mathcal{O}_{\varepsilon}(\eta n^{-1})+ o_{\eps, \delta}(n^{-1})+(1+\mathcal{O}_{\varepsilon,\eta}(\delta))\mathbb{E}_{\alpha}\bigg[F(\mathcal{C}_{n,\varepsilon}) ,\# \mathcal{C} = n \mid \# \mathcal{C} \ge (1-\varepsilon)n\bigg]}{\mathcal{O}_{\varepsilon}(\eta n^{-1})+ o_{\eps, \delta}(n^{-1})+(1+\mathcal{O}_{\varepsilon,\eta}(\delta))\Pf_{\alpha}(\#\mathcal{C} = n \mid \# \mathcal{C} \ge (1-\varepsilon)n)}.
\end{align*} 
Finally, using Fact~\ref{fact:annealed conv}, the last quantity rewrites as
\begin{align*}
 &(1+\mathcal{O}_{\varepsilon,\eta}(\delta ))\mathbb{E}_{\alpha}\bigg[F(\mathcal{C}_{n,\varepsilon}) \mid \# \mathcal{C} = n \bigg] + \mathcal{O}_{\varepsilon}(\eta) + o_{\eps, \delta}(1)\\
 &\quad = (1+\mathcal{O}_{\varepsilon,\eta}(\delta ))\mathbb{E}_{\alpha}\bigg[F(\mathcal{C}) \mid \# \mathcal{C} = n \bigg] + \mathcal{O}_{\varepsilon}(\eta) + o_{\eps, \delta}(1) + f(\varepsilon),
\end{align*}
where $|f(\varepsilon)| \downarrow 0$ as $\eps \downarrow 0$.
This leads to the bound:
\begin{align*}
\bigg | \mathbb{E}_{\mathbf{T}}\bigg[F(\mathcal{C}) \mid \# \mathcal{C} = n \bigg] - \mathbb{E}_{\alpha}\bigg[F(\mathcal{C}) \mid \# \mathcal{C} = n \bigg]\bigg | = \mathcal{O}_{\varepsilon,\eta}(\delta )+ \mathcal{O}_{\varepsilon}(\eta ) + o_{\eps, \delta}(1)+ f(\varepsilon)+ \mathcal{O}(\varepsilon^{\frac{1}{2}(1-\frac{1}{\alpha})}) . 
\end{align*}
Taking first the limit $\delta \downarrow 0$, then $\eta \downarrow 0$ and finally $\varepsilon \downarrow 0$ concludes the proof, since it is already known that $ \mathbb{E}_{\alpha}\left[F(\mathcal{C}) \mid \# \mathcal{C} = n \right] \to \mathbb{E}_{\alpha}\left[ F(\mathcal{T}_{\alpha}^{=1}) \right]$.
\end{proof}

\subsection{Proof of Proposition~\ref{events_construction}: constructing the family of events $(A_i)$}\label{sctn:the events Ai}

In order to define the sets $(A_i)_{i=0}^N$, it is useful to enhance the convergence of Theorem \ref{main_theorem} to a slightly stronger topology that includes the size of generation $k_{n, \eps}$. This is useful because the size of this generation determines the conditional probability of $\{\# \C = n\}$.

To this end we work on the space $\mathbb{K}_c \times \R_{\geq 0}$, endowed with the metric $D$ and associated topology, as introduced in Section \ref{sctn:GHP topology}.

We recall that, given $m \geq 0$, $Y_m$ denotes the number of vertices in generation $m$ of $\C$ (see Section \ref{sctn:background quenched perco}). Similarly, given $t \geq 0$, $\ell_t$ denotes the local time at level $t$ in $\Ta$ - this informally corresponds to the size of generation $t$ (see Section \ref{sctn:cont trees}). Moreover, given $\Ta$, let $k_{\eps}(\Ta) = \inf \{r \geq 0: \nu_{\alpha} (B(\rho_{\alpha}, r)) \geq 1-\eps\}$. We will usually denote this just by $k_{\eps}$ when this is unambiguous. 

We recall the definitions of $\mathcal{C}_{n,\eps}$ and $k_{n,\eps}$. First, we sample $\Cge{(1-\varepsilon)n}$. Then $k_{n,\eps}$ is defined as the minimal $k$ such that the closed ball of radius $k$ of $\Cge{(1-\varepsilon)n}$ has volume at least $(1-\eps)n$. The tree $\mathcal{C}_{n,\varepsilon}$ is then obtained by cutting $\Cge{(1-\varepsilon)n}$ at level $k_{n,\eps}$ (i.e. removing everything strictly above level $k_{n, \eps}$). Similarly, in the case of conditioning on the height, we let $\C_{H, n, \eps}$ denote $\Chgen$ cut above level $\lfloor (1-\eps) n \rfloor$.

In this subsection we prove the following stronger version of Theorem~\ref{main_theorem}.

\begin{theorem}\label{main_theorem with generation size}
Take $\gamma$ as in \eqref{eqn:annealed GHP} and fix some $t \geq 1$. Then, for almost every $\eps \in [0,1]$, for $\bPb$-almost every $\bT$, the following convergences hold in law under $\mathbb{P}_{\bT}$:
\begin{align}
\left(({\C}_{H,n,\eps}, n^{-1}d_n,{(\gamma n)^{-\frac{\alpha}{\alpha - 1}}} \nu_n,\rho_n) , \frac{Y_{\lfloor nt \rfloor}}{\gamma^{\frac{\alpha}{\alpha-1}}n^{\frac{1}{\alpha-1}}} \right) &\underset{n \to +\infty}{\overset{(d)}{\longrightarrow}} \left((B_{(1-\eps)n}(\mathcal{T}^{H\geq 1-\eps}_{\alpha}),d_{\mathcal{\mathcal{T}_{\alpha}}},\nu_{\alpha},\rho_{\alpha}), \ell_{t} \right) \label{eqn:joint conv height} \\
\left((\mathcal{C}_{n, \eps}, \gamma n^{-\left(1-\frac{1}{\alpha}\right)}d_n, n^{-1}\nu_n,\rho_n), \frac{Y_{\lfloor k_{n, \eps}\rfloor}}{\gamma n^{1/\alpha}} \right) &\underset{n \to +\infty}{\overset{(d)}{\longrightarrow}} \left((B_{k_{\eps}}(\mathcal{T}^{\geq 1-\eps}_{\alpha}),d_{\mathcal{\mathcal{T}_{\alpha}}},\nu_{\alpha},\rho_{\alpha}), \ell_{k_{\eps}}\right) \label{eqn:joint conv size}
\end{align}
with respect to the product topology (pointed Gromov-Hausdorff-Prokhorov times Borel).
\end{theorem}

\begin{remark}
Theorem \ref{main_theorem with generation size} should also be true for $t < 1$, but this is more delicate to prove (and not necessary for our argument). In addition we prove \eqref{eqn:joint conv height} for all $\eps \in [0,1]$, rather than only for almost every $\eps$. The restriction to almost every $\eps$ for \eqref{eqn:joint conv size} comes from an application of Fubini's theorem. This is sufficient for our purposes but the statement should be in fact true for all $\eps \in [0,1]$.
\end{remark}

The key to proving the stronger version is the result of \cite[Theorem 1.4]{archer2023quenched} which says that the rescaled sequence of generation sizes of $\C$ rescales to a continuous state branching process (CSBP) when conditioning the height or size of $\C$ to be large. In the finite variance case, the result of Theorem \ref{main_theorem with generation size} is essentially immediate: the limiting CSBP is almost surely continuous, which ensures that both $Y_n$ and $\ell_t$ can be approximated by rescaling the measure of an approximating annulus and hence the result follows from the convergence of $\nu_n$ to $\nu$. In the stable case the limiting CSBP has positive jumps, however the argument can be saved provided we justify that the time $k_{\eps}$ is almost surely a continuity point of the limiting CSBP (this is already known for $t=1-\eps$ but we have to be careful in the random case since a priori $Y_{k_{\eps, n}}$ is essentially a size-biased generation size).

We will prove Theorem \ref{main_theorem with generation size}\eqref{eqn:joint conv height} towards the end of this section via a sequence of lemmas at the end of this subsection: specifically combining certain known results for $\Ta^{H \geq 1}$ (Fact \ref{fact:csbp}) with the key input from \cite{archer2023quenched} (Corollary \ref{cor:Yn annulus approx}). This already contains the crux of the argument and the extension to prove \eqref{eqn:joint conv size} requires a careful justification of the fact that the time $k_{\eps}$ is almost surely a continuity point for the local time at level sets in $\Ta$. Since this is rather long and not especially enlightening we have postponed the proof of Theorem \ref{main_theorem with generation size}\eqref{eqn:joint conv size} to Appendix \ref{app:Ai}.

We first show why Theorem~\ref{main_theorem with generation size} follows from Proposition~\ref{events_construction}.

\begin{proof}[Proof of Proposition~\ref{events_construction}, given Theorem~\ref{main_theorem with generation size}]
We use the shorthand $\C_{n, \eps}$ in place of the space $\left((\mathcal{C}_{n, \eps}, \gamma n^{-\left(1-\frac{1}{\alpha}\right)}d_n, n^{-1}\nu_n,\rho_n), \frac{Y_{\lfloor k_{n, \eps}\rfloor}}{\gamma n^{1/\alpha}} \right)$ and similarly $B_{k_{\eps}}(\Ta^{\geq 1-\eps})$ in place of $ \left((B_{k_{\eps}}(\mathcal{T}^{\geq 1-\eps}_{\alpha}),d_{\mathcal{\mathcal{T}_{\alpha}}},\nu_{\alpha},\rho_{\alpha}), \ell_{k_{\eps}}\right)$. We will prove (a) to (e) of the proposition. 

By the tightness of Theorem \ref{main_theorem with generation size} (which is also known to hold for the annealed law) we can choose a compact set $\hat{K} \subset \mathbb{K}_c \times \R_{\geq 0}$ such that $\pr{\C_{n, \eps} \in \hat{K}} \wedge \Pt{\C_{n, \eps} \in \hat{K}} \geq 1-\eta$ for all $n \geq 1$. We set $A_0 = \hat{K}^c$ and let $(A_i)_{i=1}^N$ denote a finite $\eps$-partition of $\hat{K}$ (w.r.t the metric $D$ defined in Section \ref{sctn:GHP topology}: first take a finite $\eps$-cover $(B_i)_{i = 1}^N$, and then set $A_i = B_i \cap (\cup_{j=1}^{i-1} B_j)^c$). Note that if any set $A_i$ has $\Pt{B_{k_{\eps}}(\Ta^{\geq 1-\eps})) \in A_i}=0$, we can just remove it from $\hat{K}$. This proves (a), (c), and (d). For part (b), we note that the limiting law $\ell_{k_{\eps}}$ is non-zero almost surely by known results on the width and total volume of $\Ta^{\geq 1-\eps}$, which ensures that we can also assume that the final generation size is bounded away from $0$ outside of $A_0$. For part (e), it suffices to show that $\Pt{B_{k_{\eps}}(\Ta^{\geq 1-\eps})) \in \partial A_i}=0$ for all $i$, which is a well-known property of $\Ta^{\geq 1-\eps}$. %This equally gives us the result of part (f): note that if any set $A_i$ has $\Pt{B_{k_{\eps}}(\Ta^{\geq 1-\eps})) \in A_i}=0$, we can just remove it from $K$. At the cost of possibly replacing $\eta$ with $2\eta$, part (a) is still satisfied.

It is clear that the same arguments work when conditioning on the height rather than the total volume.
\end{proof}

\begin{remark}
\begin{enumerate}[(A)]
\item We have not been able to find a reference in the literature for the claim (used in the last line of the previous proof) that $\Pt{B_{k_{\eps}}(\Ta^{\geq 1-\eps}) \in \partial A_i}=0$. Rather than writing a new derivation, we remark that in fact it is not really necessary for the proof of our Theorem \ref{theorem_exactconditioningsize}: in the case that this fails, we can reallocate the boundaries $\partial A_i$ to other sets to define new sets $A_i^{\mathrm{q},n}, A_i^{\mathrm{an},n} \subset A_i^{\eps}$ such that $(A_i^{\mathrm{q},n})_{i=0}^N$ and $(A_i^{\mathrm{an},n})_{i=0}^N$ both still partition $\mathbb{K}_c \times \R$, and moreover
\[
\Pt{\mathcal{C}_{n,\varepsilon} \in A^{\mathrm{q},n}_i \mid \#\mathcal{C}\ge (1-\varepsilon)n} \sim \mathbb{P}_{\alpha}(\mathcal{C}_{n,\varepsilon} \in A^{\mathrm{an},n}_i \mid \#\mathcal{C}\ge (1-\varepsilon)n)
\]
for all $i$, and rewrite the argument using these sets. (This reallocation can be made as precise as we like by tossing extra coins if necessary.)
\item Take $K=K_{\alpha, \eps, \eta}$ as in part (b) above.  Let $(X_j)_{j \geq 1}$ be i.i.d. with $\pr{X_1 \geq x} \sim K_{\alpha} x^{-1/\alpha}$ as $x \to \infty$, where $K_{\alpha}$ is the constant in \eqref{eqn:annealed size prob convergence}. Suppose that $[r_{\min}, r_{\max}] \subset [K^{-1},K]$ and $|r_{\max} - r_{\min}| \leq \delta$. In particular this implies that $\frac{r_{\max}}{r_{\min}} \leq 1+K\delta$ and $\frac{r_{\min}}{r_{\max}} \geq 1-K\delta$.

 Let $g_{\alpha}$ be the density of the positive stable random variable arising as the limit of $m^{-1/\alpha}\sum_{j=1}^m X_j$ (see \cite[Section 50]{gnedenko1968limit} for details). Moreover $g_{\alpha}$ is well-known to be continuous and bounded away from zero on the interval $[\frac{\eps}{2K}, K \eps]$. It follows that there exists a $K'=K_{\alpha,\varepsilon,\eta}'$, depending only on $\eps $ and $K$ (and thus also on $\eta$), such that, for all sufficiently large $n$,
\begin{align}\label{eqn:stable density ratio}
1-K' \delta \le \sup_{r_1, r_2, r_3, r_4 \in [r_{\min}, r_{\max}]} \frac{(r_1)^{-\alpha} g_{\alpha}(\frac{\varepsilon - r_2n^{1/\alpha-1}}{r_1})+o(n)}{(r_3)^{-\alpha} g_{\alpha}(\frac{\varepsilon - r_4n^{1/\alpha-1}}{r_3})+o(n)} \le  1+K' \delta .
\end{align}
Part (b) will later be useful for the following reason. Later, in Section~\ref{Section_approximation_lemma}, we will take a realisation of $\C_{n, \eps}$ conditioned to be in the set $A_i$, and look at the conditional probability that the entire cluster has size exactly $n$. By a local limit theorem this will asymptotically behave like the expression in \eqref{eqn:stable density ratio} and therefore this implies that this probability is approximately constant on $A_i$, provided we chose $\delta>0$ sufficiently small.

\end{enumerate}
\end{remark}

\subsubsection{Proof of Theorem \ref{main_theorem with generation size}\eqref{eqn:joint conv height}}
In the rest of this subsection we will prove Theorem \ref{main_theorem with generation size}\eqref{eqn:joint conv height} via a series of lemmas. Towards  Theorem \ref{main_theorem with generation size}\eqref{eqn:joint conv size}, we just give an outline of the proof at the end of this section, and recall that the details are provided later in Appendix \ref{app:Ai}.

We start by recalling some known facts about $\Ta$.

\begin{fact}\label{fact:csbp} The following are true:
\begin{enumerate}[(i)]
\item 
Take any $t \geq 1$. Almost surely under the conditioning $\{\Height(\Ta) \geq 1\}$,
\[
\ell_t = \lim_{\eps \downarrow 0} \eps^{-1} \nu \left(B(\rho, t+\eps) \setminus B(\rho, t) \right).
\]
(See \cite[Equation (12)]{duquesne2005probabilistic} or \cite[Equation (1.29)]{duquesne2002random}.)
\item The mapping $r \mapsto \nu (B(\rho, r))$ is almost surely continuous. This follows from \eqref{eqn:local a def}.
\item Almost surely, the process $(\ell_t)_{t \geq 0}$ has no fixed discontinuities (see \cite[Lemma 3.3]{duquesne2005probabilistic}). In other words, for any $t \geq 0$, the process is continuous at $t$ almost surely.
\end{enumerate}
\end{fact}

We now restate the result of \cite[Theorem 1.4]{archer2023quenched}, recalling that $Y_m$ denotes the number of individuals in generation $m$ that are in the root cluster.  (We omit some superfluous details from the statement; the important point is that the limiting process $(Y_t)_{t \geq 0}$ appearing in Lemma \ref{lem:AV CSBP convergence} has the same law as $(\ell_{1+t})_{t \geq 0}$ conditionally on $\ell_1>0$.)

\begin{lemma}\label{lem:AV CSBP convergence}(\cite[Theorem 1.4]{archer2023quenched}.)
For $\bPb$-almost every T, under the conditioning $Y_n>0$ we have that the process $(n^{-\frac{1}{\alpha-1}}Y_{\lfloor n(t+1) \rfloor})_{t \geq 0}$ converges in distribution (under $\mathbb{P}_{\bT}$) to an $\alpha$-stable CSBP $(\widetilde{Y}_t)_{t \geq 0}$ (with branching mechanism given by \cite[Lemma A.1]{archer2023quenched}), and where $\widetilde{Y}_0$ is a random variable having Laplace transform given by \cite[Equation (1.3)]{archer2023quenched}. This convergence holds with respect to the Skorokhod-$J_1$ topology on the space $D([0, \infty), [0, \infty))$.
\end{lemma}

For the rest of this section we let $Y^{(n)}_{t}$ denote $n^{-1/(\alpha-1)}Y_{\lfloor n^{1-1/\alpha}t\rfloor}$ conditionally on $Y_n > 0$.

\begin{corollary}[$Y^{(n)}_t$ is well-approximated by its average over a small annulus]\label{cor:Yn annulus approx}
For any $\delta>0, t \geq 1$, there exists $\eps>0$ and $N < \infty$ such that, for all $n \geq N$ and all $\eps' \in (0,\eps)$:
\[
\Ptcond{n^{-1/(\alpha-1)}\left|Y_{\lfloor nt \rfloor} - (\eps' n)^{-1} \sum_{m= \lfloor tn \rfloor}^{\lfloor (t+\eps')n\rfloor} Y_{m}\right| > \delta}{Y_n > 0}{\bT} < \delta.
\]
\end{corollary}
\begin{proof}
Fix $\delta>0$. By Lemma \ref{lem:AV CSBP convergence}, we know that $(Y^{(n)}_{\lfloor t+1 \rfloor})_{t \geq 0}$ converges in law to $(\widetilde{Y}_t)_{t \geq 0}$, and that this limiting process is almost surely continuous at time $t$ (by Fact \ref{fact:csbp} and since $(\widetilde{Y}_t)_{t \geq 0}$ has the same law as $(\ell_{1+t})_{t \geq 0}$ conditionally on $\ell_1>0$). Hence we can find $\eps \in (0, \delta)$ such that $\sup_{s, s' \in [(t-\eps) \vee 0, t+\eps]} |\widetilde{Y}_s - \widetilde{Y}_{s'}| \leq \delta$ with probability at least $1-\delta$.

By the Skorokhod representation theorem (the space $D([0, \infty), [0, \infty))$ is separable by \cite[Theorem 12.2]{BillsleyConv}), we can assume that the convergence of $Y^{(n)}$ to $\widetilde{Y}$ is almost sure. In particular, we choose $N<\infty$ such that, for all $n \geq N$,
\[
\Pt{d_{J_1} \left((Y^{(n)}_{1+s})_{s \in [0,2t]}), (\widetilde{Y}_s)_{s \in [0,2t]})\right) < \frac{\eps}{2}} > 1-\delta.
\]
When both of the high probability events above occur, we have that $\sup_{s, s' \in [(t-\eps/2) \vee 0, t+\eps/2]} |Y^{(n)}_s - Y^{(n)}_{s'}| \leq \delta + \eps$ and hence that, for all $\eps' \in (0, \eps/2)$:
\[
\left|n^{-1/(\alpha-1)}Y_{\lfloor nt \rfloor} - (\eps' n)^{-1} \sum_{m= \lfloor tn \rfloor}^{\lfloor (t+\eps')n\rfloor} n^{-1/(\alpha-1)}Y_{m}\right| \leq \delta + \eps \leq 2 \delta.
\]
(This proves the claim with $2 \delta$ in place of $\delta$.)
\end{proof}

We now have all the ingredients to prove \eqref{eqn:joint conv height}. This is easier to prove than \eqref{eqn:joint conv size} since the result of \cite[Theorem 1.4]{archer2023quenched} is also stated conditionally on the height. Afterwards, we will adapt this proof to prove \eqref{eqn:joint conv size}.

\begin{proof}[Proof of Theorem \ref{main_theorem with generation size}\eqref{eqn:joint conv height}]
We appeal to the Skorokhod representation theorem and separability of $\mathbb{K}_c$ to assume that the convergence of Theorem \ref{main_theorem} holds almost surely on the space $(\Omega_{\bT}, \F_{\bT}, \mathbb{P}_{\bT})$.  By standard results relating to the GHP topology (see \cite[Theorem 4.11]{evans2006probability} for further details) we can almost surely find a sequence $\delta_n \downarrow 0$ and isometrically embed $(\Chgn)_{n \geq 1}$ and $\Ta$ into a common metric space $(M,D_M)$ so that
\[
d_H\left((\Chgn, \gamma^{-1} n^{-1}d_n, (\Ta, d_{\Ta})\right) \vee d_P((\gamma n)^{-\frac{\alpha}{\alpha - 1}}\nu_n, \nu) \vee D_M(\rho_n, \rho) \leq \delta_n.
\]
As a consequence, the volume of any small fixed annulus converges almost surely: more specifically, for any fixed $t \geq 1, s \geq 0$, (here balls are measured with respect to the common metric $D_M$), we have that $[B(\rho_n, t+s) \setminus B(\rho_n, t)] \cap \Chgn \subset [B(\rho, t+s) \setminus B(\rho, t)]^{2 \delta_n}$, so 
\begin{align*}
(\gamma n)^{-\frac{\alpha}{\alpha - 1}} \nu_n (B(\rho_n, t+s) \setminus B(\rho_n, t)) &\leq (\gamma n)^{-\frac{\alpha}{\alpha - 1}} \nu_n( [(B(\rho, t+s) \setminus B(\rho, t))]^{2\delta_n} \\
&\leq  \nu( [(B(\rho, t+s) \setminus B(\rho, t))]^{3\delta_n} + \delta_n.
\end{align*}
Taking $\delta_n \to 0$ gives $\limsup_{n \to \infty} (\gamma n)^{-\frac{\alpha}{\alpha - 1}} \nu_n (B(\rho_n, t+s) \setminus B(\rho_n, t)) \leq \nu ((B(\rho, t+s) \setminus B(\rho, t)))$, using the continuity of Fact \ref{fact:csbp}(ii). The lower bound is similar.

We thus fix $\delta>0$ and carry out the following procedure.
\begin{enumerate}
\item Choose $\eps$ small enough that $\Pt{|\ell_t -  \eps^{-1} \nu \left(B(\rho, t+\eps) \setminus B(\rho, t) \right) | \geq \delta} < \delta$. (This is possible by Fact \ref{fact:csbp}(i).)
\item Reduce $\eps>0$ if necessary and choose $N< \infty$ so that, for all $n \geq N$,
\[
\Pt{|\eps^{-1} n^{-\frac{\alpha}{\alpha - 1}} \nu_n (B(\rho_n, t+\eps) \setminus B(\rho_n, t)) - Y^{(n)}_t| < \delta } \geq 1-\delta.
\]
(This is possible by Corollary \ref{cor:Yn annulus approx}.)
\item Note that $\eps>0$ has been fixed by the previous two steps. Now increase $N<\infty$ if necessary so that
\[
\Pt{| (\gamma n)^{-\frac{\alpha}{\alpha - 1}} \nu_n (B(\rho_n, t+\eps) \setminus B(\rho_n, t)) - \nu ((B(\rho, t+\eps) \setminus B(\rho, t)))| < \eps \delta \text{ for all } n \geq N} \geq 1-\delta.
\]
(This is possible by the convergence of annuli shown above.)

\end{enumerate}
By the triangle inequality we thus have for all $n \geq N$ that $\Pt{|\gamma^{\frac{\alpha}{\alpha-1}}\ell_t - Y^{(n)}_t| < 3\delta } \geq 1-3\delta$, and thus we are done.
\end{proof}

The proof of \eqref{eqn:joint conv size} is very similar to that of \eqref{eqn:joint conv height}, except that we need to separately verify that $k_{n, \eps}$ converges to $k_{\eps}$ under rescaling, and that the latter is almost surely a continuity point of the limiting CSBP. In particular, rather than working under the conditioning $\{\#\C \geq (1-\eps)n\}$, we will work under the conditioning $\{\#\C \geq n/2, H \geq cn^{1-\frac{1}{\alpha}}\}$, and later restrict to the event $\{\#\C \geq (1-\eps)n\}$. By choosing $c>0$ sufficiently small, the event $\{\#\C \geq n/2, H \geq cn^{1-\frac{1}{\alpha}}\}$ is an arbitrarily good approximation of the event $\{\#\C \geq n/2\}$, but the extra conditioning on the height will enable us to apply the result of \eqref{eqn:joint conv height}. A priori, the time $t \geq 1$ appearing in \eqref{eqn:joint conv height} is fixed, but it extends immediately to hold for an independent $\textsf{Uniform}([1,t])$ variable in place of $t$. By randomising $k_{n, \eps}$ - specifically replacing it with $k_{n, \textsf{Uniform}{[0,1]}}$ where $c'>0$ is chosen arbitrarily small - we can show that $k_{n, \textsf{Uniform}{[0,1]}}$ is comparable to the original $\textsf{Uniform}([1,t])$ variable and so the result transfers to $k_{n, \textsf{Uniform}{[0,1]}}$. To conclude we would like to ``derandomise'' $k_{n, \textsf{Uniform}{[0,1]}}$, for which we apply Fubini's theorem, and which leads to the restriction of the result to almost every $\eps>0$.

The details of this proof are given in Appendix \ref{app:Ai}.

\subsection{Proof of Proposition~\ref{estimateA_i}: approximation lemmas}\label{Section_approximation_lemma}

%This section is dedicated to establishing the approximations employed in our proof strategy of Section \ref{sctn:strategy exact cond}. More precisely, we will: 
%\begin{itemize}
%\item justify the estimate appearing in \eqref{estimate_to_obtain}, for $i \geq 1$ (this is the content of Proposition \ref{estimateA_i}(a)),
%\item prove an analogous upper bound to \eqref{estimate_to_obtain} in the case $i=0$, which allows us to neglect the contribution from $A_0$ (content of Proposition \ref{estimateA_i}(b)),
%\item justify \eqref{eqn:Ceps n good approx}, that $\C_{n, \eps}$ is a good approximation of $\C_{=n}$ under the conditioning (content of  Proposition \ref{estimateA_i}(c)).
%\end{itemize}

Throughout this section, we fix $\eps, \eta, \delta>0$, the constants $K = K_{\alpha,\varepsilon,\eta}$, $N = N_{\alpha,\eps,\eta,\delta} $ and the family of events $(A_i)_{0 \le i \le N}$ as in Proposition~\ref{events_construction}. The reader should have in mind that the constants $\eps, \eta, \delta$ respect the ordering $0 < \delta \ll \eta \ll \varepsilon$ and will in the end be taken to zero in this order. Moreover we will be using big-O notation as outlined above the statement of Proposition \ref{estimateA_i}.

 To simplify notation, throughout this section, as in Proposition~\ref{events_construction}, the notation $\C_{n, \eps}$ will refer to $\left((\mathcal{C}_{n, \eps}, \gamma n^{-\left(1-\frac{1}{\alpha}\right)}d_n, n^{-1}\nu_n,\rho_n), \frac{Y_{\lfloor k_{n, \eps}\rfloor}}{\gamma n^{1/\alpha}} \right)$. This section is dedicated to proving Proposition~\ref{estimateA_i}.

The key to the proof of Proposition \ref{estimateA_i} will be Lemma \ref{lemma_exact_approximations}, a fairly general statement that allows us to separate events occurring before and after generation $k_{n, \eps}$ by conditioning on the size of generation $k_{n, \eps}$. Before stating it, we provide a technical lemma with some useful estimates.

To introduce these estimates we need some notation: Let $\mathcal{C}^1$ and $\mathcal{C}^2$ be two independent percolation clusters under the quenched measure $\mathbb{P}_{\mathbf{T}}$. For $j \in \{1,2\}$, we add a superscript $j$ when an event or random variable refers to $\mathcal{C}^j$. For any $n \ge 1$, we define the following events:
\begin{itemize}
    \item[$\bullet$] $\mathcal{H}_n$: the event where $\Height (\mathcal{C}^1 \cap \mathcal{C}^2) < n^{\frac{1}{2}-\frac{1}{2\alpha}}$. 
    \item[$\bullet$]  $\mathcal{S}_n$: the event where we have $k_{n,\varepsilon} \ge n^{\frac{1}{2}-\frac{1}{2\alpha}}$.
\end{itemize}

%To streamline the proof of Lemma \ref{lemma_exact_approximations}, we first prove a preliminary result containing some key estimates on the above quantities.

\begin{lemma}\label{eventHandS}
There exist constants $C,c > 0$ and $\beta > 0$ such that, $\mathbf{P}_{\alpha}$-almost surely, for all $n$ large enough:
\begin{enumerate}[(a)]
\item $  \Pt{\mathcal{H}_n^{c}} \le Ce^{-cn^{\beta}} $.
\item $\Pt{\mathcal{S}_n, \# \mathcal{C} = n}\le Ce^{-cn^{\beta}}$.
\item $\Pt{\mathcal{S}_n, \# \mathcal{C} \geq (1-\eps)n} = o(n^{-\frac{1}{\alpha}})$.
\end{enumerate}
\end{lemma}
\begin{proof}
\begin{enumerate}[(a)]
\item 
The first bound follows from the observation that $\mathcal{C}^1\cap \mathcal{C}^2$ is a percolation cluster with parameter $\mu^{-2}$. Indeed, the expected number of vertices of $\mathcal{C}^1\cap \mathcal{C}^2$ at level $k$ is given by $\mu^{-2k}\#\mathbf{T}_k \sim \mathbf{W} \mu^{-k}$, where $\#\mathbf{T}_k$ denotes the size of generation $k$ in $\mathbf{T}$. The conclusion therefore follows from a standard first moment method, along with the Borel-Cantelli lemma.
\item  We first claim that there exists a constant $\beta > 0$ such that the annealed probability satisfies:
\begin{align}\label{to_prove1}
\Pf_{\alpha}\left(k_{n,\varepsilon} < n^{\frac{1}{2}-\frac{1}{2\alpha}}, \# \mathcal{C} = n\right) \le e^{-n^{\beta}}.
\end{align}
If $k_{n,\varepsilon} < n^{\frac{1}{2}-\frac{1}{2\alpha}}$, then there exists $k \in \{0,\dots,  \lfloor n^{\frac{1}{2}-\frac{1}{2\alpha}} \rfloor - 1\}$ such that $Y_k \ge \frac{n/2}{n^{\frac{1}{2}-\frac{1}{2\alpha}}} = \frac{1}{2} n^{\frac{1}{2}+\frac{1}{2\alpha}}$. Moreover, on the event $\{\# \mathcal{C} = n\}$, the trees emanating from level $k$ are again independent Galton-Watson trees distributed as $\mathcal{C}$ under $\Pf_{\alpha}$, with a total mass strictly less than $n$. A union bound over $k$ therefore yields the crude upper bound:
\begin{align}\label{eqn:crude kneps}
\Pf_{\alpha}\left(k_{n,\varepsilon} < n^{\frac{1}{2}-\frac{1}{2\alpha}}, \# \mathcal{C} = n\right) \le n^{\frac{1}{2}-\frac{1}{2\alpha}}\left(1-\Pf_{\alpha}(\# \mathcal{C} \ge n)\right)^{\frac{1}{2} n^{\frac{1}{2}+\frac{1}{2\alpha}}}.
\end{align}
Choosing $\beta \in (0, \frac{1}{2}-\frac{1}{2\alpha})$ and recalling that $\Pf_{\alpha}(\# \mathcal{C} \ge n) \sim K_{\alpha}n^{-\frac{1}{\alpha}}$ (see \eqref{eqn:annealed size prob convergence}), standard computations show that \eqref{to_prove1} holds with constant exponent $\beta$.

Finally, to transfer this bound to the quenched case, we apply Markov's inequality. For any $\delta > 0$,
\[
\mathbf{P}_{\alpha}\left( \Pt{k_{n,\varepsilon} < n^{\frac{1}{2}-\frac{1}{2\alpha}}, \# \mathcal{C} = n} \ge \delta \right) \le \frac{1}{\delta} \Pf_{\alpha}\left(k_{n,\varepsilon} < n^{\frac{1}{2}-\frac{1}{2\alpha}}, \# \mathcal{C} = n\right).
\]
Setting $\delta_n = e^{-n^{\beta}/2}$, the Borel-Cantelli Lemma combined with the estimates from \eqref{to_prove1} imply that for $\mathbf{P}_{\alpha}$-almost every tree $\mathbf{T}$, and for all $n$ large enough:
\[
\Pt{k_{n,\varepsilon} < n^{\frac{1}{2}-\frac{1}{2\alpha}}, \# \mathcal{C} = n} < e^{-n^{{\beta}}/{2}}.
\]
This concludes the proof.
\item For any $n/2 \leq m \leq n^{1+\frac{\alpha-1}{4}}$, the same logic as in \eqref{eqn:crude kneps} leads to
\[
\Pf_{\alpha}\left(k_{n,\varepsilon} < n^{\frac{1}{2}-\frac{1}{2\alpha}}, \# \mathcal{C} = m\right) \le n^{\frac{1}{2}-\frac{1}{2\alpha}}\left(1-\Pf_{\alpha}(\# \mathcal{C} \ge m)\right)^{\frac{1}{2} n^{\frac{1}{2}+\frac{1}{2\alpha}}} \leq \exp\{-cn^{\frac{\alpha-1}{4\alpha}}\},
\]
so Borel-Cantelli gives that, $\bPb$-almost surely, for all $n$ large enough and all $n/2 \leq m \leq n^{1+\frac{\alpha-1}{4}}$, we have $\Pt{k_{n,\varepsilon} < n^{\frac{1}{2}-\frac{1}{2\alpha}}, \# \mathcal{C} = m}\leq \exp\{-cn^{\frac{\alpha-1}{8\alpha}}\}$. Combining with Theorem \ref{thm:ConvergenceOfSize} we have that
\begin{align*}
&\Pt{\mathcal{S}_n, \# \mathcal{C} \geq (1-\eps)n} \\
&\quad = \sup_{(1-\eps)n \leq m \leq  n^{1+\frac{\alpha-1}{4}}} \Pt{k_{n,\varepsilon} < n^{\frac{1}{2}-\frac{1}{2\alpha}}, \# \mathcal{C} = m} + \Pt{\# \C \geq  n^{1+\frac{\alpha-1}{4}}} = o(n^{-\frac{1}{\alpha}}).\qedhere
\end{align*}
%\item Clearly, by Markov's inequality, the $\bPb$-probability that there exists $(h,r,m)$ violating the event is upper bounded by $n^{-7}$. The result therefore follows by Borel-Cantelli.
\end{enumerate}
\end{proof}

The next lemma will be useful to separate what happens in our cluster $\mathcal{C}$ until generation $k_{n,\varepsilon}$ and after generation $k_{n,\varepsilon}$.

Given $\eps>0$ and conditionally on the event $\{\#\mathcal{C} \ge (1-\varepsilon)n\}$, we introduce for any $(h,r,m) \in \llbracket 1, n \rrbracket^3$ the event $\mathcal{F}(h,r,m)$, defined as $\{k_{n,\varepsilon} = h\} \cap \{Y_{h} = r\} \cap \{\#\mathcal{C}_{n,\varepsilon} = m\}$. In other words this event specifies the height, size and number of vertices of the last generation in $\mathcal{C}_{n,\eps}$.  We fix an event $D_n \in \mathcal{G}_{\bT}$ that, conditionally on $\mathcal{F}(h,r,m)$, is measurable with respect to the generations up to and including generation $h$ in $\mathcal{C}$. We fix an event $E_n \in \mathcal{G}_{\bT}$ that, conditionally on $\mathcal{F}(h,r,m)$, is measurable with respect to the generations strictly after generation $h$ in $\mathcal{C}$. Proposition~\ref{estimateA_i}(a) and (b) will then follow by taking $D_n = D_n^i = \{\mathcal{C}_{n,\varepsilon} \in A_i\}$ for $i \in \{0,\ldots,N\}$, and Proposition~\ref{estimateA_i}(c) will follow by taking $E_n$ to be the event on which $\Height(\mathcal{C}) \ge k_{n,\varepsilon}+{(\varepsilon^{\frac{1}{2}} n)}^{1-\frac{1}{\alpha}}$ (see Proposition~\ref{bound_proba_exceed_level_kneps}).

The following quantities will be useful in the proof. For any $(h,r,m) \in \llbracket 1, n \rrbracket^3$, let us define:
\begin{align}\label{definition_Sp}
\begin{split}
p_n(h,r,m {,D_n}) &= \Pt{D_n, \mathcal{F}(h,r,m) , \#\mathcal{C} \ge (1-\varepsilon)n}, \\
S_n(h,r,m{, E_n},D_n) &= \Pt{\#\mathcal{C} = n {,E_n} \mid D_n, \mathcal{F}(h,r,m), \#\mathcal{C} \ge (1-\varepsilon)n}  
\end{split}
\end{align}

We are now ready to state and prove our general lemma.

\begin{lemma}\label{lemma_exact_approximations}
    Let us fix two sequences of events $(D_n)_{n \ge 1}$ and $(E_n)_{n \ge 1}$ as above. Then, $\mathbf{P}_{\alpha}$-almost surely, there exists a random variable $F_{\mathbf{T}}(n, D_n, E_n)$ such that for all $n$ large enough:
    \begin{align}\label{eqn:lemma general main statement ideal}
    \begin{split}
   & \Pt{D_n, \# \mathcal{C} = n,  E_n  \mid \#\mathcal{C}\ge (1-\varepsilon)n} \\
   \quad &= [\Pt{D_n \mid \#\mathcal{C}\ge (1-\varepsilon)n} +o(1)] F_{\mathbf{T}}(n, D_n, E_n) + o(n^{-1}),
   \end{split}
    \end{align}
    where $F_{\mathbf{T}}(n, D_n, E_n)$ satisfies:
    \[
    \inf_{\substack{ (h,r,m): \\ p_n(h,r,m,D_n) > 0}} \Eb{S_n(h,r,m,  E_n,D_n)} \leq F_{\mathbf{T}}(n, D_n, E_n) \leq \sup_{\substack{ (h,r,m): \\ p_n(h,r,m,D_n) > 0}} \Eb{S_n(h,r,m,  E_n,D_n)}.
    \]
\end{lemma}
\begin{proof}
We will instead prove that
    \begin{align}\label{eqn:lemma general main statement}
    \begin{split}
   & \Pt{D_n, \# \mathcal{C} = n,  E_n,  {k_{n,\varepsilon} \ge n^{\frac{1}{2}-\frac{1}{2\alpha}}}, \#\mathcal{C}\ge (1-\varepsilon)n} \\
   &\qquad = \Pt{D_n, {k_{n,\varepsilon} \ge n^{\frac{1}{2}-\frac{1}{2\alpha}}}, \#\mathcal{C}\ge (1-\varepsilon)n} F_{\mathbf{T}}(n, D_n, E_n) + o(n^{-2}),
   \end{split}
    \end{align}
    which implies the result using the asymptotic of Theorem \ref{thm:ConvergenceOfSize} and Lemma \ref{eventHandS}(b,c).
    
    Since the events $D_n$ and $E_n$ are fixed in the whole proof, we abbreviate  $I_n := \llbracket n^{\frac{1}{2}-\frac{1}{2\alpha}} , n\rrbracket \times \llbracket 1, n \rrbracket^2$, $p_n(h,r,m) := p_n(h,r,m,D_n)  $,  $S_n(h,r,m): = S_n(h,r,m,  E_n,D_n)$ and $F_{\mathbf{T}}(n) := F_{\mathbf{T}}(n, D_n, E_n)$. We start by writing
    \begin{align}\label{starting_equation }
       \Pt{D_n, \#\mathcal{C}=n ,  E_n, {k_{n,\varepsilon} \ge n^{\frac{1}{2}-\frac{1}{2\alpha}}}} = \sum_{(h,r,m) \in I_n} p_n(h,r,m) S_n(h,r,m).
    \end{align}
    For any $(h,r,m) \in I_n$, define the event
    \begin{align*}
        \mathcal{A}_n(h,r,m) = \{p_n(h,r,m) \ge n^{-20}\}.
    \end{align*} 
%     \red{{(want to remove the second point to help with ??? below (14), and also measurability wrt first h levels of $\bT$, but need to update bullet point below)}}
%   where $C,c,\beta > 0$ are the constants given by Lemma~\ref{eventHandS}. 
Note that this event is measurable with respect to the $h$ first levels of $\mathbf{T}$. 
%Let us make a few observations regarding $\mathcal{A}_n(h,r,m)$:
%    \begin{itemize}
%    	\item[$\bullet$] For almost every $\mathbf{T}$ and for $n$ large enough, if $\mathcal{A}_n(h,r,m)$ occurs, then $h > n^{\frac{1}{2}-\frac{1}{2\alpha}}$. Indeed, if this is not the case, $S_n(h,r,m)$ would be bounded by:
%        \begin{align*}
%        \frac{\Pt{k_{n,\varepsilon} \le n^{\frac{1}{2}-\frac{1}{2\alpha}}, \# \mathcal{C} = n}}{p_n(h,r,m)} &= \frac{\Pt{k_{n,\varepsilon}^1 \le n^{\frac{1}{2}-\frac{1}{2\alpha}},k_{n,\varepsilon}^2 \le n^{\frac{1}{2}-\frac{1}{2\alpha}}, \# \mathcal{C}^1 = n,\# \mathcal{C}^2 = n}^{\frac{1}{2}}}{p_n(h,r,m)} \\
%        &\le \frac{\Pt{\mathcal{S}_n^{c}, \# \mathcal{C}^1 = n, \#\mathcal{C}^2 = n}^{\frac{1}{2}}}{p_n(h,r,m) } = o(n^{-20}),
%        \end{align*}
%        which contradicts the second inequality in the definition of $\mathcal{A}_n(h,r,m)$.
%        \item[$\bullet$] The event $\mathcal{A}_n(h,r,m)$ is measurable with respect to the first $h$ generations of $\mathbf{T}$. This is a direct consequence of the previous item. \red{not for final prob}
%        \item[$\bullet$] By Lemma~\ref{eventHandS}, for almost every $\mathbf{T}$ and for $n$ large enough, the last two inequalities defining $\mathcal{A}_n(h,r,m)$ are satisfied.
%    \end{itemize}
The aim of the proof is to show the following sequence of equalities:
\begin{align}\label{eqn:equalities key lem}
\begin{split}
\sum_{(h,r,m) \in I_n} p_n(h,r,m) S_n(h,r,m) &= \sum_{(h,r,m) \in I_n} p_n(h,r,m) S_n(h,r,m) \indi{\mathcal{A}_n(h,r,m)} + o(n^{-2}) \\
&= \sum_{(h,r,m) \in I_n} p_n(h,r,m) \Eb{S_n(h,r,m)} \indi{\mathcal{A}_n(h,r,m)} + o(n^{-2}) \\
&= \sum_{(h,r,m) \in I_n} p_n(h,r,m) \Eb{S_n(h,r,m)} + o(n^{-2}).
\end{split}
\end{align}
The first and third equalities are immediate since the error term in both cases is tautologically upper bounded by $n^{-17}$.
    We now claim that for $n$ large enough, for all $(h,r,m) \in  I_n$, under the event $\mathcal{A}_n(h,r,m)$ we have:
    \begin{align}\label{control_diff_quenched_annealed}
       | S_n(h,r,m)- \Eb{S_n(h,r,m)}| \le n^{-3}.
    \end{align}
    Note that the second equality in \eqref{eqn:equalities key lem} is then immediate. 
    
    To this end, fix $(h,r,m) \in I_n$ such that $\mathbf{P}_{\alpha}(\mathcal{A}_n(h,r,m)) > 0$. We express the conditional variance $\Varb{ S_n(h,r,m) \mid \mathcal{A}_n(h,r,m)}$ as:
    \begin{align}\label{eqn:var decomp}
        \Eb{S^{1}_n(h,r,m)S^{2}_n(h,r,m) \mid \mathcal{A}_n(h,r,m)} - \Eb{S_n(h,r,m) \mid \mathcal{A}_n(h,r,m)}^2,
    \end{align}
    where $S^{1}_n$ and $S^{2}_n$ denote independent (under $\mathbb{P}_{\bT}$) copies of $S_n$ defined via independent percolation clusters $\mathcal{C}^1$ and $\mathcal{C}^2$. Let $\mathcal{B}_{n}(h,r,m)$ be the event on which, for each $j \in \{1,2\}$, we have $\#\mathcal{C}^j \ge (1-\varepsilon)n$, $\mathcal{C}^j \in D_n$, and $\mathcal{F}^j_n(h,r,m)$ occurs (here a superscript $j$ denotes that an event or random variable refers to $\mathcal{C}^j$). Then, for any $\mathbf{T} \in \mathcal{A}_n(h,r,m)$, we have for all $n$ large enough (uniformly in $h,r,m$):
    \begin{align*}
        S^{1}_n(h,r,m)S^{2}_n(h,r,m) &= \Ptcond{\#\mathcal{C}^1 = n  {,\mathcal{C}^1 \in E_n}, \#\mathcal{C}^2 = n {,\mathcal{C}^2 \in E_n}}{\mathcal{B}_{n}(h,r,m)}{\mathbf{T}}{} \\
        &\le C'e^{-c'n^{\beta}} + \Ptcond{\#\mathcal{C}^1 = n  {,\mathcal{C}^1 \in E_n}, \#\mathcal{C}^2 = n {,\mathcal{C}^2 \in E_n}, \mathcal{H}_n}{\mathcal{B}_{n}(h,r,m)}{\mathbf{T}}{},
    \end{align*}
    for some universal constants $C',c' > 0$, and where $\mathcal{H}_n$ is as in Lemma \ref{eventHandS}. Furthermore, by conditioning on $\F_h$, observe that:
    \begin{align*}
         &\Eb{\Ptcond{\#\mathcal{C}^1 = n  {,\mathcal{C}^1 \in E_n}, \#\mathcal{C}^2 = n {,\mathcal{C}^2 \in E_n}, \mathcal{H}_n}{\mathcal{B}_{n}(h,r,m)}{\mathbf{T}}{} \bigg| \mathcal{A}_n(h,r,m)} 
         \end{align*}
is equal to the expectation (w.r.t. $\bPb$, and conditionally on $\mathcal{A}_n(h,r,m)$) of 
         \begin{align*}
\sum_{\substack{v_1^{1}, \ldots, v_h^{1}\\  v_1^{2}, \ldots, v_h^{2} \\
\text{disjoint} }}\prcond{v_i^j \in \C^j \forall i,j}{\mathcal{B}_{n}(h,r,m)}{\bT} \Eb{\prcond{\#\mathcal{C}^j = n  {,\mathcal{C}^j \in E_n} \forall j}{v_i^j \in C^j \forall i,j, \mathcal{B}_{n}(h,r,m)}{\bT}\mid \mathcal{A}_n(h,r,m)}.      
        % &= \Eb{S_n(h,r,m) \mid \mathcal{A}_n(h,r,m)}^2.
    \end{align*}
    Now note that the term $\Eb{\prcond{\#\mathcal{C}^j = n  {,\mathcal{C}^j \in E_n} \forall j}{v_i^j \in C^j \forall i,j, \mathcal{B}_{n}(h,r,m)}{\bT}\mid \mathcal{A}_n(h,r,m)}$ factorises over $\C^1$ and $\C^2$ since, under the conditioning, the events $\{\#\mathcal{C}^j = n,\mathcal{C}^j \in E_n\}$ are respectively measurable with respect to the subtrees of $\bT$ emanating from $v_1^j, \ldots, v_h^j$, which are disjoint. Moreover, both expectations are identical, and do not depend on the choice of $v_i^j$: in particular, we have for all terms in the sum that 
\begin{align*}
&\Eb{\prcond{\#\mathcal{C}^j = n  {,\mathcal{C}^j \in E_n} \forall j}{v_i^j \in C^j \forall i,j, \mathcal{B}_{n}(h,r,m)}{\bT}\mid \mathcal{A}_n(h,r,m)} \\
&= \Eb{\prcond{\#\mathcal{C}^1 = n  {,\mathcal{C}^1 \in E_n}}{\mathcal{B}^1_{n}(h,r,m)}{\bT}\mid \mathcal{A}_n(h,r,m)}^2.
\end{align*}    
Since this term does not depend on the choice of the $v_i^j$, it can be factorised outside of the sum, and indeed outside the entire conditional expectation. We are left with the sum over the choices of $v_i^j$ which is at most $1$, and hence we get that
\begin{align*}
&\Eb{S^{1}_n(h,r,m)S^{2}_n(h,r,m)} \\
%&\leq C'e^{-c'n^{\beta}} + \Eb{\prcond{\#\mathcal{C}^1 = n  {,\mathcal{C}^1 \in E_n}}{\mathcal{B}^1_{n}(h,r,m)}{\bT}\mid \mathcal{A}_n(h,r,m)}^2 \Eb{\sum_{\substack{v_1^{1}, \ldots, v_h^{1}\\  v_1^{2}, \ldots, v_h^{2} \\ \text{disjoint} }}\prcond{v_i^j \in \C^j \forall i,j}{\mathcal{B}_{n}(h,r,m)}{\bT}\mid \mathcal{A}_n(h,r,m)} \\
&\leq C'e^{-c'n^{\beta}} + \Eb{\prcond{\#\mathcal{C}^1 = n  {,\mathcal{C}^1 \in E_n}}{\mathcal{B}^1_{n}(h,r,m)}{\bT}\mid \mathcal{A}_n(h,r,m)}^2 =  C'e^{-c'n^{\beta}} + \Eb{S^{1}_n(h,r,m)}^2.
\end{align*}
     Combining this with \eqref{eqn:var decomp} yields:
    \[
        \Varb{ S_n(h,r,m) \mid \mathcal{A}_n(h,r,m)} \le C'e^{-c'n^{\beta}}.
    \]
    Consequently, we have, for all sufficiently large $n$:
    \begin{align*}
        \mathbf{P}_{\alpha}(\exists (h,r,m) \in I_n : \mathcal{A}_n(h,r,m) \text{ occurs} & \text{ and } |S_n(h,r,m)- \Eb{S_n(h,r,m)}| > n^{-3}) \\
        & \le n^3 \sup_{\substack{(h,r,m)\in I_n \\ \mathbf{P}_{\alpha}(\mathcal{A}_n(h,r,m)) >0 }} \frac{\Varb{ S_n(h,r,m) \mid \mathcal{A}_n(h,r,m)}}{n^{-6}} \\
        & \le C'n^{9}e^{-c'n^{\beta}}.
    \end{align*}
    By the Borel-Cantelli lemma, $\mathbf{P}$-almost surely, for $n$ large enough, for all $(h,r,m) \in I_n$, if $\mathcal{A}_n(h,r,m)$ occurs, then $|S_n(h,r,m)- \Eb{S_n(h,r,m)}| \le n^{-3}$. This establishes the second equality in \eqref{eqn:equalities key lem}, as required.
    
%    
%    Thus, going back to \eqref{starting_equation }, for $n$ large enough, we have:
%    \begin{align}\label{bound}
%        \nonumber \sum_{(h,r,m) \in I_n} p_n(h,r,m) S_n(h,r,m) & = \sum_{(h,r,m)\in I_n} \indi{\mathcal{A}_n(h,r,m)^{c}} p_n(h,r,m)S_n(h,r,m) \\
%        \nonumber &+ o(n^{-2}) + \sum_{(h,r,m) \in I_n} \indi{\mathcal{A}_n(h,r,m)} p_n(h,r,m)\Eb{S_n(h,r,m)} \\
%        &= o(n^{-2}) + \sum_{(h,r,m) \in I_n} \indi{\mathcal{A}_n(h,r,m)} p_n(h,r,m)\Eb{S_n(h,r,m)}.
%    \end{align}
%    The last equality follows from the fact that if $\mathcal{A}_n(h,r,m)$ does not occur, then either $p_n(h,r,m) \le n^{-20}$ or $S_n(h,r,m) \le n^{-20}$. 
%    Using a similar argument, we obtain:
%    \[
%    \sum_{(h,r,m) \in I_n} p_n(h,r,m)\indi{\mathcal{A}_n(h,r,m)} = o(n^{-2}) + \sum_{(h,r,m) \in I_n} p_n(h,r,m). \red{???}
%    \]
%    Combining these estimates leads to:
%    \[
%      \sum_{(h,r,m)\in I_n} p_n(h,r,m) S_n(h,r,m) = o(n^{-2}) + \sum_{(h,r,m)\in I_n} p_n(h,r,m)\Eb{S_n(h,r,m)}.
%    \]
Recalling that $\displaystyle \sum_{(h,r,m)\in I_n} p_n(h,r,m) = \Pt{D_n, \#\mathcal{C}\ge (1-\varepsilon)n,  {k_{n,\varepsilon} \ge n^{\frac{1}{2}-\frac{1}{2\alpha}}}}$, let us define
    \begin{align*}
    F_{\mathbf{T}}(n) = F_{\mathbf{T}}(n, D_n, E_n) = \frac{\sum_{(h,r,m)\in I_n} p_n(h,r,m)\Eb{S_n(h,r,m)}}{\Pt{D_n, \#\mathcal{C}\ge (1-\varepsilon)n,  {k_{n,\varepsilon} \ge n^{\frac{1}{2}-\frac{1}{2\alpha}}}}}.
    \end{align*}
Then $F_{\mathbf{T}}(n) $ clearly satisfies the conclusion of the lemma and we have established \eqref{eqn:lemma general main statement}, as required. %shown
%  \begin{align*}
%  \Pt{D_n, \# \mathcal{C} = n,  E_n,  \purple{k_{n,\varepsilon}  \ge n^{\frac{1}{2}-\frac{1}{2\alpha}}}} = \Pt{D_n , \#\mathcal{C}\ge (1-\varepsilon)n, \purple{k_{n,\varepsilon} \ge n^{\frac{1}{2}-\frac{1}{2\alpha}}}}F_{\mathbf{T}}(n) + o(n^{-2}).
%  \end{align*}
%  Dividing both sides by $\Pt{\#\mathcal{C}\ge (1-\varepsilon)n} \sim K_{\alpha}n^{-\frac{1}{\alpha}}$ \red{(by Theorem \ref{thm:ConvergenceOfSize})} concludes the proof.

\end{proof}

The first application of Lemma \ref{lemma_exact_approximations} will be to prove an intermediate result on the height difference between $\C_{=n}$ and $\C_{n, \eps}$.

\begin{proposition}\label{bound_proba_exceed_level_kneps}
$\mathbf{P}$-almost surely, as $n \to \infty$, we have
    \begin{align*}
        \Pt{\Height (\mathcal{C}) \ge k_{n,\varepsilon}+{(\varepsilon^{\frac{1}{2}} n)}^{1-\frac{1}{\alpha}}, \#\mathcal{C}=n \mid \# \mathcal{C} \ge (1-\varepsilon)n} = \mathcal{O}(\varepsilon n^{-1}).
    \end{align*}
\end{proposition}

\begin{proof}
For any $n \ge 0$, on the event $\#\mathcal{C} \ge (1-\varepsilon)n$, let us introduce the event $D_n = \Omega_{{\bT}} $ (the entire probability space) and the event $E_n = \{\Height(\mathcal{C}) \ge k_{n,\varepsilon}+{(\varepsilon^{\frac{1}{2}} n)}^{1-\frac{1}{\alpha}}\}$. Recall that the definitions of $ p_n(h,r,m,D_n)$ and $S_n(h,r,m,  E_n,D_n)$ were given in \eqref{definition_Sp}. The family of events $(D_n)_{n \ge 1}$ and $(E_n)_{n \ge 1}$ considered satisfy the hypothesis of Lemma~\ref{lemma_exact_approximations}. Thus, using Lemma~\ref{lemma_exact_approximations}, we obtain:
    \begin{align}\label{to_bound_second_part}
               &\Pt{E_n,D_n, \#\mathcal{C}=n \mid \# \mathcal{C} \ge (1-\varepsilon)n} = (1+o(1))F_{\mathbf{T}}(n)+ o(n^{-1}) ,
    \end{align}
where $F_{\mathbf{T}}(n)$ satisfies
\begin{align*}
F_{\mathbf{T}}(n) \leq \sup_{\substack{ (h,r,m): \\ p_n(h,r,m,D_n) > 0}} \Eb{S_n(h,r,m,  E_n,D_n)}.
\end{align*}

    For any $(h,r,m) \in \llbracket 1,n \rrbracket^3$ such that $p_n(h,r,m,D_n) > 0$, we have:
    \[
     \Eb{S_n(h,r,m,E_n,D_n)} = \Pf_{\alpha}\left(\sum_{k=1}^{r} X_k = n-m, \max_{1 \le k \le r} (H_k) \ge {(\varepsilon^{\frac{1}{2}} n)}^{1-\frac{1}{\alpha}}\right),
    \]
   where $(X_k, H_k)_{k\ge 1}$ are i.i.d.\ random variables distributed as $(\#\mathcal{C}, \Height(\mathcal{C}))$ under $\Pf_{\alpha}$. 
   
 We now bound this latter quantity. We fix $(h,r,m) \in \llbracket 1, n \rrbracket^3$. Let us denote by $N \ge 1$ the first index such that $H_N \ge (\varepsilon^{\frac{1}{2}} n)^{1-\frac{1}{\alpha}}$, and fix $j \ge 1$. Conditionally on $N =j$, the random variables $(X_k, H_k)_{k\ge 1}$ are again independent. Moreover, for $1 \le k \le j-1 $ the variables $(X_k,H_k)$ are conditioned to satisfy $H_k < (\varepsilon^{\frac{1}{2}} n)^{1-\frac{1}{\alpha}}$, the variable $(X_j,H_j)$ is conditioned to satisfy $H_k \ge (\varepsilon^{\frac{1}{2}} n)^{1-\frac{1}{\alpha}}$ and the variables $(X_k,H_k)$ with $k >j$ are again distributed as under $(\#\mathcal{C}, \Height (\mathcal{C}))$ under $\Pf_{\alpha}$. Now, on the event on which $ \sum_{k=1}^{r} X_k = n-m$ and $ \max_{1 \le k \le r} (H_k) \ge {(\varepsilon^{\frac{1}{2}} n)}^{1-\frac{1}{\alpha}}$ and conditionally on $\somme{\substack{k=1\\ k\neq j}}{r}{X_k}$ we have $X_j = n-m - \somme{\substack{k=1\\ k\neq j}}{r}{X_k} \in [0, \eps n]$. By conditioning on the value of $\somme{\substack{k=1\\ k\neq j}}{r}{X_k}$, we thus obtain the following bound:
\begin{align*}
 \Eb{{S}_n(h,r,m, E_n, D_n)} \le \sup_{0 \le s \le \varepsilon n }\Pf_{\alpha}(X_1 = s \mid H_1 \ge{(\varepsilon^{\frac{1}{2}} n)}^{1-\frac{1}{\alpha}} ).
\end{align*}
Now, we fix $0 \le s \le \varepsilon n $. By \cite[Theorem 2]{Igor_stable_trees_estimates}, there exist positive constants $C_1,C_2 > 0$ (independent of $n$ and $\eps$) such that 
\begin{align*}
\Pf_{\alpha}(  H_1 \ge{(\varepsilon^{\frac{1}{2}} n)}^{1-\frac{1}{\alpha}} \mid X_1 = s) \le C_1 \exp\bigg(-C_2 \varepsilon^{ \frac{\alpha-1}{4}}\bigg(\frac{n}{s}\bigg)^{\frac{\alpha-1}{2}}\bigg).
\end{align*}
Moreover, we have 
\begin{align*}
\Pf_{\alpha}( H_1 \ge{(\varepsilon^{\frac{1}{2}} n)}^{1-\frac{1}{\alpha}} ) \sim C_{\alpha}\varepsilon^{-\frac{1}{2\alpha}} n^{-\frac{1}{\alpha}} \text{ and } \Pf_{\alpha}( X_1 = s ) = \mathcal{O}\bigg( s^{-1-\frac{1}{\alpha}}\bigg),
\end{align*}
where $C_{\alpha} > 0$ (see Table~\ref{table:annealed quenched} for the first of these, and the second follows by a local limit theorem applied to the Lukasiewicz path). Putting all these equations together, we obtain 
\begin{align*}
\Pf_{\alpha}(X_1 = s \mid H_1 \ge{(\varepsilon^{\frac{1}{2}} n)}^{1-\frac{1}{\alpha}} ) &\le C_{3}\frac{s^{-\frac{\alpha+1}{\alpha}}}{\varepsilon^{-\frac{1}{2\alpha}} n^{-\frac{1}{\alpha}}}\exp\bigg(-C_2 \varepsilon^{ \frac{\alpha-1}{4}}\bigg(\frac{n}{s}\bigg)^{\frac{\alpha-1}{2}}\bigg)\\
& \le C_3 n^{-1} \bigg( \frac{s}{n}\bigg)^{-\frac{\alpha+1}{\alpha}}\exp\bigg(-C_2 \varepsilon^{ \frac{\alpha-1}{4}}\bigg(\frac{n}{s}\bigg)^{\frac{\alpha-1}{2}}\bigg).
\end{align*}
Then, let us denote by $a = \frac{\alpha+1}{\alpha}$, $b = C_2\varepsilon^{ \frac{\alpha-1}{4}}$ and $c = \frac{\alpha-1}{2}$ and $\varphi :x \mapsto x^{-a}\exp(-bx^{-c})$ for $x > 0$. A simple study of $\varphi$ shows that there exists $\varepsilon_{\alpha} > 0$ such that for $0 < \varepsilon < \varepsilon_{\alpha}$ (which we assume without loss of generality), the function is strictly increasing on $(0,\varepsilon]$. It follows that in the last equation, the right-hand side is maximal for $s = \varepsilon n$. This leads to the upper bound
\begin{align*}
C_3 n^{-1}\varepsilon^{-\frac{\alpha+1}{\alpha}}\exp\bigg(-C_2 \varepsilon^{- \frac{\alpha-1}{4}}\bigg) = \mathcal{O}(\varepsilon n^{-1}),
\end{align*}
and thus establishes that
\begin{align*}
 \Eb{{S}_n(h,r,m, E_n, D_n)}  = \mathcal{O}(\varepsilon n^{-1}).
\end{align*}
Substituting this into \eqref{to_bound_second_part}, we deduce that
\begin{align*}
 \Pt{\Height (\mathcal{C}) \ge k_{n,\varepsilon}+{(\varepsilon^{\frac{1}{2}} n)}^{1-\frac{1}{\alpha}}, \#\mathcal{C}=n \mid \# \mathcal{C} \ge (1-\varepsilon)n} = \mathcal{O}(\varepsilon n^{-1}),
\end{align*}
thus concluding the proof.
\end{proof}

We now have all the ingredients to conclude the proof of Proposition~\ref{estimateA_i}. We remind the reader to keep in mind throughout that the parameters $\eps, \eta, \delta>0$ respect the following ordering: $\delta \ll \eta \ll \varepsilon$. 

\begin{proof}[Proof of Proposition~\ref{estimateA_i}]
\begin{enumerate}[(a)]
\item We first claim that, given $\eps>0$, there exists $C_{\eps} < \infty$ such that, for all sufficiently large $n$:
\begin{align}\label{eqn:exp univ UB}
\sup_{(h,r,m) : p_n^0(h,r,m) > 0} \Eb{S_n^0(h,r,m)} \leq C_{\eps} n^{-1}.
\end{align}
To show this, we bound 
\[
\sup_{(h,r,m) : p_n^0(h,r,m) > 0} \Eb{S_n^0(h,r,m)} \leq \sup_{\substack{r>0 \\ \eps n - r \leq m' \leq \eps n}}  \Pf_{\alpha}\bigg(\sum_{k=1}^{r}{X_k} = m' \bigg).
\]
For $r \geq \eps n^{\frac{1}{\alpha}}$, it follows from a stable local limit theorem that there exists a constant $C_{\eps}$ such that for all sufficiently large $n$,
\[
\sup_{r  \geq \eps n^{\frac{1}{\alpha}}, m' \leq \eps n}  \Pf_{\alpha}\bigg(\sum_{k=1}^{r}{X_k} = m' \bigg) \leq C_{\eps} n^{-1} ||g_{\alpha}||_{\infty} 
\]
where $g_{\alpha}$ is the density of the limiting positive $\frac{1}{\alpha}$-stable random variable (again see \cite[p. 236]{gnedenko1968limit} for the details of the local limit theorem).

For $r <n^{\frac{1}{\alpha}}$, we have $m' \geq \eps n /2$ and the result therefore follows from \cite[Theorem 2.4]{bergercauchy} - note that \cite[Assumption (2.5)]{bergercauchy} is satisfied as a consequence of Gnedenko's local limit theorem - see \cite[Lemma A3(i)]{curien2015percolation} for details of how this follows from the local limit theorem (and note that the slightly stronger assumption on the offspring tails there is not necessary for the local limit theorem used in the proof). In particular, \cite[Theorem 2.4]{bergercauchy} implies that there exists $C_{\eps} < \infty$ such that
\[
\sup_{r  \leq \eps n^{\frac{1}{\alpha}}, m' \in [\eps n/2, \eps n]}  \Pf_{\alpha}\bigg(\sum_{k=1}^{r}{X_k} = m' \bigg) \leq C r (\eps n)^{-1-\frac{1}{\alpha}} \leq C_{\eps} n^{-1}.
\]
%A similar reasoning yields the bound:
%  $$\sup_{(h,r,m) : p_n^0(h,r,m) > 0} \Eb{S_n^0(h,r,m)} \le \mathcal{O}_{\varepsilon}(n^{-1}).$$
 This establishes \eqref{eqn:exp univ UB}. Consequently, for $i = 0$, using Proposition \ref{events_construction}(a) we have:
    \begin{align*}
    \Pt{D_n^0, \# \mathcal{C} = n \mid \#\mathcal{C}\ge (1-\varepsilon)n} = \mathcal{O}_{\varepsilon}(\eta n^{-1}).
    \end{align*}
    Similarly for the annealed probability, we have
\begin{align*}   
   &\mathbb{P}_{\alpha} (\mathcal{C}_{n,\varepsilon} \in A_0,\#\mathcal{C}=n  \mid \# \mathcal{C} \ge (1-\varepsilon)n) \\
   &=      \mathbb{P}_{\alpha} (\#\mathcal{C}=n  \mid \mathcal{C}_{n,\varepsilon} \in A_0, \# \mathcal{C} \ge (1-\varepsilon)n) \mathbb{P}_{\alpha} (\mathcal{C}_{n,\varepsilon} \in A_0 \mid \# \mathcal{C} \ge (1-\varepsilon)n) \\
   &\leq \sup_{\substack{r>0 \\ \eps n - r \leq m' \leq \eps n}}   \Pf_{\alpha}\bigg(\sum_{k=1}^{r}{X_k} = m' \bigg) \mathbb{P}_{\alpha} (\mathcal{C}_{n,\varepsilon} \in A_0 \mid \# \mathcal{C} \ge (1-\varepsilon)n) = \mathcal{O}_{\varepsilon}(\eta n^{-1}).
\end{align*}
Here the final bound again follows using the estimates above and Proposition \ref{events_construction}(a).
\item For $i \in \{0,\cdots,N\}$, for any $n \ge 0$, on the event $\#\mathcal{C} \ge (1-\varepsilon)n$, let us introduce the event $D_n^i = \{\mathcal{C}_{n,\varepsilon} \in A_i\} $ and the event $E_n = \Omega_{{\bT}}$. 

We abbreviate $ p^i_n(h,r,m):=p_n(h,r,m,{D_n^i})$ and $S^i_n(h,r,m) : = S_n(h,r,m,  E_n,{D_n^i})$, recalling the original definitions in \eqref{definition_Sp}. For a fixed $i$, the families of events $(D^i_n)_{n \ge 1}$ and $(E_n)_{n \ge 1}$ considered satisfy the hypotheses of Lemma~\ref{lemma_exact_approximations}. Thus by Lemma~\ref{lemma_exact_approximations} and \eqref{eqn:exp univ UB}, $\mathbf{P}_{\alpha}$-almost surely, for any $i \in \{0,\dots,N\}$, there exists a random variable $F_{\mathbf{T}}(n, A_i)$ such that:
    \begin{align}\label{recall_lemma}
    \Pt{D_n^i, \# \mathcal{C} = n \mid \#\mathcal{C}\ge (1-\varepsilon)n} = \Pt{D_n^i \mid \#\mathcal{C}\ge (1-\varepsilon)n} F_{\mathbf{T}}(n, A_i) + o_{\eps}(n^{-1}),
    \end{align}
    where $F_{\mathbf{T}}(n, A_i)$ satisfies:
    \[
    \inf_{(h,r,m) : p_n^i(h,r,m)  >  0} \Eb{S_n^i(h,r,m)} \leq F_{\mathbf{T}}(n, A_i) \leq \sup_{(h,r,m) : p_n^i(h,r,m) > 0} \Eb{S_n^i(h,r,m)}.
    \]
We claim, for $i \neq 0$, that the upper and lower bounds appearing above are very similar. In particular, for any $(h,r,m) \in \llbracket 1, n \rrbracket^3$, we have:
$$\Eb{S_n^i(h,r,m)} = \Pf_{\alpha}\bigg(\sum_{k=1}^{r}{X_k} = n -m\bigg),$$
where $(X_k)_{k \ge 1}$ denotes a family of independent random variables distributed as $\#\mathcal{C}$ under $\Pf_{\alpha}$. Recall that if $p_n^i(h,r,m) > 0$, then $m \ge (1-\varepsilon)n$ and $r \ge m - (1-\varepsilon)n$. For any $i \in \{1,\dots,n\}$, since $\mathrm{diam}_D(A_i) \leq \delta$, the admissible values for $(h,r,m)$ such that $p_n^i(h,r,m) > 0$ satisfy:
\begin{align*}
&r \in R_{n,i} := [r^{i}_{\min} n^{\frac{1}{\alpha}}, r^{i}_{\max}n^{\frac{1}{\alpha}}], \\
&m \in M_{n,i} :=  [ (1-\varepsilon)n, (1-\varepsilon)n+r^{i}_{\max}n^{\frac{1}{\alpha}}],
\end{align*}
with $r^i_{\max} - r^i_{\min} \le \delta$, and $ r^i_{\min},  r^i_{\max} \in [K_{\alpha, \eps, \eta}^{-1}, K_{\alpha, \eps, \eta}]$, where $K = K_{\alpha, \eps, \eta}>0$ is as in Proposition \ref{events_construction}(b). This yields the following bounds (note that the same bounds hold directly in the annealed case):
\begin{align}\label{encadrement_annealed_and_quenched}
    \nonumber \inf_{\substack{(h,r,m) \\ m \in M_{n,i,\varepsilon,\eta,\delta} \\ r \in  R_{n,i,\varepsilon,\eta,\delta} }} \Eb{S_n^i(h,r,m)} &\leq \Pf_{\alpha}(\#\mathcal{C} = n \mid D_n^i, \#\mathcal{C} \ge (1-\varepsilon)n)\leq \sup_{\substack{(h,r,m) \\ m \in M_{n,i} \\ r \in  R_{n,i} }} \Eb{S_n^i(h,r,m)}, \\
     \inf_{\substack{(h,r,m) \\ m \in M_{n,i} \\ r \in  R_{n,i} }} \Eb{S_n^i(h,r,m)} &\leq F_{\mathbf{T}}(n, A_i) \leq \sup_{\substack{(h,r,m) \\ m \in M_{n,i} \\ r \in  R_{n,i} }} \Eb{S_n^i(h,r,m)}. 
\end{align}
By Theorem~\ref{thm:ConvergenceOfSize}, we have $\Pf(X_1 \ge k) \sim K_{\alpha}k^{-1/\alpha}$ as $k \to \infty$.
Using \cite[Section 50]{gnedenko1968limit}, we obtain the following asymptotic which holds uniformly for $m \in M_{n,i}$ and $r \in R_{n,i}$:
\begin{align*}
\Eb{S_n^i(h,r,m)} = r^{-\alpha}g_{\alpha}\bigg(\frac{n-m}{r^{\alpha}} \bigg) + o(n^{-1})= r^{-\alpha}g_{\alpha}\bigg(\frac{\varepsilon n -r_{m,n}}{r^{\alpha}} \bigg) + o(n^{-1}),
\end{align*}
where $r_{m,n} = m-(1-\eps)n \in[0, r^i_{\max}n^{\frac{1}{\alpha}}]$ and 
where $g_{\alpha}$ is the density of the positive $\frac{1}{\alpha}$-stable random variable arising as the limit of $m^{-\alpha}\sum_{j=1}^m X_j$ (see \cite[Section 50]{gnedenko1968limit} for details). Using \eqref{eqn:stable density ratio}, we obtain 
\begin{align}\label{to_check}
 \sup_{\substack{(h,r,m) \\ m \in M_{n,i,\varepsilon,\eta,\delta} \\ r \in  R_{n,i,\varepsilon,\eta,\delta} }} \Eb{S_n^i(h,r,m)} =  \inf_{\substack{(h,r,m) \\ m \in M_{n,i,\varepsilon,\eta,\delta} \\ r \in  R_{n,i,\varepsilon,\eta,\delta} }} \Eb{S_n^i(h,r,m)}(1+\mathcal{O}_{\eps,\eta}(\delta))+o(n^{-1}) .
\end{align}
Combining \eqref{encadrement_annealed_and_quenched} with \eqref{to_check}, we deduce that
  \begin{align}\label{estimate_F}
   F_{\mathbf{T}}(n, A_i)  =\Pf_{\alpha}(\#\mathcal{C} = n \mid D_n^i, \#\mathcal{C} \ge (1-\varepsilon)n)(1+\mathcal{O}_{\varepsilon,\eta}(\delta)) + o(n^{-1}).
  \end{align}
Using Item $(e)$ of Proposition~\ref{events_construction}, we have
$$\Pt{D_n^i \mid \#\mathcal{C}\ge (1-\varepsilon)n} \sim \mathbb{P}_{\alpha}(D_n^i\mid \#\mathcal{C}\ge (1-\varepsilon)n) .$$
Combining these estimates with \eqref{recall_lemma}, we obtain that for any $i \in \{1,\dots,N\}$, for $n$ large enough
\begin{align*}
\Pt{D_n^i, \# \mathcal{C} = n \mid \#\mathcal{C}\ge (1-\varepsilon)n} = \Pf_{\alpha}(D_n^i,\#\mathcal{C} = n \mid  \#\mathcal{C} \ge (1-\varepsilon)n)(1+\mathcal{O}_{\varepsilon,\eta}(\delta)) + o_{\eps}(n^{-1}).
\end{align*}  
This proves the second part of the proposition. 
\item Using part (b), we have
    \begin{align*}
    \Pt{\#\mathcal{C}=n  \mid \# \mathcal{C} \ge (1-\varepsilon)n} &= \somme{i=0}{N}{\Pt{\mathcal{C}_{n,\varepsilon} \in A_i,\#\mathcal{C}=n  \mid \# \mathcal{C} \ge (1-\varepsilon)n}}\\
    &\ge o_{\eps}(n^{-1})+(1+\mathcal{O}_{\varepsilon,\eta}(\delta))\somme{i=1}{N}{\mathbb{P}_{\alpha} (\mathcal{C}_{n,\varepsilon} \in A_i,\#\mathcal{C}=n  \mid \# \mathcal{C} \ge (1-\varepsilon)n)}.
    \end{align*}
The sum can be rewritten as
\begin{align*}
\mathbb{P}_{\alpha} (\#\mathcal{C}=n  \mid \# \mathcal{C} \ge (1-\varepsilon)n) - \mathbb{P}_{\alpha} (\mathcal{C}_{n,\varepsilon} \in A_0,\#\mathcal{C}=n  \mid \# \mathcal{C} \ge (1-\varepsilon)n) = \Theta(n^{-1}) - \mathcal{O}_{\varepsilon}(\eta n^{-1}).
\end{align*}
Combining the last two equations and then letting $\delta \to 0$, then $\eta \to 0$ we deduce that there exists an absolute constant $C_{\alpha}' > 0$ such that for all $n$ large enough,
\begin{align*}%\label{lower_bound_proba_size_exactly_n}
\Pt{\#\mathcal{C}=n  \mid \# \mathcal{C} \ge (1-\varepsilon)n} \ge C_{\alpha}' n^{-1}.
\end{align*}
We also recall the bound from Proposition~\ref{bound_proba_exceed_level_kneps}:
    \begin{align*}
        \Pt{\Height (\mathcal{C}) \ge k_{n,\varepsilon}+(\varepsilon^{\frac{1}{2}} n)^{1-\frac{1}{\alpha}}, \#\mathcal{C}=n \mid \# \mathcal{C} \ge (1-\varepsilon)n} = \mathcal{O}(\varepsilon n^{-1}).
    \end{align*}
Putting together the two last estimates, and recalling that $F$ is Lipschitz, we deduce that
\begin{align*}
\mathbb{E}_{\mathbf{T}}\bigg[F(\mathcal{C}) \mid \# \mathcal{C} = n \bigg] =\mathbb{E}_{\mathbf{T}}\bigg[F(\mathcal{C}_{n,\varepsilon}) \mid \# \mathcal{C} = n \bigg] + \mathcal{O}(\varepsilon^{\frac{1}{2}(1-\frac{1}{\alpha})}) +\mathcal{O}(\varepsilon).
\end{align*}
    \end{enumerate}
\end{proof}

\section{Conditioning on the height}\label{section_conditioning_on_the_height}

In this section, we prove the analogue of Theorem~\ref{theorem_exactconditioningsize}, but with the height instead of the size. Since the reasoning is very similar, we provide only the main intermediate steps and detail the differences in the proof.

Conditionally on $\mathbf{T}$, for any $n \ge 0$, we recall that $\Chn$ denotes the cluster $\mathcal{C}$ conditioned to have height equal to $n$ under $\PfT$. We denote by $\mathcal{T}_{\alpha}^{H=1}$ the stable tree with parameter $\alpha$ of total height $1$. This section is dedicated to outlining the proof of the following theorem.

\begin{theorem}\label{theorem_exactconditioningheight}
Take $\gamma$ as in \eqref{eqn:annealed GHP}. Then, for $\bPb$-almost every $\bT$, the following convergence holds in law under $\mathbb{P}_{\bT}$:
\begin{align*}
(\Chn, n^{-1}d_n,{(\gamma n)^{-\frac{\alpha}{\alpha - 1}}} \nu_n,\rho_n) &\underset{n \to +\infty}{\overset{(d)}{\longrightarrow}} (\mathcal{T}_{\alpha}^{H=1},d_{\mathcal{\mathcal{T}_{\alpha}}},\nu_{\alpha},\rho_{\alpha})
\end{align*}
with respect to the pointed Gromov-Hausdorff-Prokhorov topology.
\end{theorem}

We give a lemma that provides the different inputs required in this case. Note that the other estimates of Lemma \ref{eventHandS} transfer more directly since Lemma \ref{eventHandS}(a) did not depend on the conditioning, and $k_{n, \eps}$ is now deterministic there is no need for Lemma \ref{eventHandS}(b,c).

\begin{lemma}\label{lemma_height}
Fix $\eps>0$.
\begin{enumerate}[(a)]
\item $\bPb$-almost surely, $\Pt{Y_{\lfloor (1-\eps)n \rfloor} \geq n^{10} | \Height(\C) \geq (1-\eps)n } = o(n^{-1})$.
\item We have that
 \begin{align*}
   \sup_{m \geq 1} \Pf_{\alpha}\bigg(\sum_{k=1}^{m} X_k \ge (\varepsilon^{\frac{1}{2}}n)^{\frac{\alpha}{\alpha-1}},\max_{1 \le k \le m}H_k = \varepsilon n \bigg) = \mathcal{O} (\eps n^{-1}),
    \end{align*}
    where $(X_k,H_k)_{k \ge 1}$ is a family of i.i.d. random variables distributed as $(\#\mathcal{C},\Height (\mathcal{C}))$ under $\Pf_{\alpha}$.
    \end{enumerate}
\end{lemma}
\begin{proof}
\begin{enumerate}[(a)]
\item By Markov's inequality, the probability in question is upper bounded by $n^{-10}\Et{Y_{\lfloor (1-\eps)n \rfloor}}$. Moreover, $\Eb{ \Et{Y_{\lfloor (1-\eps)n \rfloor}}} =1$ and hence, by another application of Markov's inequality and Borel-Cantelli, we have that $\Et{Y_{\lfloor (1-\eps)n \rfloor}} \leq n^2$, eventually almost surely.
\item At several times in this next proof, we will use the fact that $X_k$ and $H_k$ are positively correlated under $\mathbb{P}_{\alpha}$: in particular, for any $m_2 \geq m_1 \geq 1$, the law of $X_1$ conditioned on $H_1 < m_2$ stochastically dominates its law conditioned on $H_1 < m_1$. This can, for example, be seen from a spinal decomposition of Galton-Watson trees along their height (see \cite{geiger1998galton}).

We decompose the probability as 
\begin{align}\label{eqn:prob factor}
\prcond{\sum_{k=1}^{m} X_k \ge (\varepsilon^{\frac{1}{2}}n)^{\frac{\alpha}{\alpha-1}}}{\max_{1 \le k \le m}H_k = \varepsilon n }{\alpha} \pr{\max_{1 \le k \le m}H_k = \varepsilon n }.
\end{align}
Note that the latter probability can be bounded by:
\begin{align}\label{eqn:height equal decomp}
m\pr{H_1 = \eps n} \pr{H_1 \leq \eps n}^{m-1} \leq cC_{m,n}  \eps^{-1-\frac{1}{\alpha -1}}  \exp \{ -C_{m,n} \eps^{-\frac{1}{\alpha-1}}\}n^{-1},
\end{align}
where $C_{m,n} = {m}n^{-\frac{1}{\alpha-1}}$ (here we use the known fact that there exists $c< \infty$ such that $\pr{H_1 = N} \leq cN^{-1-\frac{1}{\alpha-1}}$ - this follows from the second line of Table \ref{table:annealed quenched} and the fact that this probability is non-increasing in $N$). In particular, when $C_{n,m} \geq \eps^{\frac{3}{4(\alpha-1)}}$, this already gives the desired upper bound, so we may assume henceforth that this is not the case.
Turning now to the first factor in \eqref{eqn:prob factor}, we can write
\begin{align*}
&\prcond{\sum_{k=1}^{m} X_k \ge (\varepsilon^{\frac{1}{2}}n)^{\frac{\alpha}{\alpha-1}}}{\max_{1 \le k \le m}H_k = \varepsilon n }{\alpha} \\
&\quad \leq \prcond{ X_1 \ge \frac{1}{3}(\varepsilon^{\frac{1}{2}}n)^{\frac{\alpha}{\alpha-1}}}{H_1 = \varepsilon n }{\alpha}  \\
&\qquad + \prcond{\sum_{k=1}^{m-1} X_k \mathbbm{1}\left\{X_k \leq  \frac{1}{3} (\varepsilon^{\frac{3}{4}}n)^{\frac{\alpha}{\alpha-1}}\right\} \geq (\varepsilon^{\frac{1}{2}}n)^{\frac{\alpha}{\alpha-1}}}{\max_{1 \le k \le m}H_k \leq \varepsilon n }{\alpha} \\
&\qquad + \prcond{\sum_{k=1}^{m-1} X_k \mathbbm{1}\left\{X_k >  \frac{1}{3} (\varepsilon^{\frac{3}{4}}n)^{\frac{\alpha}{\alpha-1}}\right\} \geq (\varepsilon^{\frac{1}{2}}n)^{\frac{\alpha}{\alpha-1}}}{\max_{1 \le k \le m}H_k \leq \varepsilon n }{\alpha}.
\end{align*}
We now treat these probabilities in turn. For the first of these, note that, by the known convergence for conditioning on the height in the annealed case, we have 
\begin{align*}
\prcond{ X_1 \ge \frac{1}{3}(\varepsilon^{\frac{1}{2}}n)^{\frac{\alpha}{\alpha-1}}}{H_1 = \varepsilon n }{\alpha} \to \pr{\nu (\Ta^{H=1}) \geq \frac{1}{3}\varepsilon^{-\frac{\alpha}{2(\alpha-1)}}} \leq \exp\{- \eps^{-\frac{\alpha^2}{2(\alpha-1)^2}} \}
\end{align*}
(this last result can be seen by combining the result of \cite[Theorem 1.8]{duquesnewangdecomp} and \cite[Proposition 5.6]{GoldHaas10}), and hence is upper bounded by $\exp\{- \eps^{-1}\}$ for all sufficiently large $n$.

For the second term, note that, using the fact that $C_{n,m} < \eps^{\frac{3}{4(\alpha-1)}}$, we can write
\begin{align*}
&\prcond{\sum_{k=1}^{m-1} X_k \mathbbm{1}\left\{X_k \leq  \frac{1}{3} (\varepsilon^{\frac{3}{4}}n)^{\frac{\alpha}{\alpha-1}}\right\} \geq (\varepsilon^{\frac{1}{2}}n)^{\frac{\alpha}{\alpha-1}}}{\max_{1 \le k \le m}H_k \leq \varepsilon n }{\alpha}  \\
&\leq \pr{\sum_{k=1}^{(\eps^{\frac{3}{4}}n)^{\frac{1}{\alpha-1}}} X_k \mathbbm{1}\left\{X_k \leq  \frac{1}{3} (\varepsilon^{\frac{3}{4}}n)^{\frac{\alpha}{\alpha-1}}\right\} \geq \eps^{-\frac{\alpha}{4(\alpha-1)}}(\varepsilon^{\frac{3}{4}}n)^{\frac{\alpha}{\alpha-1}}} \\
&\leq \pr{\sum_{k=1}^{(\eps^{\frac{3}{4}}n)^{\frac{1}{\alpha-1}}} X_k \geq \frac{1}{3}(\varepsilon^{\frac{3}{4}}n)^{\frac{\alpha}{\alpha-1}}}^{\eps^{-\frac{\alpha}{4(\alpha-1)}}} =\left(\pr{S_1 \geq \frac{1}{3}}+o(1)\right)^{\eps^{-\frac{\alpha}{4(\alpha-1)}}} %\leq \exp \{c \eps^{-\frac{\alpha}{4(\alpha-1)}}\}
\end{align*}
for all sufficiently large $n$, where $S_1$ is a $\frac{1}{\alpha}$-stable subordinator. Here the final inequality follows by applying the strong Markov property at the times $T_j := \inf \{t > T_{j-1}: \sum_{k=T_{j-1}+1}^t X_k \geq \frac{1}{3}(\varepsilon^{\frac{3}{4}}n)^{\frac{\alpha}{\alpha-1}}\}$, with $T_0=0$, and using that $\sum_{k=T_{j-1}+1}^{T_j} X_k \geq \frac{1}{3}(\varepsilon^{\frac{3}{4}}n)^{\frac{\alpha}{\alpha-1}}$, since each $X_k$ is bounded. Here we also used the fact that $X_k$ under the conditioning $H_k \leq \eps n$ is stochastically dominated by an unconditioned copy of $X_k$, as mentioned at the beginning of the proof.

Finally, for the third term, note that, again applying \cite[Theorem 1.8]{duquesnewangdecomp} , we see that, for all sufficiently large $n$, the probability $\prcond{X_k >  \frac{1}{3} (\varepsilon^{\frac{3}{4}}n)^{\frac{\alpha}{\alpha-1}}}{H_k \leq \varepsilon n }{\alpha} $ is upper bounded by
\begin{align*}
&\pr{X_k >  \frac{1}{3} (\varepsilon n)^{\frac{\alpha}{\alpha-1}}} \prcond{X_k >  \frac{1}{3} (\varepsilon^{\frac{3}{4}}n)^{\frac{\alpha}{\alpha-1}}}{X_k >  \frac{1}{3} (\varepsilon n)^{\frac{\alpha}{\alpha-1}}, H_k \leq \varepsilon n }{\alpha} \\
&\leq c(\varepsilon n)^{-\frac{1}{\alpha-1}} \frac{\prcond{ H_k \leq \varepsilon n }{X_k >  \frac{1}{3} (\varepsilon^{\frac{3}{4}}n)^{\frac{\alpha}{\alpha-1}}}{\alpha} \pr{X_k >  \frac{1}{3} (\varepsilon^{\frac{3}{4}}n)^{\frac{\alpha}{\alpha-1}}}}{\pr{X_k >  \frac{1}{3} (\varepsilon n)^{\frac{\alpha}{\alpha-1}}, H_k \leq \varepsilon n }} \leq \exp\{-c\eps^{-1}\} n^{-\frac{1}{\alpha-1}}.
\end{align*}
This implies that the number of terms contributing to the sum $\sum_{k=1}^{m-1} X_k \mathbbm{1}\left\{X_k >  \frac{1}{3} (\varepsilon^{\frac{3}{4}}n)^{\frac{\alpha}{\alpha-1}}\right\} $ is stochastically dominated by a \textsf{Binomial}($m,  \exp\{-c\eps^{-1}\} n^{-\frac{1}{\alpha-1}}$) random variable. Since we are assuming, without loss of generality, that $C_{n,m} < \eps^{\frac{3}{4(\alpha-1)}}$ this is zero with probability at least $1-\exp\{-c\eps^{-1}\}$.

Returning to \eqref{eqn:prob factor} and \eqref{eqn:height equal decomp} and combining with the above estimates, we deduce that in the case $C_{n,m} < \eps^{\frac{3}{4(\alpha-1)}}$, the probability in question can be upper bounded by 
\[
\exp\{-c\eps^{-1}\} cC_{m,n}  \eps^{-1-\frac{1}{\alpha -1}}  \exp \{ -C_{m,n} \eps^{-\frac{1}{\alpha-1}}\}n^{-1},
\]
which completes the proof.
\end{enumerate}
\end{proof}

To prove Theorem \ref{theorem_exactconditioningheight}, the idea is again to compare quenched expectations with annealed expectations. Fix $\varepsilon, \eta, \delta > 0$ and a family of events $(A_i)_{0 \le i \le N}$ given by {Proposition~\ref{events_construction}}. We consider the cluster $\Chgen$ truncated at level $\lfloor(1-\varepsilon)n\rfloor$. (Compare to the previous section, where we truncated the cluster $\mathcal{C}_{\geq (1-\eps)n}$ at the \textit{random} level $k_{n,\varepsilon}$.) Let ${\C}_{H,n,\varepsilon}$ be the tree $\C$ where the part above level $\lfloor(1-\varepsilon)n\rfloor$ is removed. We first state the proposition analogous to Proposition~\ref{estimateA_i}. 

Note that the $o_{\eps}(1)$ and $o(1)$ terms in \eqref{eqn:lemma general main statement ideal} do not appear when conditioning on the height: this is because $k_{n, \eps}$ is deterministic in this case (so does not need to be controlled separately) and because we will not exactly apply a local limit theorem (which led to the $o(n^{-1})$ term).

\begin{proposition}\label{asymptotics_height}
\begin{enumerate}[(a)]
\item
 $\mathbf{P}_{\alpha}$-almost surely
 \begin{align}\label{estimate_for_i_equal_zero_height}
    \Pt{\mathcal{C}_{n,\varepsilon} \in A_0, \Height ( \mathcal{C}) = n \mid \Height (\mathcal{C})\ge (1-\varepsilon)n} = \mathcal{O}_{\varepsilon}(\eta n^{-1}).
 \end{align}
\item $\mathbf{P}_{\alpha}$-almost surely, for any $i \in \{1,\dots,N\}$, 
\begin{align}\label{estimate_for_i_positive_height}
%\begin{split}
&\Pt{{\mathcal{C}}_{H,n,\varepsilon} \in A_i, \Height (\mathcal{C}) = n \mid \Height (\mathcal{C})\ge (1-\varepsilon)n} \\
&\quad = \Pf_{\alpha}({\mathcal{C}}_{H,n,\varepsilon} \in A_i, \Height (\mathcal{C}) = n \mid \Height (\mathcal{C})\ge (1-\varepsilon)n)(1+\mathcal{O}_{\varepsilon,\eta}(\delta))\nonumber .
%\end{split}
\end{align}
 \item  $\mathbf{P}_{\alpha}$-almost surely,
 \begin{align*}
\mathbb{E}_{\mathbf{T}}\bigg[F(\mathcal{C}) \mid  \Height (\mathcal{C})=n  \bigg] =\mathbb{E}_{\mathbf{T}}\bigg[F({\mathcal{C}}_{H,n,\varepsilon}) \mid  \Height (\mathcal{C}) = n \bigg] + \mathcal{O}(\varepsilon) .
\end{align*}
\end{enumerate}
    \end{proposition} 
\begin{proof}
We start by proving the first two items (a) and (b). The proof closely follows that of Lemma~\ref{lemma_exact_approximations}. For $i \in \{0,\dots,N\}$ and $m \ge 1$, we consider the quantities:
\begin{align*}
p_n^i(m) &= \Pt{{\mathcal{C}}_{H,n,\varepsilon} \in A_i, Y_{\lfloor(1-\varepsilon)n\rfloor} = m \mid \Height (\mathcal{C}) \ge (1-\varepsilon)n}, \\
S_n^i(m) &= \Pt{\Height (\mathcal{C}) = n \mid{\mathcal{C}}_{H,n,\varepsilon} \in A_i, Y_{\lfloor(1-\varepsilon)n\rfloor} = m , \Height (\mathcal{C}) \ge (1-\varepsilon)n}.
\end{align*}
Then, for any $i \in \{0,\dots,N\}$, we can write:
\begin{align*}
\Pt{{\mathcal{C}}_{H,n,\varepsilon} \in A_i, \Height (\mathcal{C}) = n \mid \Height (\mathcal{C})\ge (1-\varepsilon)n} = \sum_{m \ge 1} p^i_n(m)S_n^i(m).
\end{align*}
Using point (a) of Lemma~\ref{lemma_height}, for $n$ large enough, we have
\begin{align*}
 \sum_{m > n^{10}} p^i_n(m) = o(n^{-1}).
\end{align*}
Now, for $m \in \{1,\dots,n^{10}\}$, considering the event $\mathcal{A}(m)$ on which $p_n^i(m) \ge n^{-20}$ and following the same steps as in the proof of Lemma~\ref{lemma_exact_approximations}, we obtain that $\mathbf{P}_{\alpha}$-almost surely:
\begin{align*}
\Pt{{\mathcal{C}}_{H,n,\varepsilon} \in A_i, \Height (\mathcal{C}) = n \mid \Height (\mathcal{C})\ge (1-\varepsilon)n} = \sum_{m \ge 1} p^i_n(m)\mathbf{E}_{\alpha}[S_n^i(m)]+o(n^{-1}).
\end{align*}
This implies that, for all $i \in \{0,\dots,N\}$,
\begin{align*}
&\Pt{{\mathcal{C}}_{H,n,\varepsilon} \in A_i, \Height (\mathcal{C}) = n \mid \Height (\mathcal{C})\ge (1-\varepsilon)n} \\
&\qquad = \Pt{{\mathcal{C}}_{H,n,\varepsilon} \in A_i \mid \Height (\mathcal{C})\ge (1-\varepsilon)n}F_{\mathbf{T}}(n,\delta,\varepsilon,i)+o(n^{-1}), 
\end{align*}
where $F_{\mathbf{T}}(n,\delta,\varepsilon,i)$ satisfies:
\begin{align*}
\inf_{m \ge 1 : p^i_n(m) > 0}\mathbf{E}_{\alpha}[S_n^i(m)] \le F_{\mathbf{T}}(n,\delta,\varepsilon,i) \le \sup_{m \ge 1 : p^i_n(m) > 0}\mathbf{E}_{\alpha}[S_n^i(m)].
\end{align*}
We again have the upper bound of $F_{\mathbf{T}}(n,\delta,\varepsilon,0) = \mathcal{O}_{\eps}(n^{-1})$ since there exist universal constants $c_{\eps}, C_{\eps}$ and $C_{\eps}'$ such that for all $m \geq 1$,
\[
\mathbf{E}_{\alpha}[S_n^i(m)] \leq m\pr{H \leq \eps n}^{m-1}\pr{H = \lceil \eps n \rceil} \leq m C_{\eps} \exp \{-c_{\eps} m n^{-\frac{1}{\alpha-1}}\} n^{-1-\frac{1}{\alpha-1}} \leq C_{\eps}'  n^{-1},
\]
from which the result of part (a) follows by following the same logic as in the proof of Proposition~\ref{estimateA_i}(a). Similarly, when we restrict to $m \in [r^{i}_{\min} n^{\frac{1}{\alpha}}, r^{i}_{\max}n^{\frac{1}{\alpha}}]$ with $r^i_{\max} - r^i_{\min} \le \delta$ and $ r^i_{\min},  r^i_{\max} \in [K_{\alpha, \eps, \eta}^{-1}, K_{\alpha, \eps, \eta}]$, we obtain 
\begin{align*}
 \sup_{\substack{ m \in  R_{n,i,\varepsilon,\eta,\delta} }} \Eb{S_n^i(m)} =  \inf_{\substack{ m \in  R_{n,i,\varepsilon,\eta,\delta} }} \Eb{S_n^i(m)}(1+\mathcal{O}_{\eps,\eta}(\delta)),
\end{align*}
from which part (b) follows by again following the same logic as in the proof of Proposition~\ref{estimateA_i}.

Now let us prove item (c). It remains to show that, $\mathbf{P}_{\alpha}$-almost surely,
 \begin{align*}
\mathbb{E}_{\mathbf{T}}\bigg[F(\mathcal{C}) \mid  \Height (\mathcal{C})=n  \bigg] =\mathbb{E}_{\mathbf{T}}\bigg[F({\mathcal{C}}_{H,n,\varepsilon}) \mid  \Height (\mathcal{C}) = n \bigg] + \mathcal{O}(\varepsilon) .
\end{align*}
    This is the equivalent of the third asymptotic in Proposition~\ref{estimateA_i}, and the proof proceeds similarly. It consists in adapting the proof of Lemma~\ref{bound_proba_exceed_level_kneps}. The only difference is that we need to control, uniformly over $m \ge 1$, the quantity:
    \begin{align*}
    \Pf_{\alpha}\bigg(\sum_{k=1}^{m} X_k \ge (\varepsilon^{\frac{1}{2}}n)^{\frac{\alpha}{\alpha-1}},\max_{1 \le k \le m}Y_k = \varepsilon n \bigg),
    \end{align*}
    where $(X_k,Y_k)_{k \ge 1}$ is a family of i.i.d. random variables distributed as $(\#\mathcal{C},\Height (\mathcal{C}))$ under $\Pf_{\alpha}$. Using the point (b) of Lemma~\ref{lemma_height}, we obtain 
    \begin{align*}
    \mathbb{E}_{\mathbf{T}}\bigg[F(\mathcal{C}) \mid  \Height (\mathcal{C})=n  \bigg] &=\mathbb{E}_{\mathbf{T}}\bigg[F({\mathcal{C}}_{H,n,\varepsilon}) \mid  \Height (\mathcal{C}) = n \bigg] + \mathcal{O}(\varepsilon) + \mathcal{O}(\varepsilon^{\frac{\alpha}{2(\alpha-1)}})\\
    & = \mathbb{E}_{\mathbf{T}}\bigg[F({\mathcal{C}}_{H,n,\varepsilon}) \mid  \Height (\mathcal{C}) = n \bigg] + \mathcal{O}(\varepsilon),
    \end{align*}
    where the last equality follows from the fact that ${\alpha}{2(\alpha-1)} \ge 1$.

\end{proof}

We can now conclude the proof of Theorem~\ref{theorem_exactconditioningheight}.
\begin{proof}
The proof is identical to that of Theorem~\ref{theorem_exactconditioning}. Using Proposition~\ref{asymptotics_height} and Fact \ref{fact:annealed conv}, we obtain that for $n$ large enough,
\begin{align*}
\bigg | \mathbb{E}_{\mathbf{T}}\bigg[F(\mathcal{C}) \mid \Height ( \mathcal{C}) = n \bigg] - \mathbb{E}_{\alpha}\bigg[F(\mathcal{C}) \mid \Height ( \mathcal{C} )= n \bigg]\bigg | = \mathcal{O}_{\varepsilon,\eta}(\delta )+\mathcal{O}_{\varepsilon}(\eta )+ \mathcal{O}(\varepsilon) + f(\eps),
\end{align*}
where $|f(\eps)| \downarrow 0$ as $\eps \downarrow 0$. We conclude by letting $\delta \to 0$, then $\eta \to 0$ and finally $\varepsilon \to 0$ and using the fact that the same convergence is known to hold in the annealed setting.
\end{proof}

\section{Convergence of the simple random walk}\label{sctn:SRW conv}
An immediate consequence of Theorem \ref{main_theorem} is the quenched convergence of the law of a simple random walk on $\C_{=n}$ to Brownian motion on the stable tree. This latter object can be defined rigorously using the theory of resistance forms; see \cite{croydon2018introduction} for an introduction. In particular, for any metric space equipped with a so-called \textit{resistance metric} and a \textit{measure}, the general theory allows us to associate a stochastic process with this metric and measure. \textit{Brownian motion on the stable tree} can therefore be defined as the stochastic process associated with the metric-measure space $(\mathcal{T}^{=1}_{\alpha},d_{\mathcal{\mathcal{T}_{\alpha}}},\nu_{\alpha},\rho_{\alpha})$. On trees (both on discrete trees and continuum trees called \textit{real trees} - see \cite[Definition 1.1]{le2006random}, for example), the graph distance between any two points is a resistance metric. In the case of a discrete graph (such as $\C_{=n}$), we will be interested in the process associated with the graph metric and the degree measure on the vertices. This is the continuous-time stochastic process with generator 
\begin{equation*}
(\mathcal{L}f)(x) = \frac{1}{\deg x} \sum_{y \sim x}(f(y)-f(x));
\end{equation*}
in other words a continuous-time random walk on $\C_{=n}$ that has an \textsf{exponential}($1$) waiting time at each vertex at each time step, and then moves to a uniformly chosen neighbour, and continues to evolve independently in this way. Due to concentration of the sums of these exponential waiting times, it is elementary to show that this stochastic process has the same scaling limit as a discrete time simple random walk on $\C_{=n}$. Moreover, letting $\deg$ denote the degree measure on vertices, it is also straightforward to verify that
\[
\dGHP{(\mathcal{C}_{=n}, \gamma n^{-\left(1-\frac{1}{\alpha}\right)}d_n, n^{-1}\nu_n,\rho_n) , (\mathcal{C}_{=n}, \gamma n^{-\left(1-\frac{1}{\alpha}\right)}d_n, \frac{1}{2}n^{-1} \deg ,\rho_n) } \leq  \gamma n^{-\left(1-\frac{1}{\alpha}\right)}.
\]
The main result of \cite{croydon2018scaling} asserts (under some mild conditions) that, if a sequence of (resistance) metric-measure spaces converges to a limit, then the associated stochastic processes also converge in law. This allows us to deduce that, for $\bPb$-almost every tree, the law of a simple random walk on $\C_{=n}$ converges under rescaling to the law of Brownian motion on $\Ta^{= 1}$. The result is slightly awkward to state rigorously. Applying the Skorokhod representation theorem leads to the formulation of Corollary \ref{cor:SRW convergence}, i.e. quenched convergence of the \textit{quenched} law of the random walk. An alternative approach is to use the framework developed in \cite{Khezeli}, using the extended topology defined in Section \ref{sctn:GHP topology}, which allows us to state the quenched convergence of the \textit{annealed} law, defined as follows:  for a stochastic process $X^K$ on a random state-space $K$, equipped with a metric $d_K$, measure $\mu_K$ and distinguished point $\rho_K$, we define the associated annealed law of $X^K$ started from $\rho_K$ to be the probability measure on $\widetilde{\mathbb{K}}_c$ given by
\[\mathbb{P}_K\left(\cdot\right):=\int P^K_{\rho_K}\left((K,d_K,\mu_K,\rho_K,X^K)\in \cdot\right)\mathbb{P}\left(d(K,d_K,\mu_K,\rho_K)\right),\]
where $\mathbb{P}$ is the probability measure under which $(K,d_K,\mu_K,\rho_K)$ is selected, and, for a particular realisation of $K$, $P^K_{\rho_K}$ is the law of $X^K$ started from $\rho_K$. To state the theorem, we recall the definition of the space $\widetilde{\mathbb{K}}_c$ given in Section \ref{sctn:GHP topology}.

\begin{corollary}
As $n \to \infty$, the annealed laws of 
\[
\left(\mathcal{C}_{=n},\gamma n^{-\left(1-\frac{1}{\alpha}\right)}d_n, n^{-1}\nu_n, \rho_n, \left(X^{(n)}_{\lfloor \gamma^{-1} n^{\frac{2\alpha-1}{\alpha}} t \rfloor}\right)_{t \geq 0} \right)
\]
converge weakly as probability measures on $\widetilde{\mathbb{K}}_c$ to the annealed law of
\[
(\mathcal{T}_{\alpha}^{= 1},d_{\mathcal{\mathcal{T}_{\alpha}}},\nu_{\alpha},\rho_{\alpha}, (B_t)_{t \geq 0})
\]
with respect to the extended GHP topology on $\widetilde{\mathbb{K}}_c$ defined in Section \ref{sctn:GHP topology}.
\end{corollary}

\appendix

\section{Technical proposition}

\begin{proposition}\label{technical_proposition}
    Let $X$ be a random variable with finite first moment. Let $(X_{m,n})_{m,n \ge 0}$ be a sequence of random variables such that for any $n \ge 0$, the family $(X_{m,n})_{m \ge 0}$ is i.i.d and $(X_{m,n})_{n \ge 0}$ converges in distribution to $X$. We also assume that the family $(X_{m,n})_{m,n \ge 0}$ has finite moments and $\mathbb{E}[|X_{m,n}|] \underset{n \to +\infty}{\longrightarrow} \mathbb{E}[|X|] < \infty$. Then for any sequence $k_n \underset{n \to +\infty}{\rightarrow} +\infty$ and any family of random variables $(T_n)_{n\ge 0}$ such that $\displaystyle k_n^{-1} T_n \underset{n \to +\infty}{\overset{(\mathbb{P})}{\rightarrow}}1$, we have 
    \begin{align*}
        k_n^{-1} \displaystyle\somme{m=1}{T_n}{X_{m,n}}\underset{n \to +\infty}{\overset{(\mathbb{P})}{\rightarrow}}  \mathbb{E}[X].
    \end{align*}
\end{proposition}
\begin{proof}
\noindent Fix $(k_n)_{n \ge 0}$ a sequence as in the proposition. Let us first treat the case where $T_n = k_n$. Using the Skorokhod representation theorem, we can construct a probability space $(\Omega,\mathcal{F},\mathbb{P})$ and a family of random variables $(\widetilde{X}_{1,n})_{m \ge 0}$ such that $\widetilde{X}_{1,n} \underset{n \to +\infty}{\overset{a.s}{\longrightarrow}} \widetilde{X}_1$ and such that $\widetilde{X}_{1,n} \overset{(d)}{=} X_{1,n}$ and $\widetilde{X}_{1} \overset{(d)}{=} X_1$. Then, up to choosing $\Omega$ bigger, we introduce $((\widetilde{X}_{m,n})_{n \ge 1},\widetilde{X}_{m})_{m \ge 1}$, a countable number of independent copies of $((\widetilde{X}_{1,n})_{n \ge 1},\widetilde{X}_{1})$. In particular, for any $n \ge 1$, we have $(\widetilde{X}_{m,n})_{m \ge 1} \overset{(d)}{=}(X_{m,n})_{m\ge 1}$, thus in the rest of the proof we assume that $(X_{m,n})_{m,n\ge 1}$ satisfy the same assumptions than $(\widetilde{X}_{m,n})_{m,n \ge 1} $. Indeed, we have $ k_n^{-1}\displaystyle \somme{m=1}{k_n}{X_{m,n}} \overset{(d)}{=} k_n^{-1}\displaystyle \somme{m=1}{k_n}{\widetilde{X}_{m,n}}$, thus proving the result for $(\widetilde{X}_{m,n})$ or $(X_{m,n})$ is the same.\\

For any $m \ge 1$, since $\widetilde{X}_{m,n}\overset{a.s}{\underset{ n \to +\infty}{\longrightarrow}} \widetilde{X}_m $ and $\mathbb{E}[|\widetilde{X}_{m,n}|] \underset{n \to +\infty}{\longrightarrow} \mathbb{E}[|\widetilde{X}_m|]$, by the Scheffé's Lemma $\widetilde{X}_{m,n}{\overset{L^1}{\longrightarrow}} \widetilde{X}_m $. Then, we can write 
\begin{align*}
    k_n^{-1}\displaystyle \somme{m=1}{k_n}{\widetilde{X}_{m,n}} = k_n^{-1}\displaystyle \somme{m=1}{k_n}{\widetilde{X}_{m}} + k_n^{-1}\displaystyle \somme{m=1}{k_n}{(\widetilde{X}_{m,n}-\widetilde{X}_{m})}.
\end{align*}
The law of large numbers gives that $ k_n^{-1}\displaystyle \somme{m=1}{k_n}{\widetilde{X}_{m}} $ converges almost surely to $\mathbb{E}[\widetilde{X}_1]$. Moreover,

\begin{align*}
\mathbb{E}\bigg[k_n^{-1} \bigg |\displaystyle \somme{m=1}{k_n}{\widetilde{X}_{m,n}-\widetilde{X}_{m}}\bigg|\bigg] \le k_n^{-1} \somme{m=1}{k_n}{\mathbb{E}[|\widetilde{X}_{m,n}-\widetilde{X}_{m}|]} = \mathbb{E}[|\widetilde{X}_{1,n}-\widetilde{X}_{1}|] \to 0.
\end{align*}
We deduce that
\begin{align*}
     k_n^{-1} \displaystyle \somme{m=1}{k_n}{\widetilde{X}_{m,n}}\underset{n \to +\infty}{\overset{(\mathbb{P})}{\rightarrow}}\mathbb{E}[\widetilde{X}_1] = \mathbb{E}[X].
\end{align*}
To prove the general case when $\displaystyle k_n^{-1}  T_n \underset{n \to +\infty}{\overset{(\mathbb{P})}{\rightarrow}}1$, we simply need to prove the following convergence:
\begin{align*}
k_n^{-1} \displaystyle \somme{m=1}{k_n}{X_{m,n}} - k_n^{-1}  \displaystyle \somme{m=1}{T_n}{X_{m,n}}\underset{n \to +\infty}{\overset{(\mathbb{P})}{\rightarrow}}0.
\end{align*}
Fix $\varepsilon > 0$. For any $\delta  > 0$, we have
\begin{align*}
    \Pf\left(\left|\displaystyle \somme{m=1}{k_n}{X_{m,n}} - \displaystyle \somme{m=1}{T_n}{X_{m,n}}\right| \ge k_n \varepsilon\right) &\le \Pf\left(\bigg|\frac{T_n}{k_n}-1\bigg|\ge \delta \right)\\
    &+\Pf\left(\displaystyle \frac{\somme{m=(1-\delta )k_n}{(1+\delta)k_n}{|X_{m,n}|}}{k_n} \ge \varepsilon\right)\\
    &\le \Pf\left(\bigg|\frac{T_n}{k_n}-1\bigg|\ge \delta \right) + \frac{2\delta M}{\varepsilon},
\end{align*}
where the second inequality follows from the Markov inequality and $M > 0$ is a bound for the first moment of the variables $(X_{m,n})$. By our assumption on $T_n$, it follows that
\begin{align*}   \limsup_{n}\Pf\left(\left|\displaystyle \somme{m=1}{k_n}{X_{m,n}} - \displaystyle \somme{m=1}{T_n}{X_{m,n}}\right| \ge k_n \varepsilon\right) \le \frac{2\delta M}{\varepsilon}.
\end{align*}
The left-hand side does not depend on $\delta $ so letting $\delta \to 0$ we deduce the desired result. This concludes the proof.
\end{proof}

\section{Proof of Theorem \ref{main_theorem with generation size}\eqref{eqn:joint conv size}}\label{app:Ai}

The purpose of this section is to establish \eqref{eqn:joint conv size}, following the strategy outlined at the end of Section \ref{sctn:the events Ai}. in particular, we first verify in Lemma \ref{lem:kneps convergence} that $k_{n, \eps}$ converges to $k_{\eps}$ under rescaling. We then apply this in Proposition \ref{prop:joint conv uniform height} to verify that quantities in \eqref{eqn:joint conv size} converge when we instead condition on the height and moreover when we consider $k_{n, U}$ for a uniform random variable $U\sim$ \textsf{Uniform}($[0,1]$). We then transfer this back to $k_{n, c'}$ for almost every $c' \in [0,1]$ using Fubini's theorem (Corollary \ref{cor:unif to deter}), then switch the conditioning event from the height to the total volume in Proposition \ref{prop:Ai almost every}, which allows us to then restrict to the event $\{\#{\C}_{\geq n/2} \geq (1-c')n\}$ for $c' \in [0,1]$ (Corollary \ref{cor:size 1-c'}), which immediately implies \eqref{eqn:joint conv size}.

Note that one may hope to avoid such a lengthy argument by directly arguing that $\ell$ must be continuous at $k_{c'}$, and that an analogous property holds in the discrete setting. The subtlety is that the random variable $k_{c'}$ is highly dependent on the process $(\ell_t)_{t \geq 0}$ and so the result of Fact \ref{fact:csbp}(iii) cannot be applied in this way.

We start with a useful lemma.

\begin{lemma}\label{lem:kneps convergence}
\begin{enumerate}
\item Fix $\eps>0$, and suppose that
\begin{align*}
(\mathcal{C}_{\geq n(1-\eps)}, \gamma n^{-\left(1-\frac{1}{\alpha}\right)}d_{n}, {n}^{-1}\nu_{n},\rho_{n}) &\underset{n \to +\infty}{\longrightarrow} (\mathcal{T}^{\geq 1-\eps}_{\alpha},d_{\mathcal{\mathcal{T}_{\alpha}}},\nu_{\alpha},\rho_{\alpha})
\end{align*}
holds almost surely on the space $(\Omega_{\bT}, \F_{\bT}, \mathbb{P}_{\bT})$ (recall that this is under the conditioning $\# \C \geq (1-\eps) n$). Then we also have that
\[
\gamma n^{-\left(1-\frac{1}{\alpha}\right)}k_{n, \eps} \to k_{\eps}%(\Ta^{\geq 1-\eps})
\]
almost surely on $(\Omega_{\bT}, \F_{\bT}, \mathbb{P}_{\bT})$.
\item The same is true under conditioning the height to be at least  $n_c := c n^{1-1/\alpha}$, for any fixed $c>0$: in this case, if
\begin{align*}
({\C}_{H \geq n_c}, \gamma n^{-\left(1-\frac{1}{\alpha}\right)}d_n, n^{-1}\nu_n,\rho_n) &\underset{n \to +\infty}{\longrightarrow} (\mathcal{T}^{H \geq c\gamma}_{\alpha},d_{\mathcal{\mathcal{T}_{\alpha}}},\nu_{\alpha},\rho_{\alpha})
\end{align*}
almost surely, then also
\[
\gamma n^{-\left(1-\frac{1}{\alpha}\right)}k_{n, \eps} \mathbbm{1}\{k_{n, \eps} < \infty\} \to k_{\eps}%(\Ta^{H \geq c\gamma^{-1}})
\mathbbm{1}\{k_{\eps} < \infty\}
\]
almost surely.
\end{enumerate}
\end{lemma}
\begin{proof}
We prove the first point. By the Skorokhod representation theorem and standard results about GHP embeddings (see \cite[Theorem 4.11]{evans2006probability} for further details) we can a.s. find a sequence $\delta_n \downarrow 0$ and isometrically embed $(\C_{\geq n(1-\eps)})_{n \geq 1}$ and $\Ta^{\geq 1-\eps}$ into a common metric space $(M,D_M)$ so that
\[
d_H\left(\left(\mathcal{C}_{\geq n(1-\eps)}, \gamma n^{-\left(1-\frac{1}{\alpha}\right)}d_n\right), (\Ta^{\geq 1-\eps}, d_{\Ta})\right) \vee d_P(n^{-1}\nu_n, \nu) \vee D_M(\rho_n, \rho) \leq \delta_n
\]
Assume this is the case and now suppose that $K_2 > K_1 > k_{\eps}$. By the triangle inequality, we have that $B(\rho, K_1)^{\delta_n} \cap \C_{\geq n(1-\eps)} \subset B(\rho_n, K_2) \cap \C_{\geq n(1-\eps)}$ (balls are measured with respect to $D_M$) for all sufficiently large $n$. It follows that
\[
1 \leq \nu (B(\rho, K_1)) \leq n^{-1}\nu_n (B(\rho, K_1)^{\delta_n}) + \delta_n \leq n^{-1}\nu_n ( B(\rho_n, K_2) ) + \delta_n.
\]
Taking $\delta_n \downarrow 0$ we deduce that $\limsup_{n \to \infty}\gamma n^{-\left(1-\frac{1}{\alpha}\right)} k_{n, \eps} \leq K_2$ and thus also that
\[
\limsup_{n \to \infty} \gamma n^{-\left(1-\frac{1}{\alpha}\right)} k_{n, \eps} \leq k_{\eps}
\]
(by taking $K_2 \downarrow k_{\eps}$). A similar argument shows that  $\liminf_{n \to \infty} \gamma n^{-\left(1-\frac{1}{\alpha}\right)} k_{n, \eps} \geq k_{\eps}$.

The result when conditioning on the height follows by exactly the same arguments.
\end{proof}

This is useful to prove the following proposition.

\begin{proposition}\label{prop:joint conv uniform height}
Fix $c \in (0, 1)$. Let $U \sim \textsf{Uniform} (0,1)$, independently of everything else. 
\begin{align*}
&\left(({\C}_{H \geq n_c}, \gamma n^{-(1-\frac{1}{\alpha})}d_{n_c}, n^{-1} \nu_{n_c},\rho_{n_c}) , \frac{Y_{ k_{n,U}}\mathbbm{1}\{n_c \leq k_{n,U} < \infty\}}{\gamma n^{\frac{1}{\alpha}}} \right) \\
&\qquad \underset{n \to +\infty}{\overset{(d)}{\longrightarrow}} \left((\mathcal{T}^{H\geq c \gamma}_{\alpha},d_{\mathcal{\mathcal{T}_{\alpha}}},\nu_{\alpha},\rho_{\alpha}), \ell_{k_{U}}\mathbbm{1}\{c\gamma \leq k_{U} < \infty\} \right).
\end{align*}
\end{proposition}
\begin{proof}
We first note the trivial extension of \eqref{eqn:joint conv height}: take $T > 1$ and sample $U \sim \textsf{Uniform}([1,T])$, independently of all other random variables. Then \eqref{eqn:joint conv height} holds with $t=U$.

We now work on the space $[0, 1] \times \Omega_{\bT}$, and we note that, on this space, $\mathbbm{1}\{k_{n,{U}} \geq n_c\}$ converges to $\mathbbm{1}\{k_{U} \geq c\gamma \}$, jointly with the metric-measure space convergence in the above statement, as a trivial consequence of Lemma \ref{lem:kneps convergence}. Hence it suffices to work on the event where these indicators are both $1$.

We will show that, when we restrict to some high probability event, the law of $k_{n,U}$ is absolutely continuous with respect to that of $U$.

Pick $T< \infty$ such that
\[
\Ptcond{\Height ({\C}_{H \geq n_c}) \geq T n^{1-\frac{1}{\alpha}}}{{k_{n,U} \geq n_c}}{\bT} < \delta
\]
for all $n \geq 1$ (this is possible by Theorem \ref{main_theorem}; note that the conditioning event has uniformly positive probability). We moreover consider the following good event:
\[
E_{K,c,T} := \left\{ K^{-1} \leq n^{-1/\alpha}\inf_{m \in [c n^{1-\frac{1}{\alpha}}, T n^{1-\frac{1}{\alpha}}]} Y_m \leq n^{-1/\alpha} \sup Y_m \leq K\right\} \cap \{\Height ({\C}_{H \geq n_c}) < T n^{1-\frac{1}{\alpha}}\}.
\]
It follows from Lemma \ref{lem:AV CSBP convergence} that we can also pick $K< \infty$ such that $\Ptcond{E_{K,c,T}}{k_{n,U} \geq n_c}{} \geq 1-2\delta$ for all sufficiently large $n$.

Then, thanks to our initial observation, and the Skorokhod representation theorem (the product of two Polish spaces is Polish), if we take $U \sim \textsf{Uniform}([1,T])$, independently of all other random variables, we can assume that
\begin{align}\label{eqn:conv as uniform gen}
\left(({\C}_{H \geq n_c}, \gamma n^{-(1-\frac{1}{\alpha})}d_{n_c}, n^{-1} \nu_{n_c},\rho_{n_c}) , \frac{Y_{n_cU}}{\gamma n^{\frac{1}{\alpha}}} \right) &\underset{n \to +\infty}{{\longrightarrow}} \left((\mathcal{T}^{H\geq c \gamma}_{\alpha},d_{\mathcal{\mathcal{T}_{\alpha}}},\nu_{\alpha},\rho_{\alpha}), \ell_{c\gamma U} \right)
\end{align}
almost surely on the space $(\Omega_{\bT}, \F_{\bT}, \mathbb{P}_{\bT})$.
%Throughout the rest of the proof we condition on the event $\{k_{n,{1}} \geq n_c\}$. Note that this event has uniformly positive probability. 
In particular, by our assumption of almost sure convergence, we can find $N<\infty$ such that 
\begin{equation}\label{eqn:YU lU close}
\Ptcond{\left|\frac{Y_{\lfloor n_{c}U \rfloor}}{\gamma n^{\frac{1}{\alpha}}}- \ell_{c\gamma U} \right| > \delta}{k_{n,U} \geq n_c}{\bT}< \frac{1}{K^3T}\delta
\end{equation}
for all $n \geq N$.

Moreover, on the event $E_{K,c,T}$, we have that 
\begin{equation}\label{eqn:abs cont U k}
\Ptcond{k_{n,{U}}=m}{E_{K,c,T}, k_{n,U} \geq n_c}{\bT} \leq {K^3T} \Pt{\lfloor n_cU \rfloor=m}.
\end{equation}
For $m \geq 1$ set $m^* = m + \textsf{Uniform}[0,1]$, where the latter random variable is completely independent of everything else (this will be more convenient for the coupling with $U$). Combining \eqref{eqn:YU lU close} and \eqref{eqn:abs cont U k} gives (for all sufficiently large $n$):
\begin{equation*}
\Ptcond{\left|\frac{Y_{k_{n,{U}}}}{\gamma n^{\frac{1}{\alpha}}}- \ell_{\gamma n^{-(1-\frac{1}{\alpha})}k_{n,{U}}^*} \right| > \delta}{E_{K,c,T}, k_{n,U} \geq n_c}{\bT}< \delta.
\end{equation*}
It follows from Lemma \ref{lem:kneps convergence} that $\gamma n^{-(1-\frac{1}{\alpha})}k_{n,{U}}^* \to k_{U}$ almost surely, jointly with \eqref{eqn:conv as uniform gen} (recall that we are working on the probability space $[0,1] \times \Omega_{\bT}$, so once we've sampled $U$ on $[0, 1]$, we can work pointwise on the space $[0, 1]$, so we can assume that $U$ is fixed to be constant almost everywhere on $(\Omega_{\bT}, \F_{\bT}, \mathbb{P}_{\bT})$ and apply the result of Lemma \ref{lem:kneps convergence}). We would further like that $\ell_{\gamma n^{-(1-\frac{1}{\alpha})}k_{n,{U}}^*} \to \ell_{k_{U}}$. This follows by a similar argument: it is known (by Fact \ref{fact:csbp}(iii)) that $\ell$ is almost surely continuous at $c\gamma U$. By taking the limit in \eqref{eqn:abs cont U k} we deduce the same for $k_{U}$ and thus it follows that
\begin{equation*}
\Ptcond{\left|\frac{Y_{k_{n,{U}}}}{\gamma n^{\frac{1}{\alpha}}}- \ell_{k_{U}} \right| > 2\delta}{E_{K,c,T}, k_{n,U} \geq n_c}{\bT}< 2\delta.
\end{equation*}
for all sufficiently large $n$. Since $E_{K, c, T}$ occurred with probability at least $1-2\delta$ we can then remove this event from the conditioning provided we increase the right hand side to $4 \delta$. Since $\delta>0$ was arbitrary this proves that the desired convergence holds almost surely on the space $(\Omega_{\bT}, \F_{\bT}, \mathbb{P}_{\bT})$ and on the event $\{k_{n,U} \geq n_c\}$, and the result follows.
\end{proof}

An application of Fubini's theorem gives the following.

\begin{corollary}\label{cor:unif to deter}
Fix $c \in (0, 1)$. For almost every $c' \in [0,1]$, the following holds:
\begin{align*}
&\left(({\C}_{H \geq n_c}, \gamma n^{-(1-\frac{1}{\alpha})}d_{n_c}, n^{-1} \nu_{n_c},\rho_{n_c}) , \frac{Y_{ k_{n,{c'}}}\mathbbm{1}\{n_c \leq k_{n,{c'}} < \infty\}}{\gamma n^{\frac{1}{\alpha}}} \right) \\
&\qquad \underset{n \to +\infty}{\overset{(d)}{\longrightarrow}} \left((\mathcal{T}^{H\geq c \gamma}_{\alpha},d_{\mathcal{\mathcal{T}_{\alpha}}},\nu_{\alpha},\rho_{\alpha}), \ell_{k_{{c'}}}\mathbbm{1}\{c\gamma \leq k_{{c'}} < \infty\} \right).
\end{align*}
\end{corollary}

In fact this extends to all $c' \in [0,1]$ using the fact that the the limiting process $\ell$ has no
fixed discontinuities and that the mapping $c' \mapsto k_{c'}$ is continuous. But the above statement is sufficient for our needs.

\begin{proposition}\label{prop:Ai almost every}
Fix $c \in (0, 1)$. For almost every $c' \in [0,1]$, the following holds:
\begin{align*}
\left(({\C}_{\geq n/2}, \gamma n^{-(1-\frac{1}{\alpha})}d_{n}, n^{-1} \nu_{n},\rho_{n}) , \frac{Y_{ k_{n,{c'}}}\mathbbm{1}\{k_{n,{c'}} < \infty\}}{\gamma n^{\frac{1}{\alpha}}} \right) &\underset{n \to +\infty}{\overset{(d)}{\longrightarrow}} \left((\mathcal{T}^{\geq 1/2}_{\alpha},d_{\mathcal{\mathcal{T}_{\alpha}}},\nu_{\alpha},\rho_{\alpha}), \ell_{k_{{c'}}}\mathbbm{1}\{k_{{c'}} < \infty\} \right).
\end{align*}
\end{proposition}
\begin{proof}
Fix $\bT$ and $c', \delta>0$ and assume that the statement of Corollary \ref{cor:unif to deter} holds for $c'$. By Theorem \ref{main_theorem}\eqref{eqn:first geq conv size}, we can find $c>0$ such that $\limsup_{n \to \infty} \Pt{\Height ({\C}_{\geq n/2}) < c n^{1-\frac{1}{\alpha}} \text{ or } k_{n, c'} \leq n_c} \leq \delta$ and $\eta := \lim_{n \to \infty} \Ptcond{ \#\C \geq n/2}{\Height (\C) \geq c n^{1-\frac{1}{\alpha}}}{\bT} > 0$ (this limit exists by Theorem \ref{main_theorem}\eqref{eqn:first geq conv height}).

For any closed set $A \in \F_{\bT}$, we thus have 
\begin{align*}
|\Ptcond{\C \in A}{\# \C \geq n/2}{\bT} - \Ptcond{\C \in A, k_{n, c'} > n_c}{\# \C \geq n/2, \Height (\C) \geq c n^{1-\frac{1}{\alpha}}}{\bT} | \leq 2\delta.
\end{align*}
%and moreover 
%\begin{align*}
%\Ptcond{\C \in A}{\# \C \geq n, \Height (\C) \geq c n^{1-\frac{1}{\alpha}}}{\bT} = \frac{\Ptcond{\C \in A, \# \C \geq n}{ \Height (\C) \geq c n^{1-\frac{1}{\alpha}}}{\bT}}{\Ptcond{\# \C \geq n}{ \Height (\C) \geq c n^{1-\frac{1}{\alpha}}}{\bT}}
%\end{align*}
%We know that the denominator converges to $\eta = \Ptcond{\nu (\Ta) \geq 1}{\Height (\Ta) \geq c}{\bT}$.
Let a superscript of $(n)$ denote the following rescaled version of $\C$, and similarly in the continuum:
\begin{align*}
\C^{(n)} &= \left((\mathcal{C}, \gamma n^{-\left(1-\frac{1}{\alpha}\right)}d, n^{-1}\nu,\rho_{\C}), \frac{Y_{\lfloor k_{n, c'}\rfloor}\mathbbm{1}\{n_c \leq k_{n,{c'}} < \infty\}}{\gamma n^{1/\alpha}} \right), \\
\mathcal{T}^{H\geq c \gamma}_{\alpha} &= \left((\mathcal{T}^{H\geq c \gamma}_{\alpha},d_{\mathcal{\mathcal{T}_{\alpha}}},\nu_{\alpha},\rho_{\alpha}), \ell_{k_{{c'}}}\mathbbm{1}\{c\gamma \leq k_{{c'}} < \infty\} \right),
\end{align*}
Note that we have added an extra lower bound in the indicators in both cases, compared with the statement of the proposition. As above, we have for any closed $A$ that 
\begin{align*}
&\qquad \Ptcond{\C^{(n)} \in A, k_{n, c'} > n_c}{\# \C \geq n/2, \Height (\C) \geq c n^{1-\frac{1}{\alpha}}}{\bT} \\
&= \frac{\Ptcond{\C^{(n)} \in A, k_{n, c'} > n_c, n^{-1} \#\C \geq 1/2}{ \Height (\C) \geq c n^{1-\frac{1}{\alpha}}}{\bT}}{\Ptcond{n^{-1} \#\C \geq 1/2}{ \Height (\C) \geq c n^{1-\frac{1}{\alpha}}}{\bT}}.% \\
%&= \frac{\Pt{\widetilde{C}_{n_{c}}^{(n_{c}^{\frac{\alpha}{\alpha-1}})} \in A^{\gamma, c}, (\gamma n_{c})^{-\frac{\alpha}{\alpha-1}} \nu_{n_{c}}(\widetilde{C}_{n_{c}}) \geq (\gamma c)^{-\frac{\alpha}{\alpha-1}}}}{\Pt{(\gamma n_{c})^{-\frac{\alpha}{\alpha-1}} \nu_{n_{c}}(\widetilde{C}_{n_{c}}) \geq (\gamma c)^{-\frac{\alpha}{\alpha-1}}}}
\end{align*}
By Corollary \ref{cor:unif to deter}, we can take limits of both the numerator and the denominator to deduce that, as $n \to \infty$, this converges to
\begin{align*}
\frac{\Pt{\Ta^{H \geq \gamma c} \in A, k_{c'} > c\gamma, \nu(\Ta^{H \geq \gamma c}) \geq 1/2}}{\Pt{\nu(\Ta^{H \geq \gamma c}) \geq 1/2}} &= \Ptcond{\Ta^{H \geq \gamma c} \in A, k_{c'} > c\gamma}{\nu(\Ta^{H \geq \gamma c}) \geq 1/2}{\bT} \\
&= \Ptcond{\Ta \in A, k_{c'} > c\gamma}{\nu(\Ta) \geq 1/2, {\Height (\Ta) \geq \gamma c}}{\bT} \\
&= \frac{\Ptcond{\Ta \in A,  k_{c'} > c\gamma, {\Height (\Ta) \geq \gamma c}}{\nu(\Ta) \geq 1/2}{\bT}}{\Ptcond{{\Height (\Ta) \geq \gamma c}}{\nu(\Ta) \geq 1/2}{\bT}}.
\end{align*}
Standard properties of the stable trees (in particular that $\Height (\Ta)$ and $k_{c'}$ are almost surely non-zero under the conditioning $\nu (\Ta) \geq 1/2$ - e.g. see \cite[Theorem 1.8 and Equation (72)]{duquesnewangdecomp} for the result for the height) imply that this expression converges to $\Ptcond{\Ta \in A}{\nu(\Ta) \geq 1/2}{\bT}$ as $c \downarrow 0$. The result then follows by the Portmanteau theorem since we can take $\delta>0$ and $c>0$ arbitrarily small.
\end{proof}

To deduce the desired result, one then just has to note that, for fixed $c' \in [0, \frac{1}{2})$, the cluster ${\C}_{\geq (1-c')n}$ has the law of ${\C}_{\geq n/2}$ but under an additional conditioning on $\{ \# {\C}_{\geq n/2} \geq (1-c')n\}$, which has uniformly positive probability as $n \to \infty$, and which converges to the event $\{\nu (\Ta^{\geq 1/2}) \geq 1-c'\}$ jointly with the convergence of Corollary \ref{cor:unif to deter}.

\begin{corollary}\label{cor:size 1-c'}
For almost every $c' \in [0,1/2]$, the following holds:
\begin{align*}
\left(({\C}_{\geq (1-c')n}, \gamma n^{-(1-\frac{1}{\alpha})}d_{n}, n^{-1} \nu_{n},\rho_{n}) , \frac{Y_{ k_{n,{c'}}}\mathbbm{1}\{k_{n,{c'}} < \infty\}}{\gamma n^{\frac{1}{\alpha}}} \right) &\underset{n \to +\infty}{\overset{(d)}{\longrightarrow}} \left((\mathcal{T}^{\geq 1-c'}_{\alpha},d_{\mathcal{\mathcal{T}_{\alpha}}},\nu_{\alpha},\rho_{\alpha}), \ell_{k_{{c'}}}\mathbbm{1}\{k_{{c'}} < \infty\} \right).
\end{align*}
\end{corollary}

Finally, one can remove the indicators above since both are one almost surely under the above conditioning, in order to deduce \eqref{eqn:joint conv size}.

\printbibliography

\end{document}